# BASIC NEUTROSOPHIC ALGEBRAIC STRUCTURES AND THEIR APPLICATION TO FUZZY AND NEUTROSOPHIC MODELS

**Dr. W. B. Vasantha Kandasamy**
**Florentin Smarandache**

# Basic Neutrosophic Algebraic Structures and Their Application to Fuzzy and Neutrosophic Models

**W. B. Vasantha Kandasamy**
Department of Mathematics
Indian Institute of Technology, Madras
Chennai – 600036, India
e-mail: **vasantha@iitm.ac.in**
web: **http://mat.iitm.ac.in/~wbv**

**Florentin Smarandache**
Department of Mathematics
University of New Mexico
Gallup, NM 87301, USA
e-mail: **smarand@gallup.unm.edu**

## CONTENTS

## Preface

Study of neutrosophic algebraic structures is very recent. The introduction of neutrosophic theory has put forth a significant concept by giving representation to indeterminates. Uncertainty or indeterminacy happen to be one of the major factors in almost all real-world problems. When uncertainty is modeled we use fuzzy theory and when indeterminacy is involved we use neutrosophic theory. Most of the fuzzy models which deal with the analysis and study of unsupervised data make use of the directed graphs or bipartite graphs. Thus the use of graphs has become inevitable in fuzzy models. The neutrosophic models are fuzzy models that permit the factor of indeterminacy. It also plays a significant role, and utilizes the concept of neutrosophic graphs. Thus neutrosophic graphs and neutrosophic bipartite graphs plays the role of representing the neutrosophic models. Thus to construct the neutrosophic graphs one needs some of the neutrosophic algebraic structures viz. neutrosophic fields, neutrosophic vector spaces and neutrosophic matrices. So we for the first time introduce and study these concepts. As our analysis in this book is application of neutrosophic algebraic structure we found it deem fit to first introduce and study neutrosophic graphs and their applications to neutrosophic models.

This book is organized into four chapters. In Chapter One we introduce some of the basic neutrosophic algebraic structures essential for the further development of the other chapters. Chapter Two recalls basic graph theory definitions and results which has interested us and for which we give the neutrosophic analogues. In this chapter we give the application of graphs in fuzzy models. An entire section is devoted for this purpose. Chapter Three introduces many new neutrosophic concepts in graphs and applies it to the case of neutrosophic cognitive maps and neutrosophic relational maps. The last section of this chapter clearly illustrates how the neutrosophic graphs are utilized in the neutrosophic models. The final chapter gives some problems about neutrosophic graphs which will make one understand this new subject.

The authors are grateful to Ilanthenral for formatting the book and drawing all the graphs and diagrams in this book. Our thanks are also due to Meena Kandasamy and Dr. Kandasamy for helping us in all possible ways in the preparation of this book.

Finally, we dedicate this book to the social revolutionary Periyar (lit. Great One) E V Ramasamy (1879-1973) on his 125th birth anniversary. He valiantly fought against Brahmin supremacy for the emancipation of the oppressed peoples of the downtrodden castes. This present book was a side-project while we were both working towards a fuzzy and neutrosophic understanding of Periyar's ideology and approach.

W.B.VASANTHA KANDASAMY
FLORENTIN SMARANDACHE

## Chapter One

# Introduction to Some Neutrosophic Algebraic Structure

In this chapter, we define some new neutrosophic algebraic structures like neutrosophic fields, neutrosophic spaces and neutrosophic matrices and illustrate them with examples. For these notions are used in the definition of neutrosophic graphs and its applications to neutrosophic cognitive maps which is dealt in the later chapters of this book.

Throughout this book by '*I* ' we denote the indeterminacy of any notion / concept / relation. That is when we are not in a position to associate a relation between any two concepts then we denote it as an indeterminacy or when a concept cannot be defined we use the symbol *I*.

This chapter has two sections. In section one we define neutrosophic fields and in section two we define neutrosophic vector spaces and illustrate them by examples.

## 1.1 Neutrosophic fields

In this book we assume all fields to be real fields of characteristic 0 all vector spaces are taken as real spaces over reals and we denote the indeterminacy by '*I* ' as i will make a confusion as it denotes the imaginary value, viz $i^2 = -1$ that is $\sqrt{-1} = i$. The indeterminacy $I$ is such that $I \cdot I = I^2 = I$.

**DEFINITION 1.1.1:** *Let K be the field of reals. We call the field generated by* $K \cup I$ *to be the neutrosophic field for it involves the indeterminacy factor in it. We define* $I^2 = I$, $I + I = 2I$ *i.e.,* $I + \ldots + I = nI$, *and if* $k \in K$ *then* $k.I = kI$, $0I = 0$. *We denote the* *neutrosophic field by K(I) which is generated by K $\cup$ I that is K(I) = $\langle K \cup I \rangle$. ($\langle K \cup I \rangle$ denotes the field generated by K and I.*

***Example 1.1.1:*** Let R be the field of reals. The neutrosophic field of reals is generated by $\langle R \cup I \rangle$ i.e. R($I$) clearly R $\subset$ $\langle R \cup I \rangle$.

***Example 1.1.2:*** Let Q be the field of rationals. The neutrosophic field of rational is generated by Q $\cup$ $I$ denoted by Q($I$).

**DEFINITION 1.1.2:** *Let K(I) be a neutrosophic field we say K(I) is a prime neutrosophic field if K(I) has no proper subfield which is a neutrosophic field.*

***Example 1.1.3:*** Q($I$) is a prime neutrosophic field where as R($I$) is not a prime neutrosophic field for Q($I$) $\subset$ R ($I$).

It is very important to note that all neutrosophic fields used in this book are of characteristic zero. Likewise we can define neutrosophic subfield.

**DEFINITION 1.1.3:** *Let K(I) be a neutrosophic field, P $\subset$ K(I) is a neutrosophic subfield of P if P itself is a neutrosophic field. K(I) will also be called as the extension neutrosophic field of the neutrosophic field P.*

## 1.2 Neutrosophic Vector spaces

Now we proceed on to define neutrosophic vector spaces, which can be defined over fields or neutrosophic fields. We can define two types of neutrosophic vector spaces one when it is a neutrosophic vector space over ordinary field other being neutrosophic vector space over neutrosophic fields. To this end we have to define neutrosophic group under addition.

**DEFINITION 1.2.1:** *We know Z is the abelian group under addition. Z(I) denote the additive abelian group generated by the set Z and I, Z(I) is called the neutrosophic abelian group under '+'.*

*Thus to define basically a neutrosophic group under addition we need a group under addition. So we proceed on to define neutrosophic abelian group under addition.*

Suppose G is an additive abelian group under '+'. $G(I) = \langle G \cup I \rangle$, additive group generated by G and $I$, $G(I)$ is called the neutrosophic abelian group under '+'.

***Example 1.2.1:*** Let Q be the group under '+'; $Q(I) = \langle Q \cup I \rangle$ is the neutrosophic abelian group under addition; '+'.

***Example 1.2.2:*** R be the additive group of reals, $R(I) = \langle R \cup I \rangle$ is the neutrosophic group under addition.

***Example 1.2.3:*** $M_{n\times m}(I) = \{(a_{ij}) \mid a_{ij} \in Z(I)\}$ be the collection of all $n \times m$ matrices under '+'; $M_{n\times m}(I)$ is a neutrosophic group under '+'.

Now we proceed on to define neutrosophic subgroup.

**DEFINITION 1.2.2:** *Let G(I) be the neutrosophic group under addition. $P \subset G(I)$ be a proper subset of G(I). P is said to be the neutrosophic subgroup of G(I) if P itself is a neutrosophic group i.e. $P = \langle P_1 \cup I \rangle$ where $P_1$ is an additive subgroup of G.*

***Example 1.2.4:*** Let $Z(I) = \langle Z \cup I \rangle$ be a neutrosophic group under '+'. $\langle 2Z \cup I \rangle = 2Z(I)$ is the neutrosophic subgroup of $Z(I)$.

In fact $Z(I)$ has several neutrosophic subgroups.

Now we proceed on to define the notion of neutrosophic quotient group.

**DEFINITION 1.2.3:** *Let $G(I) = \langle G \cup I \rangle$ be a neutrosophic group under '+', suppose P(I) be a neutrosophic subgroup of G (I) then the neutrosophic quotient group*

$$\frac{G(I)}{P(I)} = \{a + P(I) \mid a \in G(I)\}.$$

***Example 1.2.5:*** Let $Z(I)$ be a neutrosophic group under addition, Z the group of integers under addition, $P = 2Z(I)$ is a neutrosophic subgroup of $Z(I)$, the neutrosophic quotient group

$$\frac{Z(I)}{P} = \{a + 2Z(I) \mid a \in Z(I)\} = \{(2n+1) + (2n+1)I \mid n \in Z\}.$$

Clearly $\frac{Z(I)}{P}$ is a group. For P = 2Z ($I$ ) serves as the additive identity. Take a, b $\in \frac{Z(I)}{P}$. If a,b $\in$ Z($I$ ) \ P then two possibilities occur.

a + b is odd times $I$ or a + b is odd or a + b is even times $I$ or even if a + b is even or even times $I$ then a + b $\in$ P, if a + b is odd or odd times $I$, a + b $\in \frac{Z(I)}{P = 2Z(I)}$.

It is easily verified that P acts as the identity and every element in $\frac{Z(I)}{P}$, that is for every

$$a + 2Z\,(I) \in \frac{Z(I)}{2Z(I)}$$

*has a unique inverse.*

Now we proceed on to define the notion of neutrosophic vector spaces over fields and then we define neutrosophic vector spaces over neutrosophic fields.

**DEFINITION 1.2.4:** *Let G(I) by an additive abelian neutrosophic group. K any field. If G(I) is a vector space over K then we call G(I) a neutrosophic vector space over K.*

Now we give the notion of strong neutrosophic vector space.

**DEFINITION 1.2.5:** *Let G(I) be a neutrosophic abelian group. K(I) be a neutrosophic field. If G(I) is a vector space over K(I) then we call G(I) the strong neutrosophic vector space.*

**THEOREM 1.2.1:** *All strong neutrosophic vector space over K(I) are a neutrosophic vector space over K; as K $\subset$ K(I).*

*Proof:* Follows directly by the very definitions.

Thus when we speak of neutrosophic spaces we mean either a neutrosophic vector space over K or a strong neutrosophic vector space over the neutrosophic field K(*I*). By basis we mean a linearly independent set which spans the neutrosophic space.

Now we illustrate with an example.

***Example 1.2.6:*** Let $R(I) \times R(I) = V$ be an additive abelian neutrosophic group over the neutrosophic field *R(I).* Clearly V is a strong neutrosophic vector space over *R(I).* The basis of V are {(0,1), (1,0)}.

***Example 1.2.7:*** Let $V = R(I) \times R(I)$ be a neutrosophic abelian group under addition. V is a neutrosophic vector space over R. The neutrosophic basis of V are {(1,0), (0,1), (*I*,0), (0,*I*)}, which is a basis of the neutrosophic vector space V over R.

A study of these basis and its relations happens to be an interesting form of research.

**DEFINITION 1.2.6:** *Let G(I) be a neutrosophic vector space over the field K. The number of elements in the neutrosophic basis is called the neutrosophic dimension of G(I).*

**DEFINITION 1.2.7:** *Let G(I) be a strong neutrosophic vector space over the neutrosophic field K(I). The number of elements in the strong neutrosophic basis is called the strong neutrosophic dimension of G(I).*

We denote the neutrosophic dimension of G(*I*) over K by $N_k$ (dim) of G (*I*) and that the strong neutrosophic dimension of G (*I*) by $SN_{K(I)}$ (dim) of G(*I*).

Now we define the notion of neutrosophic matrices.

**DEFINITION 1.2.8:** *Let* $M_{nxm} = \{(a_{ij}) \,/\, a_{ij} \in K(I)\}$*, where K (I), is a neutrosophic field.We call* $M_{nxm}$ *to be the neutrosophic matrix.*

***Example 1.2.8:*** Let Q(*I* ) = ⟨Q ∪ *I*⟩ be the neutrosophic field.

$$M_{4x3} = \begin{pmatrix} 0 & 1 & I \\ -2 & 4I & 0 \\ 1 & -I & 2 \\ 3I & 1 & 0 \end{pmatrix}$$

is the neutrosophic matrix, with entries from rationals and the indeterminacy *I*.

We define product of two neutrosophic matrices and the product is defined as follows:

Let

$$A = \begin{pmatrix} -1 & 2 & -I \\ 3 & I & 0 \end{pmatrix}_{2\times 3} \quad \text{and} \quad B = \begin{pmatrix} I & 1 & 2 & 4 \\ 1 & I & 0 & 2 \\ 5 & -2 & 3I & -I \end{pmatrix}_{3x4}$$

$$AB = \begin{bmatrix} -6I+2 & -1+4I & -2-3I & I \\ -4I & 3+I & 6 & 12+2I \end{bmatrix}_{2x4}$$

(we use the fact $I^2 = I$).

Let $M_{n\times n} = \{(a_{ij}) \mid (a_{ij}) \in Q(I)\}$, $M_{n\times n}$ is a neutrosophic vector space over Q and a strong neutrosophic vector space over $Q(I)$.

To define Neutrosophic graphs and Neutrosophic Cognitive Maps we direly need the notion of Neutrosophic Matrices. We use square neutrosophic matrices for Neutrosophic Cognitive Maps and use rectangular neutrosophic matrices for Neutrosophic Relational Maps.

Chapter Two

# SOME BASIC RESULTS ON GRAPH THEORY AND THEIR APPLICATIONS TO FUZZY MODELS

In this chapter we recall some of the basic results from Graph Theory. This is mainly done to make this book a self contained one. We just give only the definition and results on graph theory which we have used. We have taken the results from [2-6, 16, 23, 27]. It is no coincidence that graph theory has been independently discovered many times, since it may quite properly be regarded as an area of applied mathematics.

Euler (1707-1782) became the father of graph theory. In 1847 Kirchoff developed the theory of trees, in order to solve the system of simultaneous linear equations, which give the current in each branch and each circuit of an electric network. Thus in effect Kirchoff replaced each electrical network by its underlying graph and showed that it is not necessary to consider every cycle in the graph of an electric network separately in order to solve the system of equations.

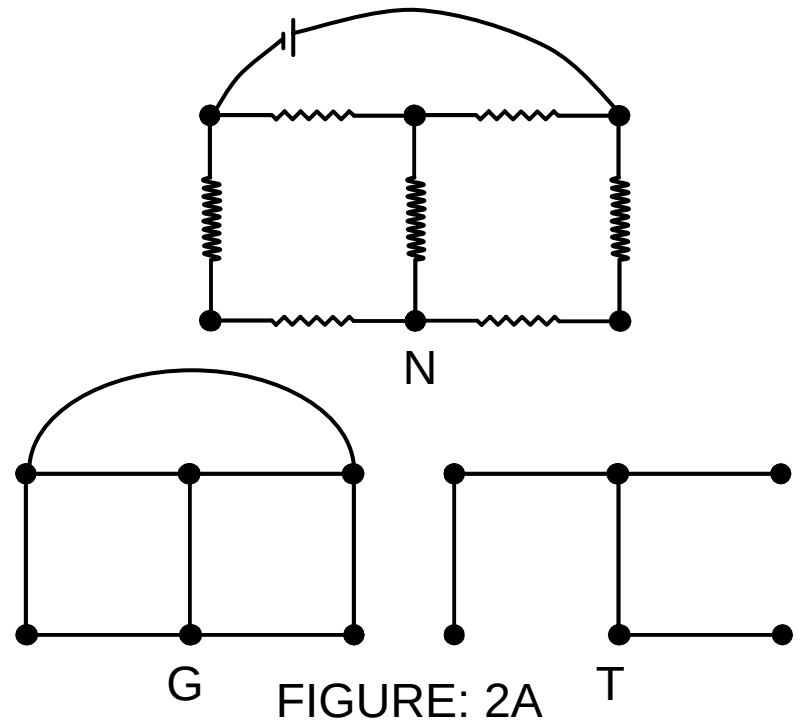

FIGURE: 2A

Instead he pointed out by a simple but powerful construction, which has since become the standard procedure that the independent cycles of a graph determined by any of its spanning trees will suffice. A contrived electrical network N, its under lying graph G and a spanning tree T are shown in the following figure [22-23].

In 1857 Cayley discovered the important class of graphs called trees by considering the changes of variables in the differential calculus. Later he was engaged in enumerating the isomers of the saturated hydrocarbons $C_n$ $H_{2n+2}$, with a given number n of carbon atoms as shown in the Figure.

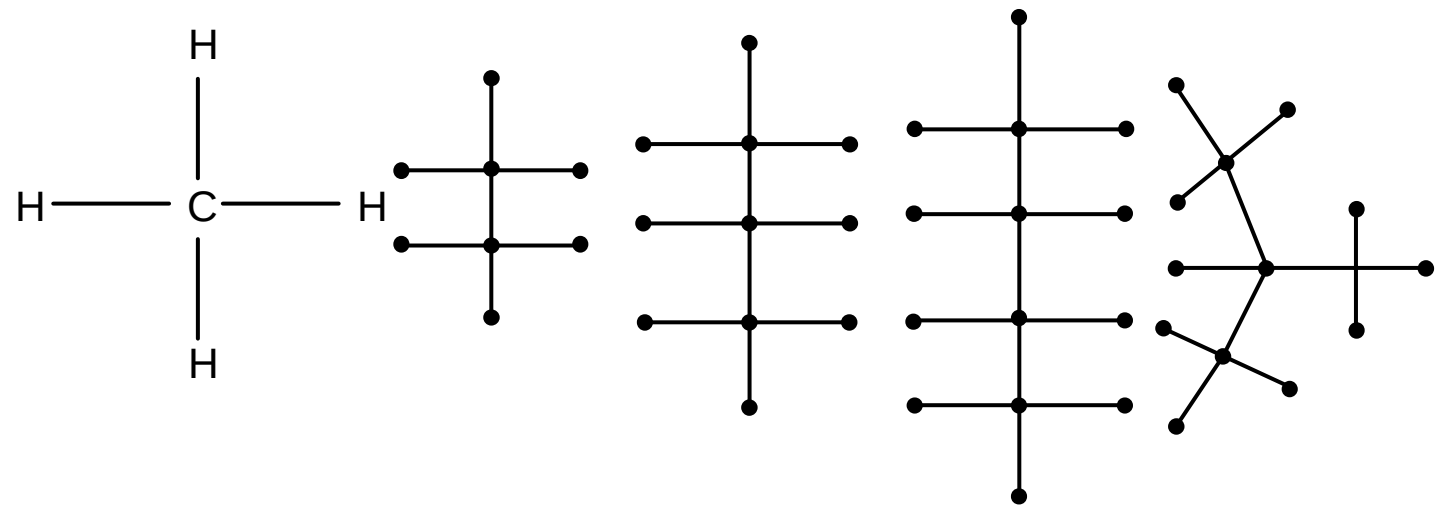

FIGURE: 2.B

Cayley restated the problem abstractly, find the number of tress with p points in which every point has degree 1 or 4. He did not immediately succeed in solving this and so he altered the problem until he was able to enumerate; rooted trees (in which one point is distinguished from the others), with point of degree at most 4 and finally the chemical problem of trees in which every point has degree 1 or 4.

Jordan in 1869 independently discovered trees as a purely mathematical discipline and Sylvester 1882 wrote that Jordan did so without having any suspicion of its bearing on modern chemical doctrine.

The most famous problem in graph theory and perhaps in all of mathematics is the celebrated four color conjecture. The remarkable problem can be explained in five minutes by any mathematician to the so called man in the street. At the end of the explanation both will understand the problem but neither will be able to solve it.

The following quotation from the historical article which state the Four color conjecture and describe its role.

Any map on a plane or the surface of a sphere can be colored with only four colors so that no two adjacent countries have the same color.

Each country must consists of a single connected region and adjacent countries are those having a boundary line (not merely a single point) in common.

The conjecture has acted as a catalyst in the branch of mathematics known as combinatorial topology and is closely related to the currently fashionable field of graph theory. Although the computer oriented proof of [2, 29] settled the conjecture in 1976 and has stood a test of time, a theoretical proof of the four colour problem is still to be found.

Lewin the psychologist proposed in 1936 that the life span of an individual be represented by a planar map. In such a map the regions would represent the various activities of a person such as his work environment, his home and his hobbies. It was pointed out that Lewin was actually dealing with graphs as indicated by Figure 2.C.

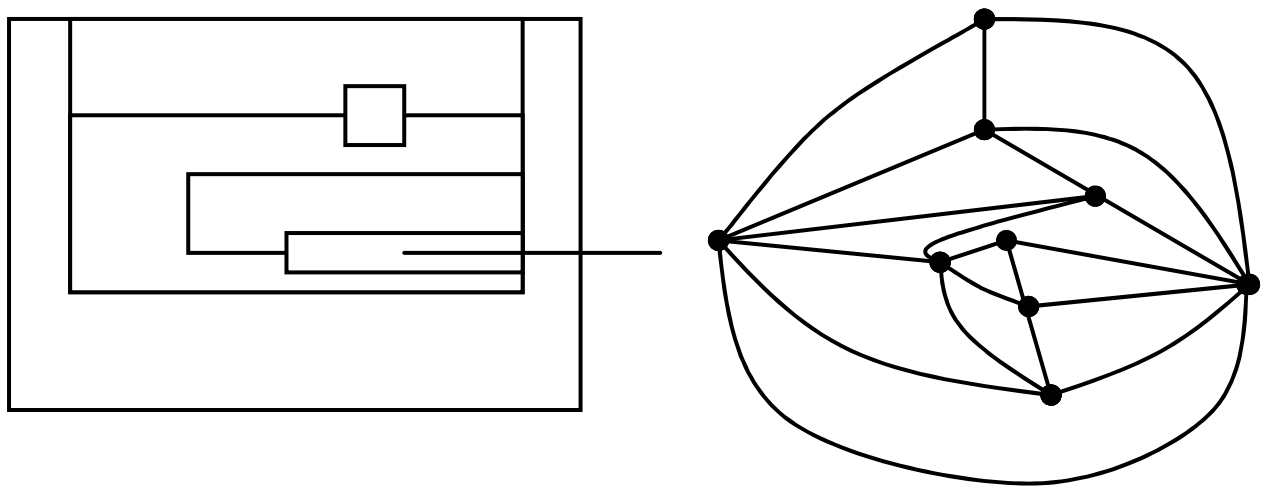

FIGURE: 2.C

This viewpoint led the psychologists at the Research center for Group Dynamics to another psychological interpretation of a graph in which people are represented by points and interpersonal relations by lines. Such relations include love, hate, communication and power. In fact it was precisely this approach which led the author to a personal discovery of graph theory, aided and abetted by psychologists L. Festinger and D. Cartwright.

The world of Theoretical physics discovered graph theory for its own purposes more than once. In the study of statistical mechanics by Uhlenbeck the points stand for molecules and two adjacent points indicate nearest neighbor interaction of some physical kind, for example magnetic attraction or repulsion. In a similar interpretation by Lee and Yang the points stand for small cubes in Euclidean space where each cube may or may not be occupied by a molecule.

Then two points are adjacent whenever both spaces are occupied. Another aspect of physics employs graph theory rather as pictorial device. Feynmann proposed the diagram in which the points represent physical particles and the lines represent paths of the particles after collisions.

The study of Markov chains in probability theory involves directed graphs in the sense that events are represented by points and a directed line from one point to another indicates a positive probability of direct succession of these two events. This is made explicit in which Markov chain is defined as a network with the sum of the values of the directed lines from each point equal to 2. A similar representation of a directed graph arises in that part of numerical analysis involving matrix inversion and the calculation of eigen values.

A square matrix is given preferable sparse and a directed graph is associated with it in the following ways. The points denote the index of the rows and columns of the given matrix and there is a directed line from point i to point j whenever the i, j entry of the matrix is nonzero. The similarity between this approach and that for Markov chains in immediate.

Thus finally in the 21st century the graph theory has been fully exploited by fuzzy theory. The causal structure of fuzzy cognitive maps from sample data [33-34] mainly uses the notion of fuzzy signed directed graphs with feedback. Thus the use of graph theory especially in the field of applications of fuzzy theory is a grand one for most of analysis of unsupervised data are very successfully carried out by the use of Fuzzy Cognitive Maps (FCMs) which is one of the very few tools which can give the hidden pattern of the dynamical system.

The study of Combined Fuzzy Cognitive Maps (CFCMs) mainly uses the concept of digraphs.

The directed graphs or the diagraphs are used in the representation of Binary relations on a single set. Thus one of the forms of representation of a fuzzy relation R(X, X) is represented by the digraph.

Let X = {1, 2, 3, 4, 5} and R (X, X) the binary relation on X defined by the following membership matrix

| | 1 | 2 | 3 | 4 | 5 |
|---|---|---|---|---|---|
| 1 | 0.2 | 0 | 0 | 0.1 | 0.6 |
| 2 | 0 | 0.8 | 0.4 | 0 | 0 |
| 3 | 0.3 | 0 | 0.9 | 0.2 | 0 |
| 4 | 0 | 0 | 0.2 | 0.9 | 0 |
| 5 | 0 | 0.8 | 0 | 0.5 | 0 |

The related graph for the binary relation.

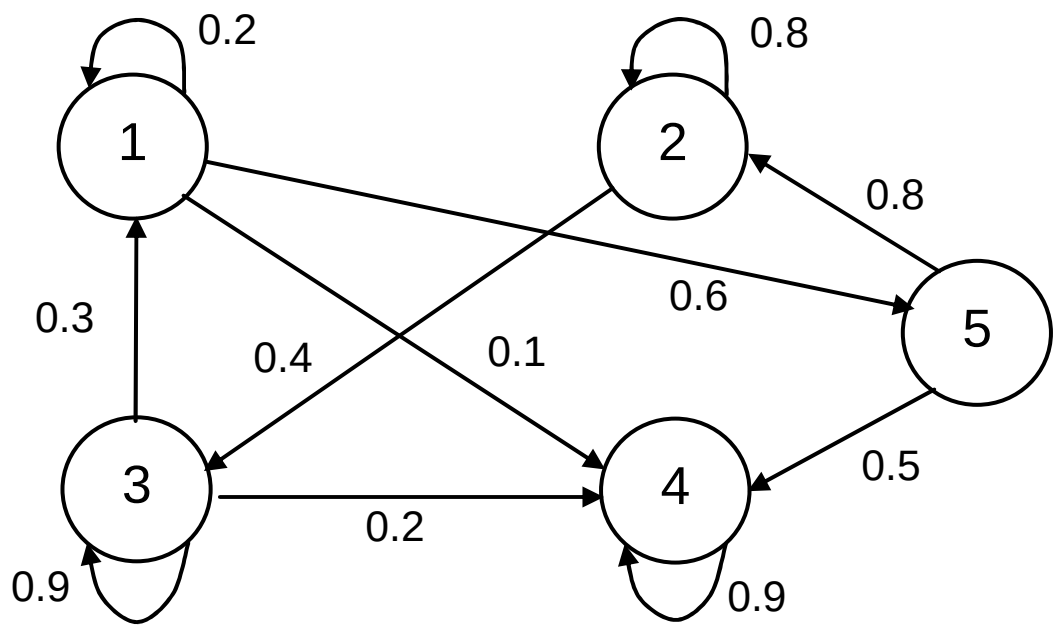

FIGURE: 2.D

Further in the description of the fuzzy compatibility relations also the graphs were used.

The notion of graph theory was used in describing the fuzzy ordering relations. In fact one can say that the graph theory method was more simple and an easy representation even by a lay man.

Thus we can say whenever the data had a fuzzy matrix representation it is bound to get the digraph representation. Also in the application side graph theory has been scrupulously used in the description of automaton and semi automaton i.e., in finite machines we do not study in this direction in this book.

We devote an entire section in this chapter to see applications of graphs in fuzzy models.

This chapter has eight sections. The first seven sections give some basic results on graphs theory to make this book a self contained one. The last section contains a complete application of graphs to fuzzy models.

## 2.1 Some basics on Graphs

In this section we recall the definition of graphs and some of its properties. [5, 23]

**DEFINITION 2.1.1:** *A graph G consists of a finite non empty set V = V (G) of p points. (vertex, node, junction O-simplex elements) together with a prescribed set X of q unordered pairs of distinct points of V.*

*Each pair x = {u, v} of points in X is a line (edge, arc, branch, 1-simplex elements) of G and x is said to join u and v.*

*We write x = uv and say that u and v are adjacent points (some times denoted as u adj v); point u and line x are incident with each other as arc v and x. E(G) will denote the edges or lines of G.*

*If two distinct lines x and y are incident with a common point, then they are adjacent lines. A graph with p points and q lines is called a (p, q) graph.*

*Clearly (1, 0) graph is trivial. A graph is represented always by a diagram and we refer to it as the graph.*

The graph in figure 2.1.1 is totally disconnected.

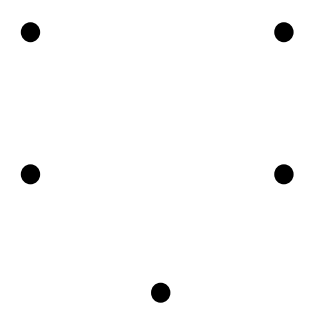

FIGURE: 2.1.1

The graph in figure 2.1.2 is disconnected.

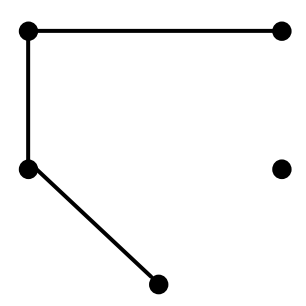

FIGURE: 2.1.2

The graph with four lines in figure 2.1.3 is a path.

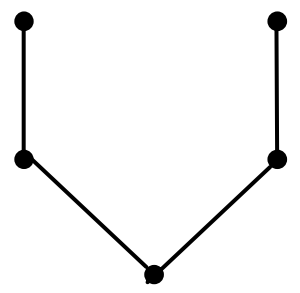

FIGURE: 2.1.3

The graph in figure 2.1.4 is a connected graph.

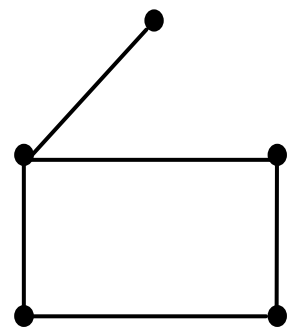

FIGURE: 2.1.4

The graph in figure 2.1.5 is a complete graph.

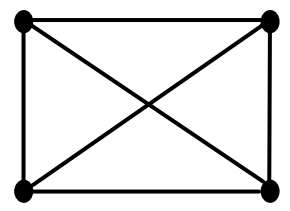

FIGURE: 2.1.5

The graph in figure 2.1.6 with four lines is a cycle.

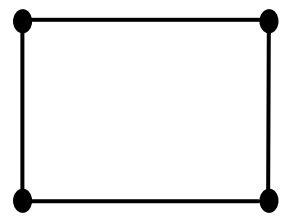

FIGURE: 2.1.6

It is important to note that in a graph if any two lines intersect it is not essential that their intersection is a point of the graph.

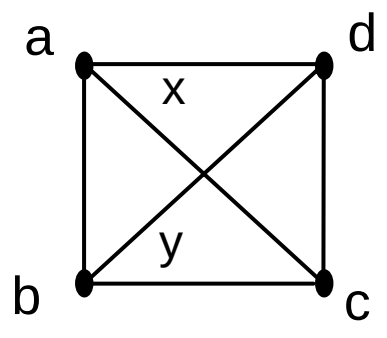

FIGURE: 2.1.7

i.e., the lines x and y intersect in the diagram their intersection is not a point of the graph.

Recall in the graph in figure 2.1.8.

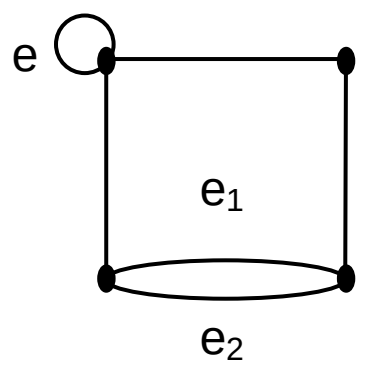

FIGURE: 2.1.8

e is loop and $\{e_1, e_2\}$ is a set of multiple edges. Thus a graph with loops and multiple edges will be known by some authors as pseudo graphs. But we shall specify graphs with multiple edges and graphs with loops distinctly. A graph is simple if it has no loops and multiple edges.

**DEFINITION 2.1.2:** *A graph is called finite if both V(G) and E(G) are finite. A graph that is not finite is called infinite.*

We shall in this book deal only with graphs which are finite. N(G) and m(G) are the number of vertices and edges of the graph G, respectively. The number n(G) is called the order of G and m(G) is the size of G.

A graph is said to be labeled if its n vertices are distinguished from one another by labels such as $v_1, v_2, \ldots, v_n$. A graph with 3 vertices $v_1, v_2, v_3$ is a labelled as follows:

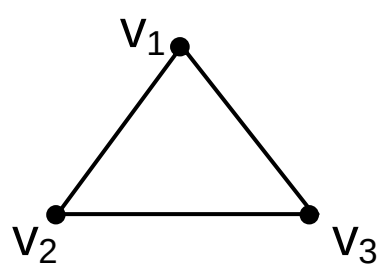

FIGURE: 2.1.9

Thus a graph G is represented as an ordered triple G = (V(G), E(G), $I_G$) where V(G) the vertices or nodes of G E(G) the edges or lines of G. If for any edge e of G $I_G(e) = \{u, v\}$ we write $I_G(e) = uv$.

**DEFINITION 2.1.3:** *A graph isomorphism between two graphs G = (V(G), E(G), $I_G$) and H = (V(H), E(H), $I_H$) written $G \cong H$ is a pair ($\phi$, $\theta$) where*

$$\phi : V(G) \rightarrow V(H) \text{ and}$$
$$\theta : E(G) \rightarrow E(H)$$

*are bijections with the property that* $I_G(e) = \{u, v\}$ *if and only if* $I_H(\theta(e)) = \{\phi(u), \phi(v)\}$. *If* $(\phi, \theta)$ *is a graph isomorphism the pair of inverse mappings* $(\phi^{-1}, \theta^{-1})$ *is also a graph isomorphism.*

*Note that the bijection* $\phi$ *satisfies the condition that u and v are end vertices of an edge e of G if and only if* $\phi(u)$ *and* $\phi(v)$ *are end vertices of the edge* $\phi(e)$ *in H.*

*If graphs G and H are simple a bijection*

$$\phi: V(G) \rightarrow V(H)$$

*such that u and v are adjacent in G if and only if* $\phi(u)$ *and* $\phi(v)$ *are adjacent in H induces a bijection*

$$\theta : E(G) \rightarrow E(H)$$

*satisfying the condition that* $I_G(e) = \{u, v\}$ *if and only if* $I_H(\theta(e)) = \{\phi(u), \phi(v)\}$. *Hence* $\phi$ *itself is referred to as an isomorphism in the case of simple graphs G and H. Thus if G and H are simple graphs an isomorphism from G to H is a bijection* $\phi: V(G) \rightarrow V(H)$ *such that u and v are adjacent in G if and only if* $\phi(u)$ *and* $\phi(v)$ *are adjacent in H.*

*Hence* $uv \in E(G)$ *implies* $\phi(u)\,\phi(v) \in E(H)$.

**DEFINITION 2.1.4:** *A simple graph H is said to be complete if every pair of distinct vertices of G are adjacent in G.*

Any two complete graphs each of a set of n vertices are isomorphic, each such graph is denoted by $K_n$

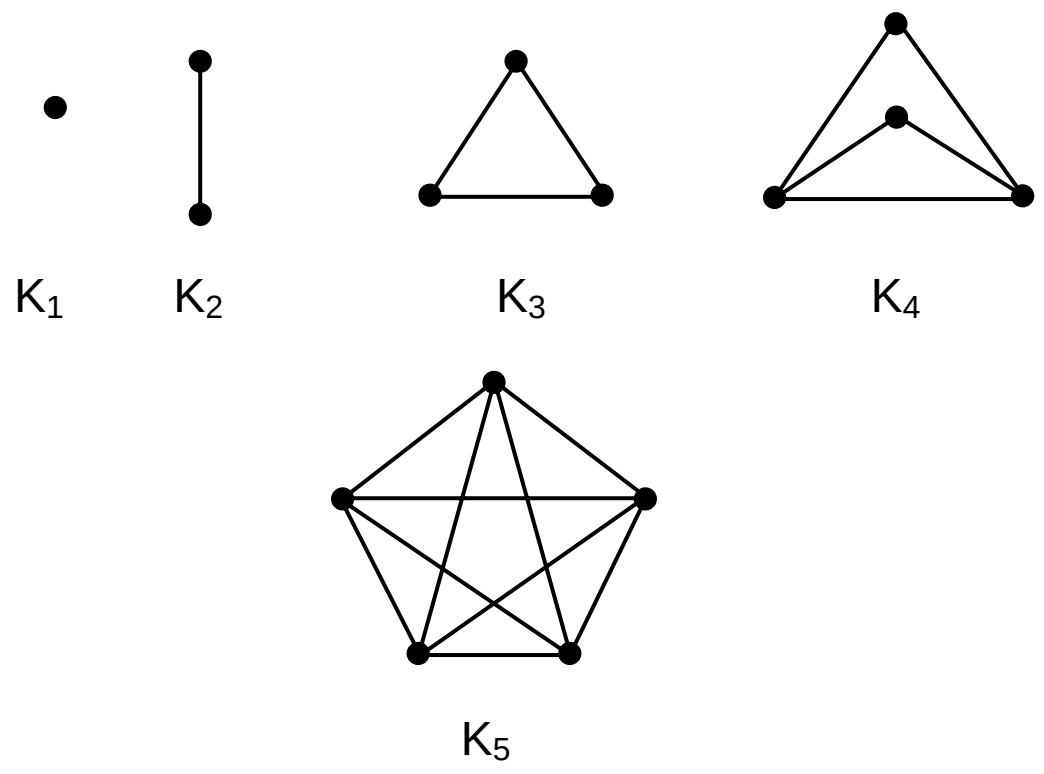

FIGURE: 2.1.10

A simple graph with n vertices can have atmost $\binom{n}{2}$ = n (n – 1)/2 edges. $K_n$ has the maximum number of edges among all simple graphs with n vertices.

We may have graphs which has no edges at all such graphs are called totally disconnected graphs

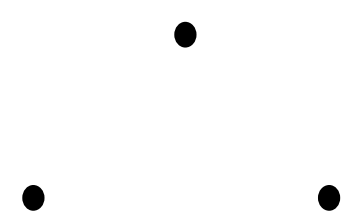

FIGURE: 2.1.11

A totally disconnected graph with 3 vertices is shown in figure 2.1.11.

Thus we have any simple graph G with n vertices have

$$0 \leq |E(G)| \leq \frac{n(n-1)}{2}$$

where $|E(G)|$ is the number edges of G.

A trivial graph is a graph with a singleton set with no edges.

**DEFINITION 2.1.5:** *A graph is bipartite if its vertex set can be partitioned into two nonempty subsets X and Y such that each edge of G has one end in X and the other end in Y. The pair (X, Y) is called a bipartition of the bipartite graph. The bipartite graph G with bipartition (X, Y) is denoted by G(X, Y). A simple bipartite graph is complete if each vertex of X is adjacent to all the vertices of Y.*

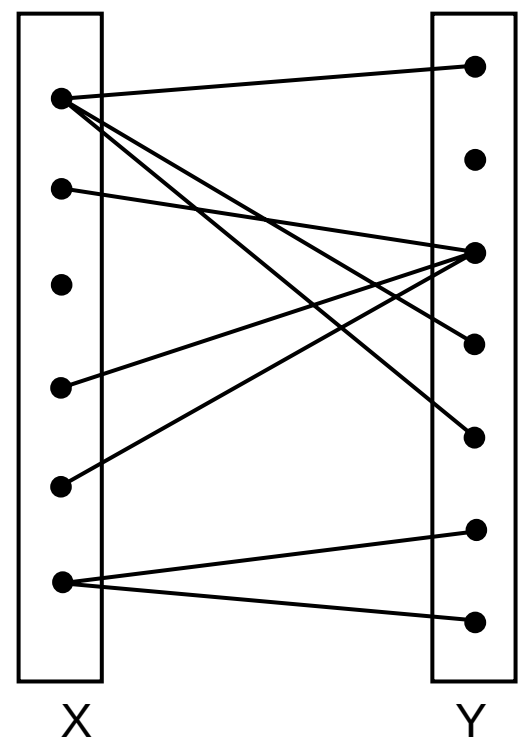

FIGURE: 2.1.12

A bipartite Graphs. The bipartite graph G (X, Y) is complete denoted by $K_{t,s}$ if |X| = t and |Y| = s.

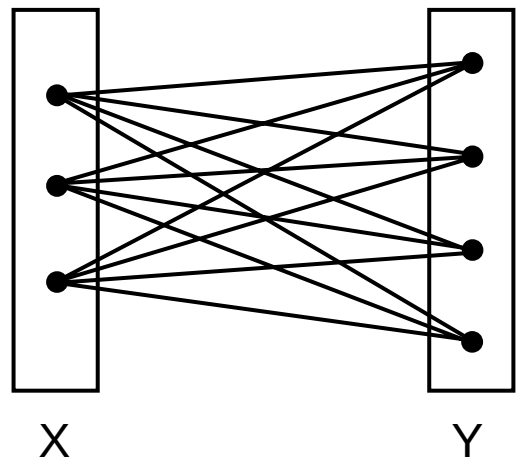

FIGURE: 2.1.13

The graph $K_{3,4}$ A compete bipartite graph $K_{1,s}$ is called a star.
The star graph. $K_{1,8}$

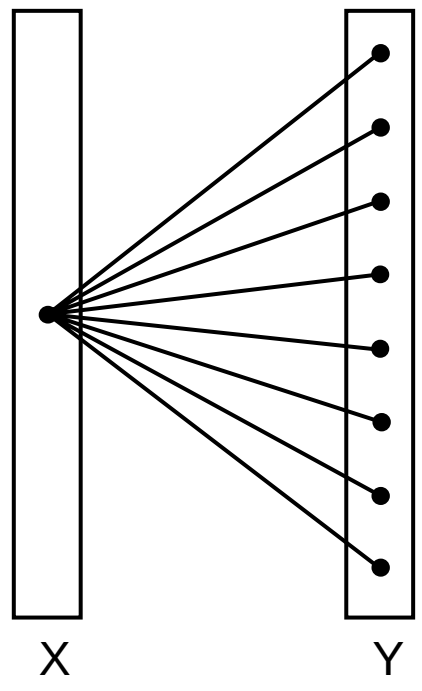

FIGURE: 2.1.14

**DEFINITION 2.1.6:** *Let G be a simple graph. Then the complement $G^c$ of G is defined by taking $V(G^c) = V(G)$ and making two vertices u and v adjacent in $G^c$*

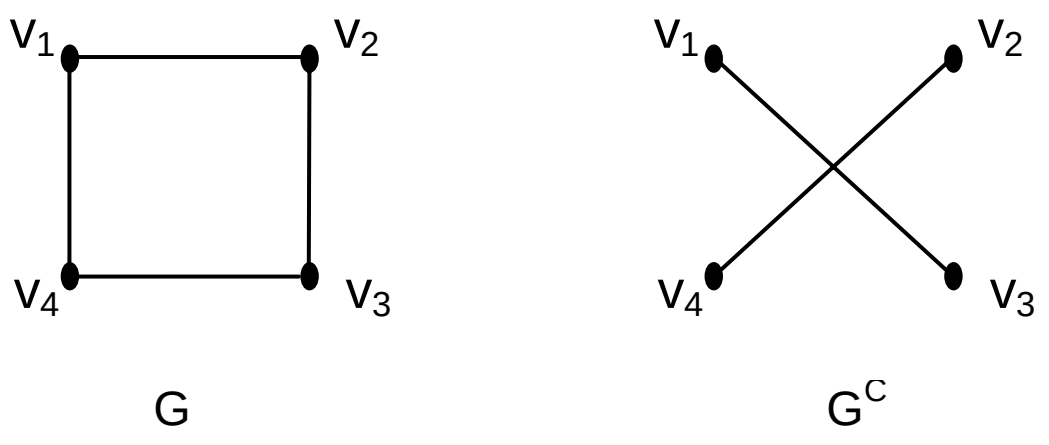

FIGURE: 2.1.15

*$G^c$ is also a simple graph and $(G^c)^c = G$.*

**DEFINITION 2.1.7:** *A simple graph G is self complementary if $G \cong G^c$.*

**DEFINITION 2.1.8:** *A graph H is called the subgraph of G if $V(H) \subseteq V(G)$, $E(H) \subseteq E(G)$ and $I_H$ is the restriction of $I_G$ to E(H). If H is a subgraph of G then G is said to be a supergraph of H. A subgraph H of a graph G is a proper subgraph of G if either $V(H) \neq V(G)$ or $E(H) \neq E(G)$.*

A subgraph H of G is said to be an induced subgraph of G if each edge of G having its ends in V(H) is also an edge of H. A subgraph H of G is a spanning subgraph of G if V(H) = V(G). The induced subgraph of G with vertex set $S \subseteq V(G)$ is called the subgraph of G induced by S and is denoted by G(S). Let E' be a subset of E and let S denote the subset of V consisting of all end vertices in G of edges in E'.

Then the graph (S, E' $I_G \mid$ E' ) is the subgraph of G induced by the edge set E' of G. It is denoted by G(E' ).

**DEFINITION 2.1.9:** *A clique of G is a complete subgraph of G. A clique of G is a maximal clique of G if it is not properly contained in another clique of G.*

**DEFINITION 2.1.10:** *Deletion of vertices and edges in a graph. Let G be a graph. X be proper subset of the vertex set V and E' a subset of E. The subgraph G[V \ X] is said to be obtained from G by deletion of X this subgraph is denoted by G\X. If X = {u}, G\X is simply denoted G-u. The spanning subgraph of G with edge set E/E' is the subgraph obtained from G by deleting the edge subset E'. This subgraph is denoted by G-E' whenever E' = {e}, G-E' is simply denoted by G-e.*

It is important to note if a vertex is deleted from graph G all the edges incident to it are deleted where as if an edge is deleted from a graph G it does not affect the vertices of G.

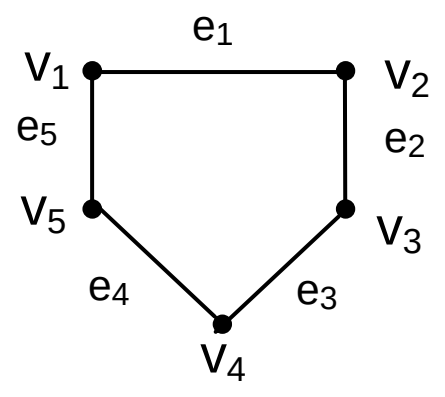

FIGURE: 2.1.16

If $e_1$ is removed in figure 2.1.16.

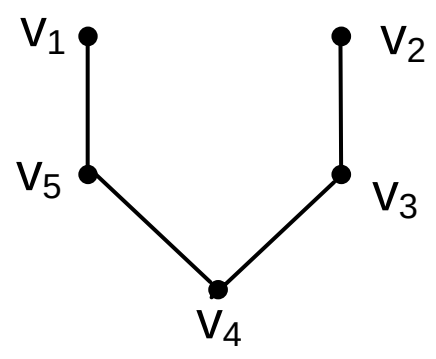

FIGURE: 2.1.17

G \ {$e_1$} is the subgraph given in figure 2.1.17.

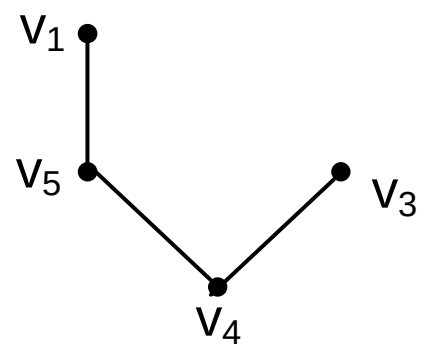

FIGURE: 2.1.18

Subgraph G – {$v_2$} is given in figure 2.1.18.

**DEFINITION 2.1.11:** *Let G be a graph and v $\in$ V. The number of edges incident at v in G is called the degree (or valency) of the vertex v in G and is denoted by $d_G(v)$, or simply d(v) when G requires no explicit reference. A loop of v is to be counted twice in computing the degree of v. The minimum (respectively maximum) of the degrees of the vertices of a graph G is denoted by $\delta(G)$ or $\delta$(respectively $\Delta$ (G) or $\Delta$). A graph G is called k-regular if every vertex of G has degree k. A graph is said to be regular if it is k-regular for some nonnegative integer k. A 3-regular graph is called a cubic graph.*

**DEFINITION 2.1.12:** *A spanning 1-regular subgraph of G is called a 1-factor or a perfect matching of G.*

**DEFINITION 2.1.13:** *A vertex of degree 0 is known as an isolated vertex of G. A vertex of degree 1 is called a pendant vertex of G, where as the unique edge of G incident to such a vertex of G is a pendant edge of G. A sequence formed by the degrees of vertices of G is called a degree sequence of G.*

The first famous theorem due to Leonard Euler (1707-1783) gives:

"The sum of the degrees of the vertices of a graph is equal to twice the number of its edges".

Clearly if d = ($d_1, \ldots, d_n$) is degree sequence of G then $\sum_{i=1}^{n} d_i = 2m$ where n and m are the order and size of G respectively.

In any graph G the number of vertices of odd degrees is even. Now we proceed onto recall briefly the notion of path and connectedness in graphs.

## 2.2. More Properties on Graphs

One of the elementary but interesting property about graphs is connectedness and disconnectedness.

**DEFINITION 2.2.1:** *A walk of a graph G is an alternating sequence of points and lines $v_0$, $e_1$, $v_1, \ldots, v_{n-1}$ $e_n$ $v_n$ beginning and ending with the points / vertices, in which each line is incident with the two points immediately proceeding and following it. This walk joins $v_0$ and $v_n$ and may also be denoted by $v_0$ $v_1, \ldots, v_n$. It is some times called $v_0$ -$v_n$ walk. It is closed if $v_0 = v_n$ and is open other wise.*

*A walk is called a trial if all the edges appearing in the walk are distinct and a path if all the points are distinct.*

If the walk is closed then it is a cycle provided its n points are distinct and $n \geq 3$. In the labeled graph G given in the following figure 2.2.1.

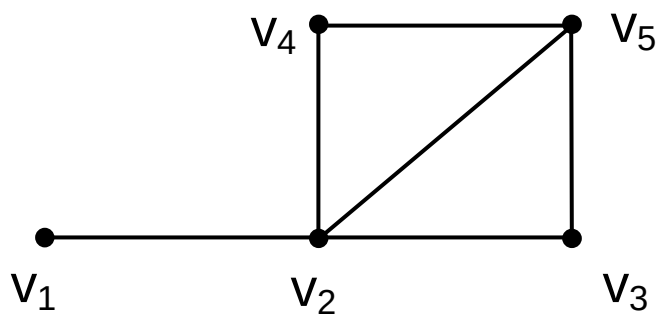

FIGURE: 2.2.1

$v_1$ $v_2$ $v_5$ $v_2$ $v_3$ is a walk which is not a trial and $v_1$ $v_2$ $v_5$ $v_4$ $v_2$ $v_3$ is a trial which is not a path; $v_1$ $v_2$,$v_5$ $v_4$ is a path and $v_2$ $v_4$,$v_5$ $v_2$ is a cycle.

We denote by $C_n$ the graph consisting of a cycle with n points and by $P_n$ a path with n-points $C_3$ is called a triangle.

A graph is connected if every pair of points are joined by a path. A maximal connected subgraph of G is called a connected component or simply a component of G. Thus a disconnected graph has at least two components. The graphs in the following figure have eight components.

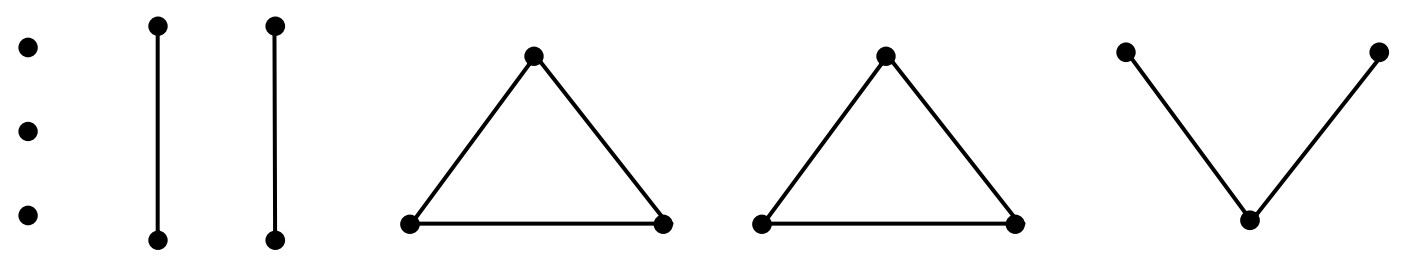

FIGURE: 2.2.2

The length of a walk $v_0$ $v_1$… $v_n$ is n the number of occurrences of lines in it. The girth of the graph G denoted by g(G) is the length of a shortest cycle if any in G; the circumference c(G) is the length of any longest cycle. These terms are undefined or has no meaning if G has no cycles.

The distance d(u, v) between two points u and v in G is the length of a shortest path joining them if any, otherwise d (u, v) = ∞. In a connected graph distance is a metric i.e., for all points u, v and w

1. d(u, v) ≥ 0 with
   d (u, v) = 0 if and only if u = v
2. d (u, v) = d(v, u)
3. d (u, v) + d (v, ω ≥ d (u, ω))
4. A shortest u - v path is often called a geodesic. The diameter d(G) of a connected graph G is the length of any longest geodesic. The graph G given below:

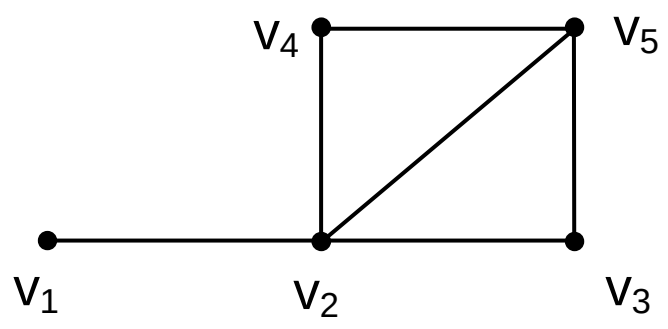

FIGURE: 2.2.3

has girth G = 3, circumference c(G) = 4 and diameter d(G) = 2.

The square $G^2$ of a graph G has V (G) = $V(G^2)$ with u, v adjacent in $G^2$ whenever $d(u, v) \leq 2$ in G. The powers $G^3$, $G^4$,…, of G are defined similarly.

***Result:*** For any graph with six points G or $G^2$ contains a triangle.

A graph G is called locally connected if for every vertex v of G, the neighbor set of v in G, $N_G$ (v) is connected. (The vertex u is a neighbor of v in G if uv is an edge of G and $u \neq v$. The set of all neighbors of v is the open neighborhood of v or the neighbor set of v and is denoted by N(v); the set $N(v) \cup \{v\} = N[v]$ is the closed neighborhood of v in G.) A cycle is odd or even according as the length of the cycle is odd or even.

***Result:*** A graph is bipartite if and only if it contains no odd cycles.

**DEFINITION 2.2.2:** *An automorphism of a graph G is an isomorphism of G into itself.*

Thus two simple graphs G and H are isomorphic if and only if there exists a bijection $\theta : V(G) \to V(H)$ such that uv is an edge of G if and only if $\phi(u)\, \phi(v)$ is an edge of H. In this case $\phi$ is called an isomorphism of G onto H. $\tau(G)$ will denote the set of automorphisms of G and $\tau(G)$ is a group.

***Result:*** The set $\tau(G)$ of all automorphisms of a simple graph G is a group with respect to the composition 'o' of mappings as the group operation.

***Result:*** If G is any simple graph $\tau(G) = \tau(G^c)$.

We now proceed on to recall the definition of line graph. Let G be a loopless graph. The graph L(G) is constructed in the following way.

The vertex set L(G) is in one to one correspondence with the edge set of G and two vertices of L(G) are joined by an edge if and only if the corresponding edges of G are adjacent in G. The graph L(G) is called the line graph or the edge graph of G. The following figure gives conversion of the graph G to the graph L(G).

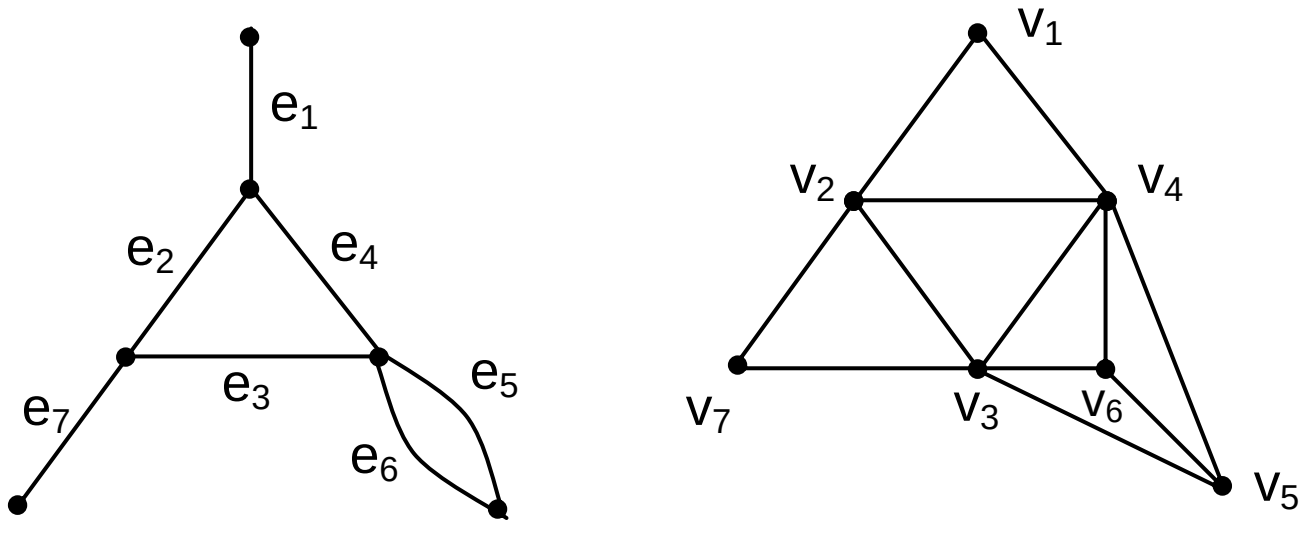

FIGURE: 2.2.4

Some simple properties of the line graph L(G) of a graph G is as follows:

1. G is connected if and only if L(G) is connected.
2. If H is a subgraph of G then L(H) is a subgraph of L(G).
3. The edges incident of a vertex of G give rise to a maximal complete subgraph of L(G).
4. If e is an edge of G joining u and v then the degree of e in L(G) is the same as the number of edges of G adjacent to e in G.

***Result:*** The line graph of a simple graph G is a path if and only if G is a path.

***Result:*** If the simple graphs $G_1$ and $G_2$ are isomorphic then $L(G_1)$ and $L(G_2)$ are isomorphic.

***Result:*** Let G and G' be simple connected graphs with isomorphic line graphs. Then G and G' are isomorphic unless one of them is $K_{13}$ and the other is $K_3$.

A graph H is called a forbidden subgraph for a property P of graphs if it satisfies the following condition. If a graph G has property P, then G cannot contain an induced subgraph isomorphic to H.

***Result:*** If G is a line graph, then $K_{1,\,3}$ is a forbidden subgraph of G.

Now we proceed on to define union of two graphs.

**DEFINITION 2.2.3:** *Let $G = (V(G), (E(G))$ where $V(G) = V_1(G_1) \cup V_2(G_2)$ and $E = E_1 \cup E_2$ i.e., $E(G) = E_1(G_1) \cup E_2(G_2)$ is called the union of the two graphs*

$G_1 = (V_1(G_1), E_1(G_1))$ *and*

$G_2 = (V_2(G_2), E_2(G_2))$ *and is denoted by* $G_1 \cup G_2$

*when $G_1$ and $G_2$ are vertex disjoint $G_1 \cup G_2$ is denoted by $G_1 + G_2$ and is called the sum of the graphs $G_1$ and $G_2$.*

We just recall the definition of intersection of two graphs.

**DEFINITION 2.2.4:** *If $G_1 = (V_1(G_1), E_1(G_1))$ and $G_2 = (V_2(G_2), E_2(G_2))$. If $V_1(G_1) \cap V_2(G_2) \neq \phi$ the graph $G = G_1 \cap G_2 = (V(G), E(G))$ where $V(G) = V_1(G_1) \cap V_2(G_2)$ and $E = E_1 \cap E_2$ is the intersection of $G_1$ and $G_2$ and is written as $G_1 \cap G_2$.*

Now for two vertex disjoint graphs $G_1$ and $G_2$ we define join of the two graphs.

**DEFINITION 2.2.5:** *Let $G_1 = (V_1(G_1), E_1(G_1))$ and $G_2 = (V_2(G_2), E_2(G_2))$. Then the join $G_1 \vee G_2$ of $G_1$ and $G_2$ is the super graph of $G_1 + G_2$ in which each vertex of $G_1$ is adjacent to every vertex of $G_2$.*

The following figure illustrates the graph $G_1 \vee G_2$.

Take $G_1$ and G2 as in figure 2.2.5 and figure 2.1.16.

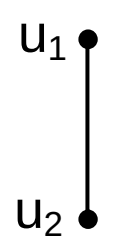

FIGURE: 2.2.5

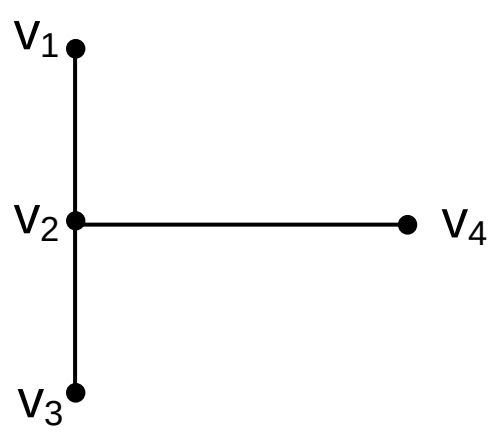

FIGURE: 2.2.6

$G_1 \vee G_2$:

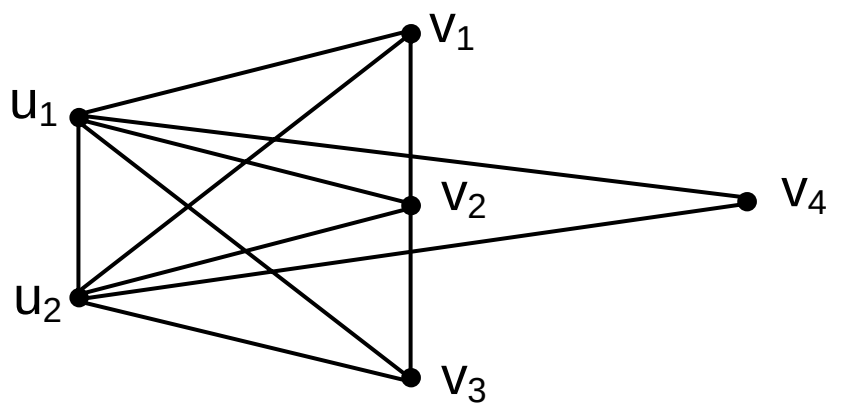

FIGURE: 2.2.7

Now we proceed onto recall the definition of wheel $W_n$ .

If $G_1 = K_1$ and $G_2 = C_n$ then $G_1 \vee G_2$ is called the wheel $W_n$. $W_5$ is shown in the following figure:

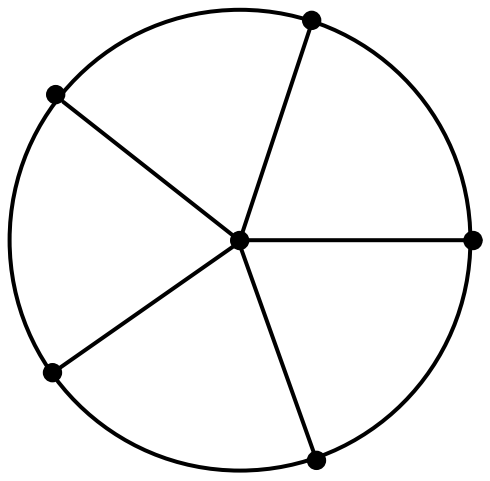

FIGURE: 2.2.8

Now we recall the definition of the Cartesian product of $G_1 \times G_2$ of two graphs $G_1$ and $G_2$.

**DEFINITION 2.2.6:** *The Cartesian product $G_1 \times G_2$ of two graphs $G_1$ and $G_2$ is the simple graph with $V_1 \times V_2$ as its vertex set and two vertices $(u_1, v_1)$ and $(u_2, v_2)$ are adjacent in $G_1 \times G_2$ if and only if either $u_1 = u_2$ and $v_1$ is adjacent to $v_2$ in $G_2$ or $u_1$ is adjacent to $u_2$ in $G_1$ and $v_1 = v_2$.*

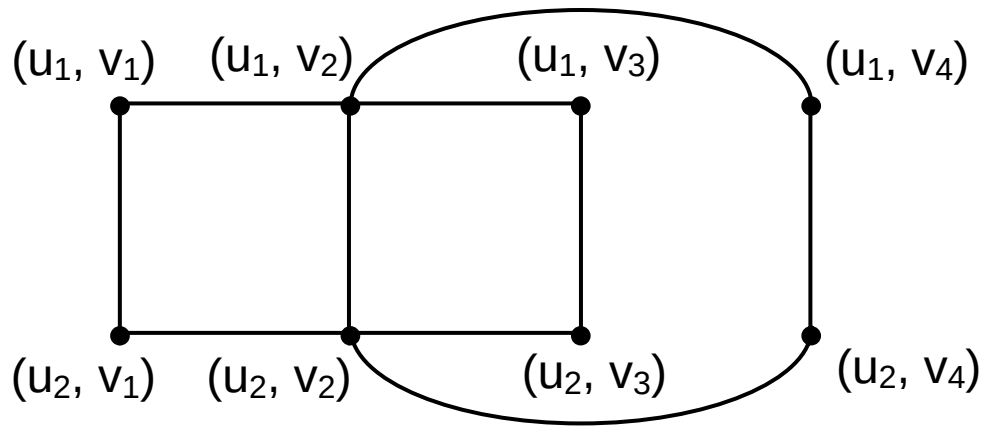

$G_1 \times G_2$

FIGURE: 2.2.9

## 2.3 Directed Graphs

In this section we define the notion of directed graphs and recall some of the properties, which we will be using for defining the notion of their neutrosophic analogue.

**DEFINITION 2.3.1:** *A directed graph D is an ordered triple (V(D), A(D), $I_D$) where V(D) is a nonempty set called the set of vertices of D; A(D) is a set disjoint from V(D), called the set of arcs of D and $I_D$ is an incidence map that associates with each arc of D an ordered pair of vertices of D. If a is an arc of D, and $I_D$ (a) = (u, v), u is called the tail of a and v is the head of a . The arc a is said to join v and u, u and v are called the ends of a. A directed graph is also called a digraph.*

With each digraph D, we can associate a graph G (written G(D) when reference of D is needed) on the same vertex set as follows: corresponding to each arc of D there is an edge arc of G with the same end. This graph G is called the underlying graph of the digraph D. Thus every digraph D defines a unique graph G. Conversely given any graph G we can obtain a digraph from G by specifying for each edge of G an order of its end. Such a specification is called an orientation of G.

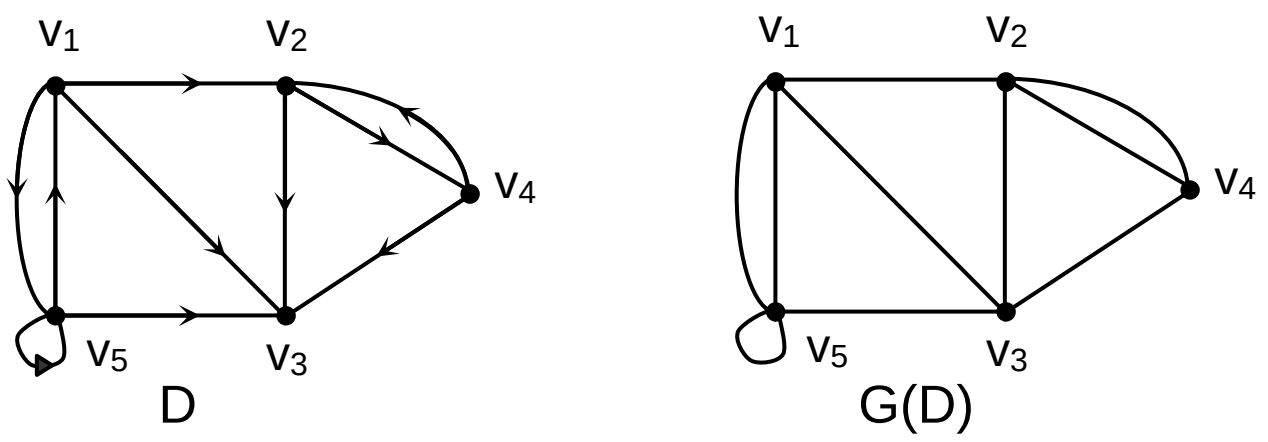

FIGURE: 2.3.1

Digraph D and its underlying graph G. A digraph D' is a subdigraph of a digraph D if V(D') $\subset$ V(D), A(D') $\subseteq$ A(D) and $I_{D'}$ is the restriction of $I_D$ to A(D').

A directed walk joining the vertex $u_0$ to the vertex $u_k$ in D is an alternating sequence w = $u_0$ $a_1$ $u_1$ $a_2$ $u_2$… $a_k$ $u_k$ , $1 \leq i \leq k$, with $a_i$ incident out of $u_{I-1}$ and incident into $u_i$. Directed trails, directed paths, directed cycles and induced subdigraphs are defined analogously as for graphs.

A vertex v is reachable from a vertex u of D if there is a directed path in D from u to v. Two vertices of D are disconnected if each is reachable from the other.

It is easily verified that disconnection is an equivalence relation on the vertex set of D and if the equivalence classes are $V_1, V_2, \ldots, V_w$ the subdigraphs of D induced by $V_1, V_2, \ldots, V_w$ are called dicomponents of D.

A digraph is diconnected if it has exactly one dicomponent. A diconnected digraph is also called a strong digraph.

A digraph is strict if its underlying graph is simple. A digraph D is symmetric if whenever (u, v) is an arc of D, then (u, v) is also an arc of D.

A digraph D is a tournament if its underlying graph is a complete graph. Thus in a tournament for every pair of distinct elements u and v either (u, v) or (v, u) but not both is an arc of D.

Tournaments on 3 or four vertices.

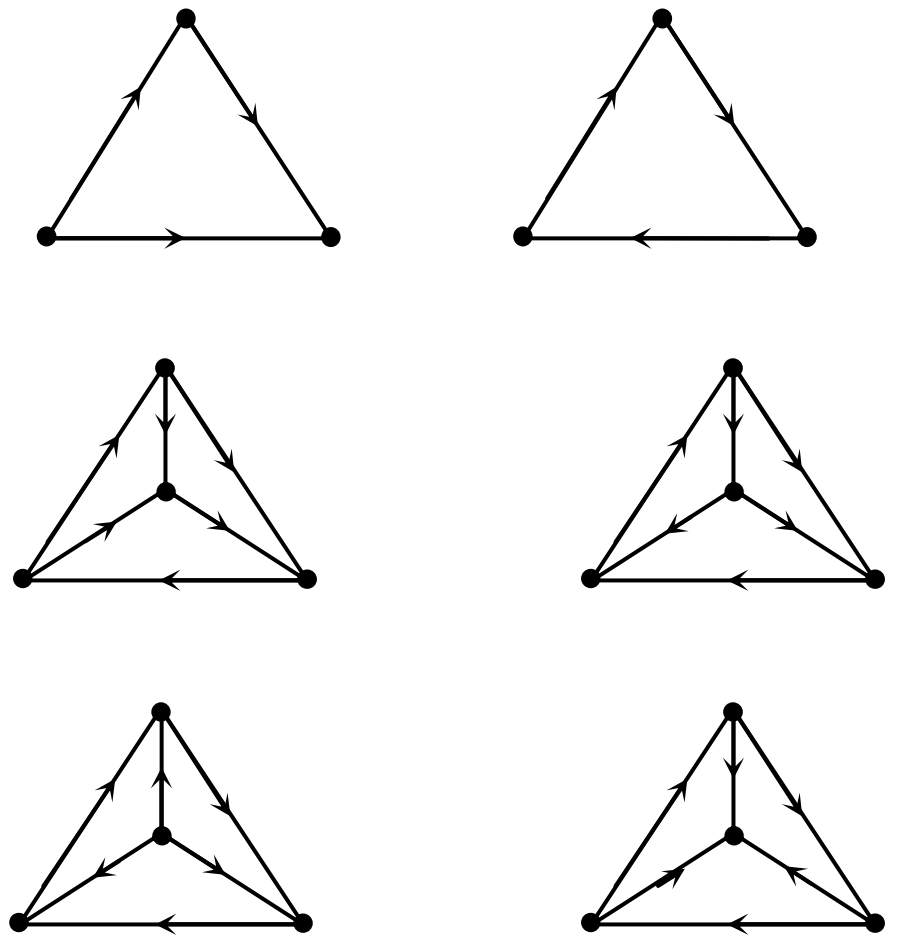

FIGURE: 2.3.2

***Result:*** Every tournament contains a directed Hamilton path. (A directed Hamilton path is a spanning directed path).

Several interesting results in this direction can be had from [43].

**DEFINITION 2.3.2:** *A K- partite graph $K \geq 2$ is a graph G in which V(G) is partitioned into K nonempty subsets $V_1, V_2, \ldots, V_k$ such* *that the induced subgraphs $G[V_1]$, $G[V_2]$,…, $G[V_k]$ are all totally disconnected.*

It is said to be complete if for $i \neq j$ each vertex of $V_i$ is adjacent to every vertex of $V_j$, $1 \leq j$, $j \leq K$. A K-partite tournament is an ordered complete K-partite graph.

The interested reader is requested to refer [23, 43].

Now we proceed on to recall the definition of connectivity. The connectivity of a graph is a measure of its connectedness. Some connected graphs are connected rather loosely in the sense that the deletion of a vertex or an edge from the graph destroys the connectedness of the graph. There are graphs in the other extreme as well such as the complete graph $K_n$, $n \geq 2$ which remain connected even after removal of all but one vertex.

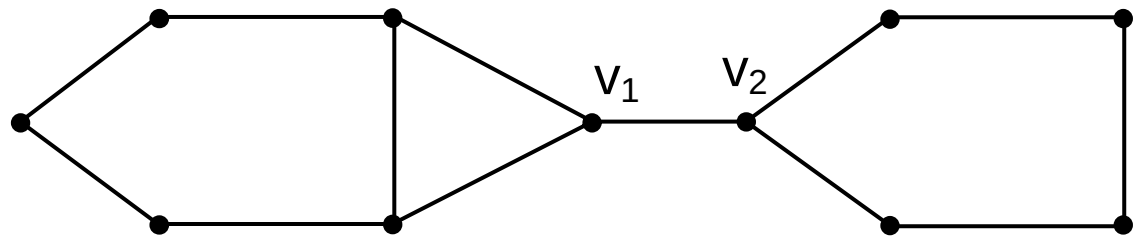

FIGURE: 2.3.3

Removal or deletion of a vertex $v_1$ or $v_2$ from the above graph destroys the connectedness of the graph.

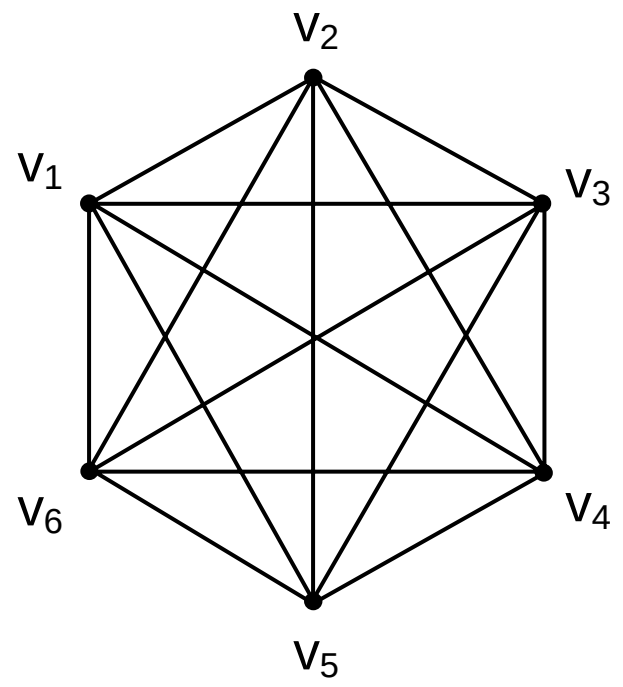

FIGURE: 2.3.4

Even removal of all the five vertices from $K_6$ it remains connected.

Now we proceed on to recollect the definitions of vertex cuts and edge cuts.

**DEFINITION 2.3.3:** *A subset V' of the vertex set V(G) of a connected graph G is a vertex cut of G if G – V' is disconnected, it is K-vertex cut if* $|V'| = K$. *V' is then called a separating set of vertices of G.*

A vertex v of G is a cut vertex of G if {v} is a vertex cut of G.

***Example 2.3.1:*** In the graph

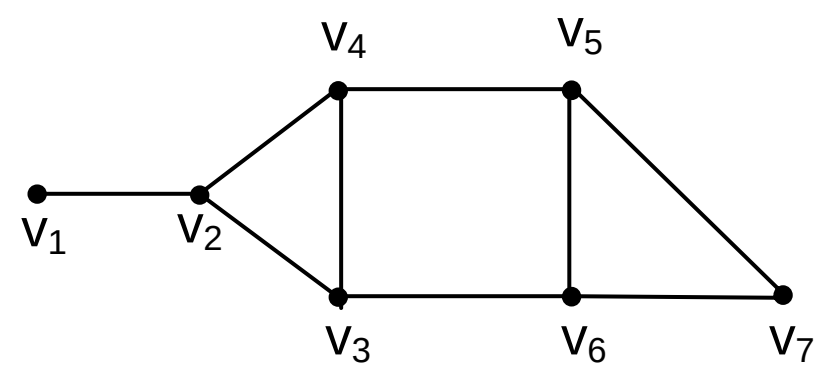

FIGURE: 2.3.5

The set $\{v_2\}$, $\{v_4, v_3\}$ and $\{v_6, v_5\}$ are vertex cuts.

**DEFINITION 2.3.4:** *Let G be a nontrivial connected graph with vertex set V and let S be a non empty subset of V. For* $\overline{S} = V \setminus S$, *let* $\{S, \overline{S}\}$ *denote the set of all edges of G that have one end vertex in S and the other in* $\overline{S}$. *A set of edges of G of the form* $\{S, \overline{S}\}$ *is called an edge cut of G. An edge e is a cut edge of G, if {e} is an edge cut of G. An edge cut of cardinality K is called a K-edge cut of G.*

***Example 2.3.2:*** Consider the graph given by the following figure:

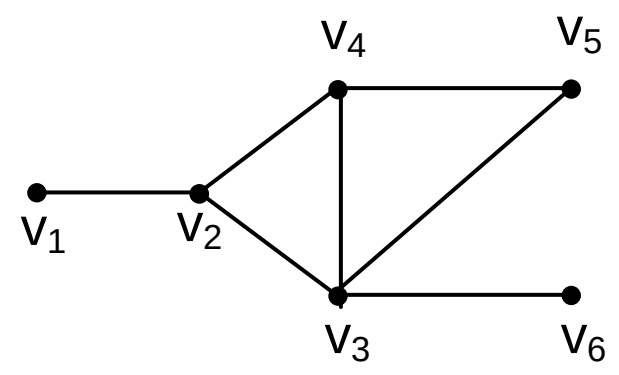

FIGURE: 2.3.6

$v_1 v_2$ and $v_3 v_6$ are cut edges.

Several results in this direction can be had from [3, 5, 23].

**DEFINITION 2.3.5:** *For a nontrivial connected graph G having a pair of non adjacent vertices, the minimum K for which there exists a K- vertex cut is called the vertex connectivity or simply the connectivity of G; it is denoted by K(G) or simply κ (Kappa) when G is understood. If G is trivial or disconnected K(G) is taken to be zero where as if G contains* $K_n$ *as a spanning subgraph, K(G) is taken to be n-1. A set of vertices or edges of a connected graph G is said to disconnect the graph if its deletion results in a disconnected graph.*

When a connected graph G (on n vertices) does not contain $K_n$ as a spanning subgraph, K is the connectivity of G if there exists a set of K vertices of G whose deletion results in a disconnected subgraph of G while no set of (K – 1) or fewer vertices has this property.

Several results can be had from [3, 35]. Interested readers are requested to refer [23].

Now we proceed on to recall the notion of edge connectivity.

**DEFINITION 2.3.6:** *The edge connectivity of a connected graph G is the smallest K for which there exists a K-edge cut. The edge connectivity of a trivial or disconnected graph is taken to be 0. The edge connectivity of G is denoted by* $\lambda$ *(G). If* $\lambda$ *is the edge connectivity of a connected graph G, there exists a set of* $\lambda$ *edges whose deletion results in a disconnected graph and no subset of edges of G of size less than* $\lambda$ *has this property. A graph is r-connected if* $K(G) \geq r$*, G is r-edge connected if* $\lambda (G) \geq r$.

***Result:*** For any loop less connected graph G, $K(G) \leq \lambda (G) \leq \delta(G)$.

**DEFINITION 2.3.7:** *A family of two or more paths in a graph G is said to be internally disjoint if no vertex of G is an internal vertex of more than one path in the family.*

Several interesting results and examples in this direction can be had from any book on graph theory. [23]

**DEFINITION 2.3.8:** *A graph G is non separable if it is non trivial connected and has no cut vertices. A block of the graph is the maximal non separable subgraph of G. If G has no cut vertex, G itself is a block.*

Consider the graph, which has 5 blocks.

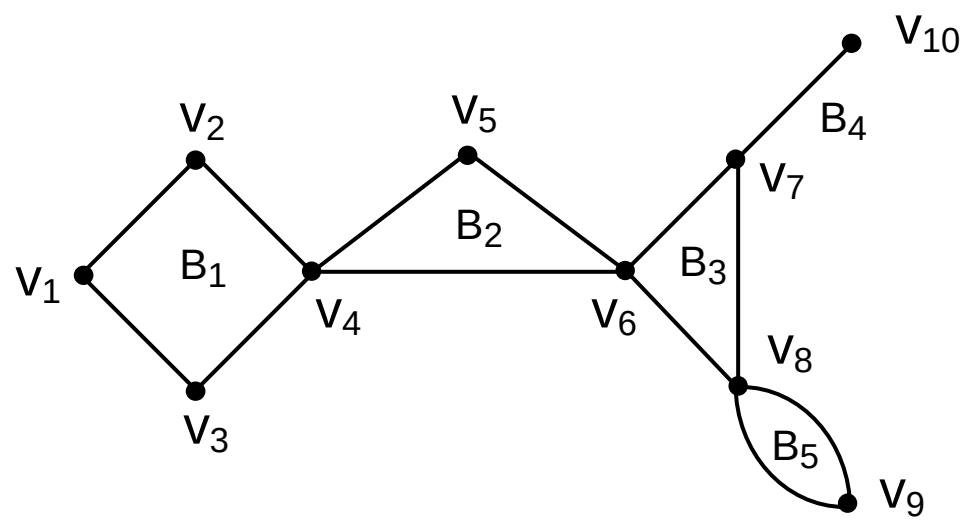

FIGURE: 2.3.7

The blocks of the graph are

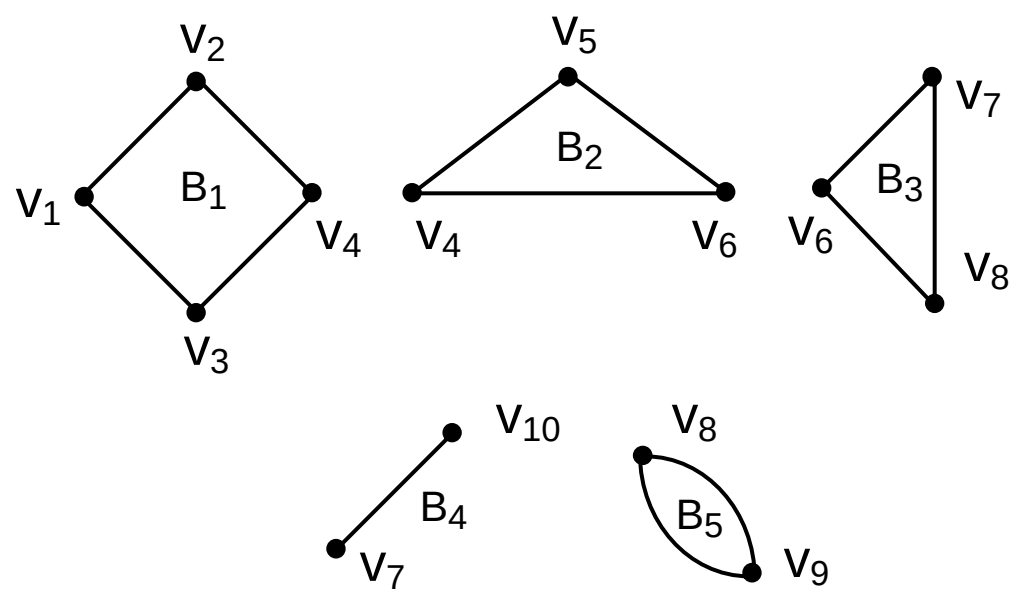

FIGURE: 2.3.8

**DEFINITION 2.3.9:** *Let G be a connected simple graph containing at least two disjoint cycles. Then the cyclical edge connectivity of G is defined to be the minimum number of edges of G whose deletion results in a graph having two components each containing a cycle. It is denoted by* $\lambda_C$ *(G).*

## 2.4 Trees

Graphs also derive their names from the diagrams. A tree is one such graph formally a connected graph without cycles is defined as a tree. A graph without cycles is called an acyclic graph or a forest. So each component of a forest is a tree. A forest may consists of just a single tree.

There is only one tree with one point, two point, three points

FIGURE: 2.4.1

two trees with four points

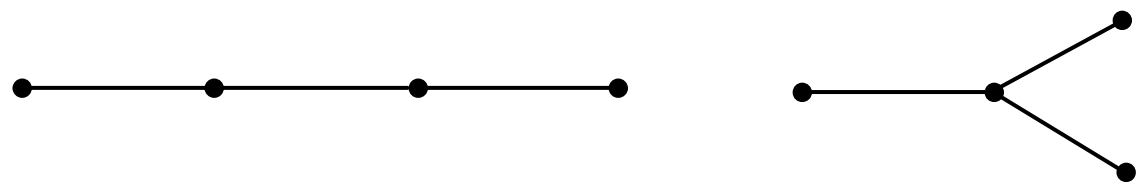

FIGURE: 2.4.2

There are 11 trees with 7 points

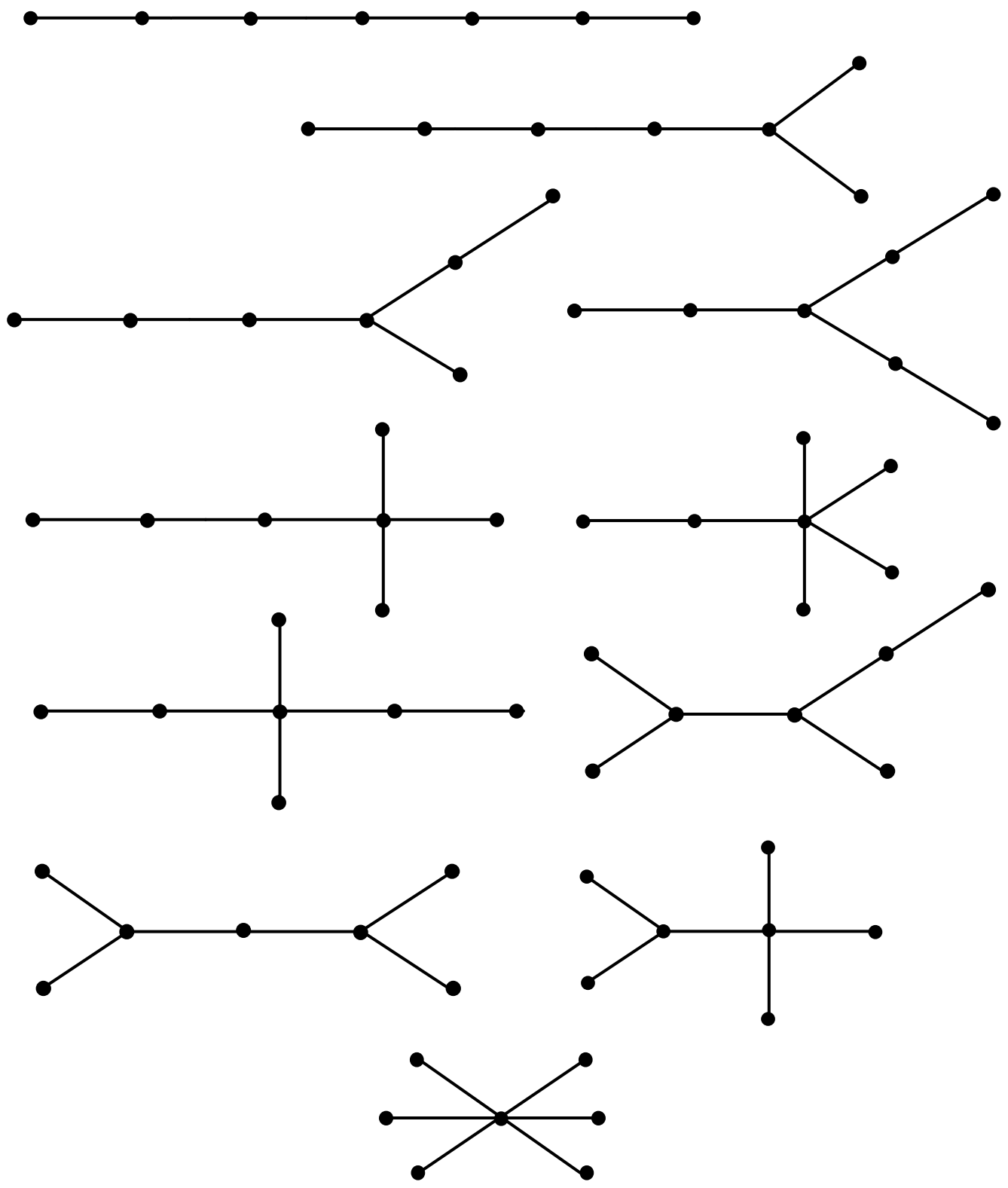

FIGURE: 2.4.3

There are twenty three trees with eight points. A branch at a point u of a tree T is a maximal subtree containing u as an end point.

A spanning subgraph of a graph which is also a tree is called a spanning tree of the graph. A graph G and two of its spanning trees $T_1$ and $T_2$ are shown in Figure 2.4.4:

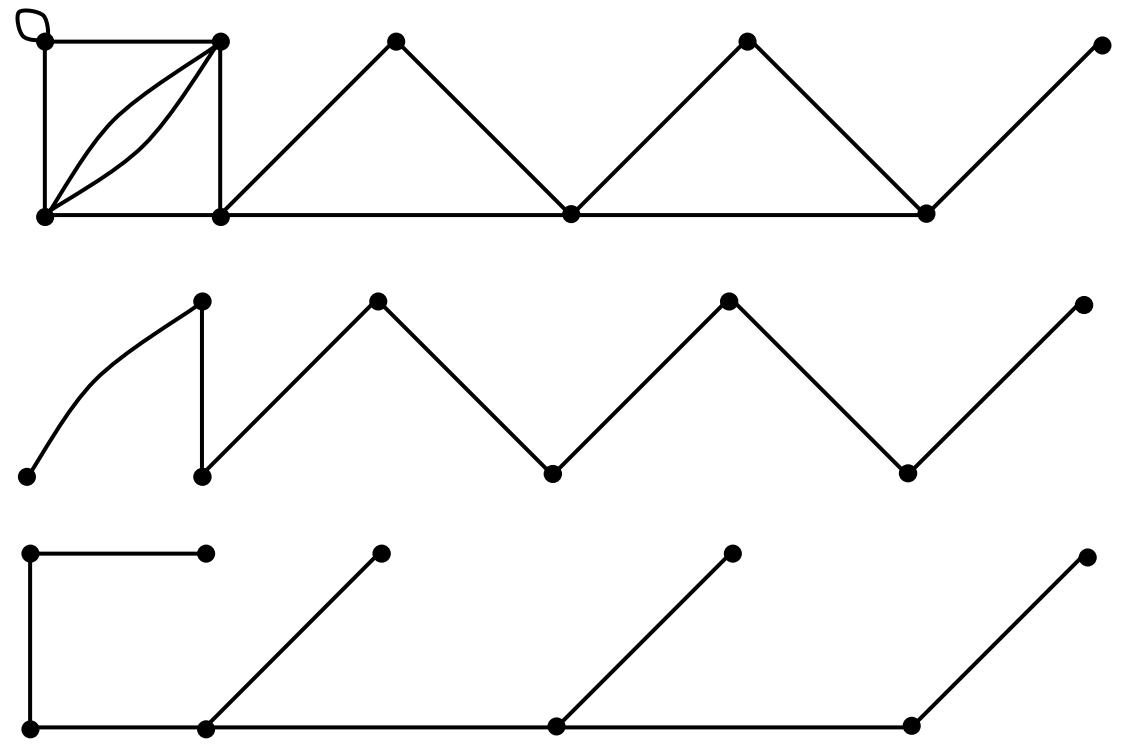

FIGURE: 2.4.4

Many interesting results in this direction can be had from [23]. A branch at a vertex u of a tree T is a maximal subtree containing u as an end vertex. Hence the number of branches at u is d (u). The weight of a vertex u of T is the maximum number of edges in any branch at u.

A vertex v is a controid vertex of T if v has minimum weight. The set of all controid vertices is called the controid of T. Counting the number of spanning trees in a graph occurs as a natural problem in many branches of science. Spanning trees were used by Kirchoff to generate a "cycle basis" for the cycles in the graphs of electrical networks. Now we consider enumeration of the spanning trees of graphs. The number of spanning trees of a connected labeled graph G will be denoted by $\tau$ (G). There is a recursive formula for $\tau$ (G). Before we establish this formula we shall define the concept of edge contraction in graphs.

**DEFINITION 2.4.1:** *An edge e of a graph G is said to be contracted, if it is deleted from G and its ends are identified. The resulting graph is denoted by G.e.*

If e is not a loop of G, then n(G, e) = n(G) – 1. m(G, e) = m(G) – 1 and w(G, e) = w(G). For a loop e, m(G, e) = m(G) – 1 and n(G, e) = n (G) and w (G, e) = w (G).

As our main aim is to study the graphs with indeterminate vertex and (or) indeterminate edges i.e more specifically on neutrosophic vertex graphs and neutrosophic edge graphs and doubly or strongly neutrosophic graphs we restrict ourselves only with a very few properties.

Several results can be had from [22, 35].

## 2.5 Independent – Sets and matchings

In this section all graphs considered are loopless. We recall some basic definitions of K-regular graphs and other related notions.

**DEFINITION 2.5.1:** *A subset S of the vertex set V of a graph G is called independent if no two vertices of S are adjacent in G.* $S \subset V$ *is a maximum independent set of G if G has no independent set S' with* $|S'| > |S|$.

A maximal independent set of G is an independent set that is not a proper subset of another independent set of G.

**DEFINITION 2.5.2:** *A subset K of V is called a covering of G if every edge of G is incident with atleast one vertex of K. A covering K is minimum if there is no covering K' of G such that* $|K'| < |K|$; *it is minimal if there is no covering* $K_1$ *of G such that* $K_1$ *is a proper subset of K.*

**DEFINITION 2.5.3:** *The number of vertices in a maximum independent set of G is called the independent number of G and is denoted by* $\alpha$ *(G). The number of vertices in a minimum covering of G is the covering number of G and is denoted by* $\beta$ *(G).*

Now we proceed on to recall the definition of edge independent set.

**DEFINITION 2.5.4:** *A subset M of the edge set E of a loopless graph G is called independent if no two edges of M are adjacent in G. A matching in G is a set of independent edges.*

An edge covering of G is a subset L of E such that every vertex of G is incident to some edge of L. Hence an edge covering of G exists if and only if $\delta > 0$. A matching M of G is maximum if G has no matching $M^1$ with $|M^1| > |M|$. M is maximal if G has no matching $M^1$ strictly containing M. $\alpha^1$ (G) is the cardinality of a maximum matching and $\beta^1$ (C) is the size of a minimum edge covering of G.

A set S of vertices of G is said to be a matching M of S if every vertex of S is incident to some edge of M: a vertex v of G is M-unsaturated if it is not M–saturated Herschel graph.

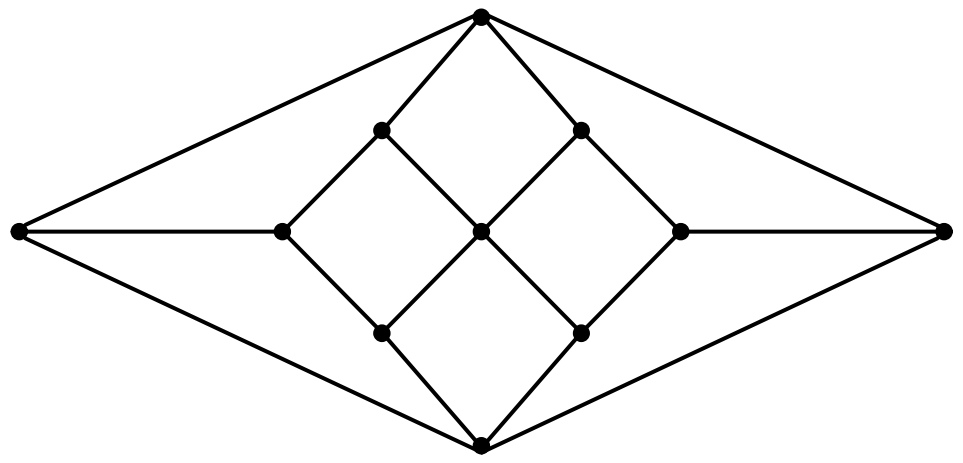

HERSCHEL GRAPH

FIGURE: 2.5.1

**DEFINITION 2.5.5:** *A matching of a graph G is a set of independent edges of G. If e = uv is an edge of a matching M of G the end vertices u and v of e are said to be matched by M.*

If $M_1$ and $M_2$ are matching of G the edge subgraph defined by $M_1 \Delta M_2$ the symmetric differences of $M_1$ and $M_2$ is a subgraph H of G whose components are paths or even cycles of G in which the edges alternate between $M_1$ and $M_2$.

**DEFINITION 2.5.6:** *An M-augmenting path in G is a path in which the edges alternate between E \ M and M and its end vertices are M-unsaturated. An M-alternating path in G is a path whose edges alternate between E \ M and M.*

Now we proceed on to give the definition of factor of a graph for more about these concepts please refer [23].

**DEFINITION 2.5.7:** *A factor of a graph G is a spanning subgraph of G. A K-factor of G is a factor of G that is K-regular. Thus a 1-factor of G is a matching that saturates all the vertices of G. For this reason a 1-factor of G is called a perfect matching of G. A two factor of G is a factor of G that is disjoint union of cycles of G. A graph G is K-factorable if G is an edge disjoint union of K-factors of G.*

Several nice and important results can be had from any text book on Graph Theory.

Let A = ($a_{ij}$) be a binary matrix of size p by q.

From a bipartite graph G with bipartition (X, Y) where X and Y are sets of cardinality p and q respectively say X = {$\upsilon_1, \upsilon_2$, …, $\upsilon_p$} and Y = {$\omega_1$, $\omega_2$, …, $\omega_q$}. Make $\upsilon_I$ adjacent to $w_j$ in G if and only if $a_{ij}$ = 2. Then an entry 1 in A corresponds to an edge of G and two independent 1's correspond to two independent edges of G. Further each vertex of G corresponds to a line of A. Thus the matrix version of Konigs theorem is actually a restatement of Konig's theorem. Existence of a System of Distinct Representatives (SDR) for a family of subsets of a given finite set.

**DEFINITION 2.5.8:** *Let F = $\{A_\alpha / \alpha \in J\}$ be a family of sets. A SDR for the family F is a family of elements $\{x_\alpha / \alpha \in J\}$ such that $x_\alpha \in A_\alpha$ for every $\alpha \in J$ and $x_\alpha \neq x_\beta$ whenever $\alpha \neq \beta$.*

Halls theorem on the existence of an SDR is given without proof.

**THEOREM 2.5.1:** *(Hall's Theorem on the existence of an SDR). [58] Let $\tau = \{A_i / 1 \leq i \leq r\}$ be a family of finite sets. Then $\tau$ is an SDR if and only if the union of any K, $1 \leq K \leq r$, members of $\tau$ contains at least K-elements.*

Several results exploiting these properties can be had from [3, 25, 58].

Just for the sake of completeness we recall the notions of perfect matchings and the Tutte matrix.

**DEFINITION 2.5.9:** *Let G = (V, E) be a simple graph of order n and let V = $\{v_1, v_2, \ldots, v_n\}$. Let $\{x_{ij} / 1 \leq i < j \leq n\}$ be a set of indeterminates. Then the Tutte matrix of G is defined to be the n by n matrix T = ($t_{ij}$)*
*where*

$$t_{ij} = \begin{cases} x_{ij} & \text{if } v_i v_j \in E(G),\ i < j \\ -x_{ij} & \text{if } v_i v_j \in E(G),\ i > j \\ 0 & \text{otherwise} \end{cases}$$

*Thus T is a skew symmetric matrix of order n.*

***Example 2.5.1:*** Let G be the graph given by the following figure

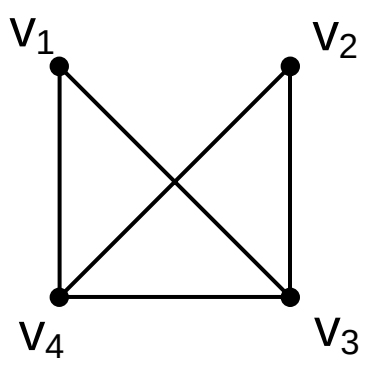

FIGURE: 2.5.2

The Tutte matrix of the graph G is

$$T = \begin{bmatrix} t_{11} & t_{12} & t_{13} & t_{14} \\ t_{21} & t_{22} & t_{23} & t_{24} \\ t_{31} & t_{32} & t_{33} & t_{34} \\ t_{41} & t_{42} & t_{43} & t_{44} \end{bmatrix}$$

$$= \begin{bmatrix} 0 & 0 & x_{13} & x_{14} \\ 0 & 0 & x_{23} & x_{24} \\ -x_{13} & -x_{23} & 0 & x_{34} \\ -x_{14} & -x_{24} & -x_{34} & 0 \end{bmatrix}$$

We just state Tutte theorem the proof of which is left for the reader as an exercise.

**THEOREM [3, 23]:** *Let G be a simple graph with Tutte matrix T. Then G has a 1-factor if and only if det T $\neq$ 0.*

The following graph gives Tutte matrix for which $G_1$ has a 1-factor.

The graph of G is as follows:

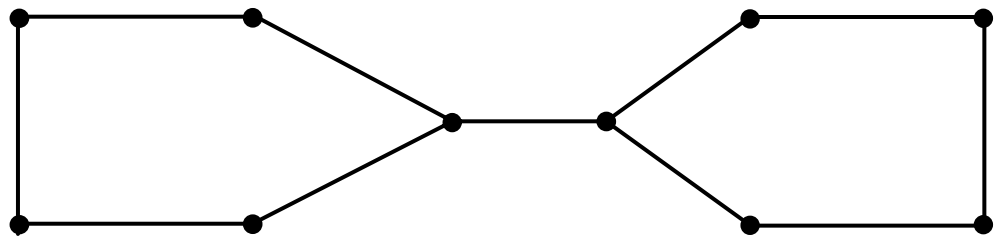

$G_1$

FIGURE: 2.5.3

For more results refer [21, 41].

## 2.6 Eulerian and Hamiltonian Graphs

The study of Eulerian graphs started in the $18^{th}$ century and that of Hamiltonian graphs in the $19^{th}$ century. These graphs have rich structures and hence their study is a very fertile field of research for graph theorists.

**DEFINITION 2.6.1:** *An Euler trail in a graph G is a spanning trail in G that contains all the edges of G. An Euler tour of G is a closed Euler trail of G. G is called Eulerian in G if G hus an Euler tour.*

It was Euler who first considered these graphs and hence their name.

***Example 2.6.1:*** The following graph G is an Eulerian graph.

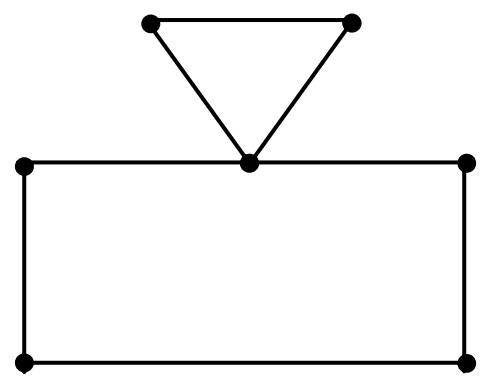

FIGURE: 2.6.1

***Example 2.6.2:*** This gives a Eulerian graph with 12 vertices.

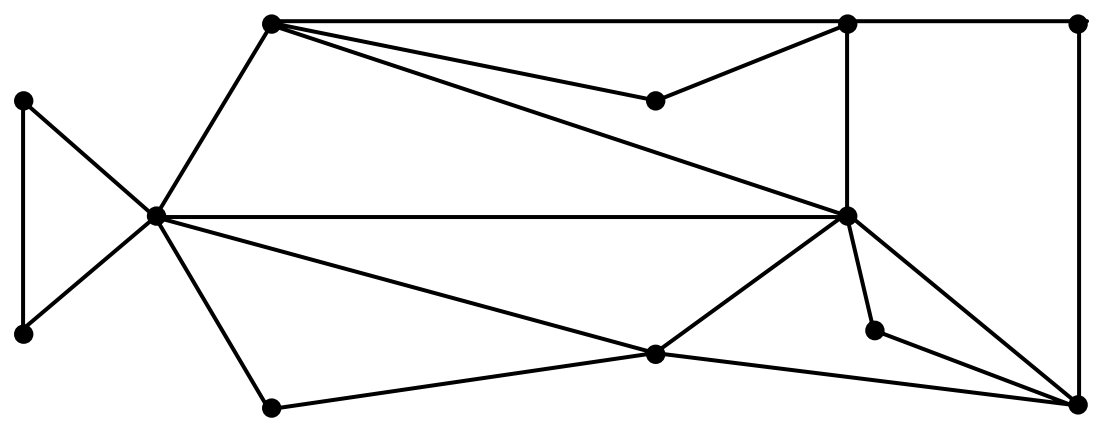

FIGURE: 2.6.2

Euler was the one to show in the year 1736 that the famous Konigsberg bridge problem has no solution. The following theorem which gives equivalence of three conditions is left for the reader as an exercise to prove.

**THEOREM 2.6.1:** *For a connected graph G the following statements are equivalent*

*i.* *E is Eulerian*
*ii.* *The degree of each vertex of G is an even positive integer*
*iii.* *G is an edge disjoint union of cycles*

The reader is expected to prove the following theorems.

**THEOREM 2.6.2:** *A connected graph is Eulerian if and only if it admits a cycle decomposition.*

The following example gives a Eulerian graph with edge e belonging to three cycles.

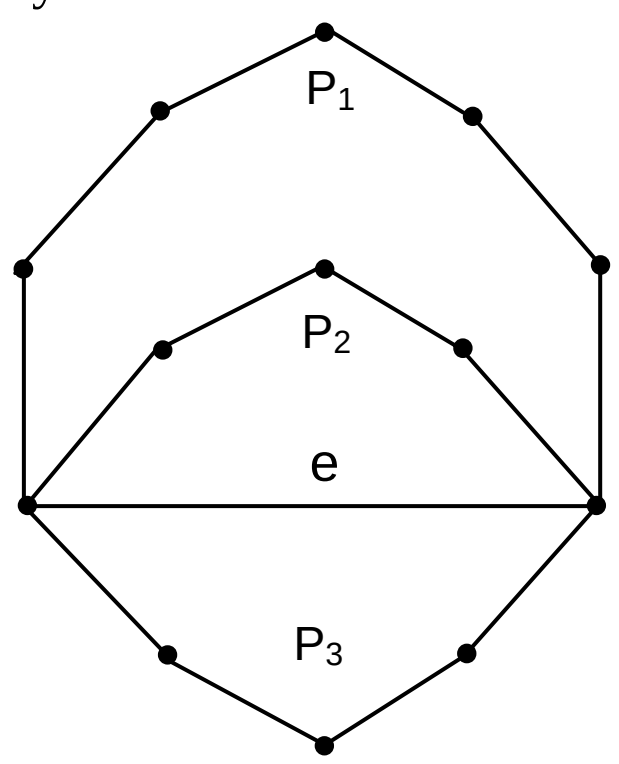

FIGURE: 2.6.3

Finally the author is expected to prove.

**THEOREM 2.6.3:** *A graph is Eulerian if and only if each edge e of G belongs to an odd number of cycles of G.*

**THEOREM 2.6.4:** *A graph is Eulerian if and only if it has an odd number of cycle decompositions.*

Now we just recall the definition of Eulerian digraph.

An Eulerian trial in a digraph D is a closed spanning walk in which each arc of D occurs exactly once. A digraph is eulerian if it has such a trial. Just we have a weak digraph D, is eulerian if and only if every point of D has equal indegree and outdegree.

**THEOREM 2.6.5:** *In an Eulerian digraph the number of Eulerian trails is*

$$C \prod_{i=1}^{p} (d_i - 1)\,!$$

*where $d_i$ = id ($v_i$) and C is the common value of all the cofactors.*

Several interesting properties in this direction can be had from [23].

Now we proceed on to recall the definition of Hamiltonian graphs.

**DEFINITION 2.6.2:** *A graph is called Hamiltonian if it has a spanning cycle.*

These graphs were first studied by Sir William Hamilton a mathematician. A spanning cycle of a graph G when it exists is often called a Hamilton cycle or a Hamiltonian cycle of G.

**DEFINITION 2.6.3:** *A graph G is called traceable if it has a spanning path of G. A spanning path of G is also called a Hamilton path of G.*

Hamilton introduced these graphs in 1859. The following graph is a Hamiltonian graph.

***Example 2.6.3:***

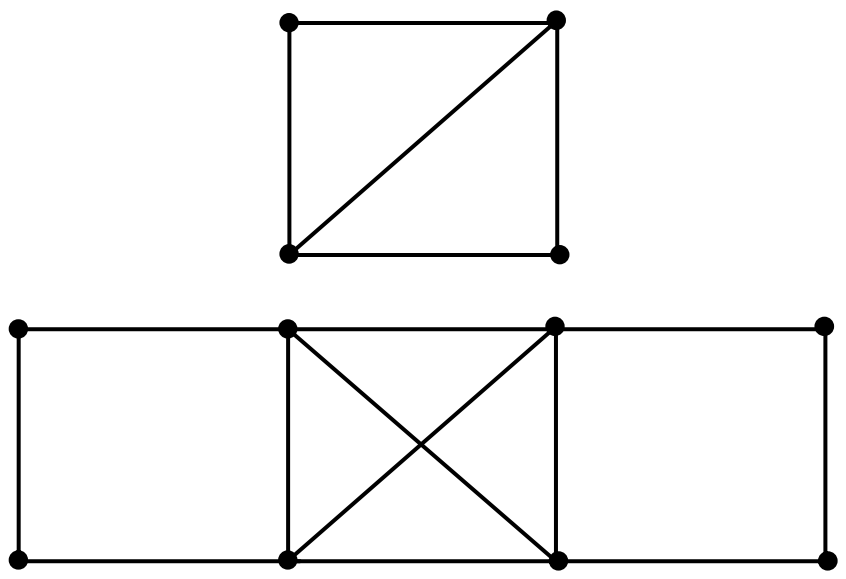

FIGURE: 2.6.4

The first figure is an Hamiltonian graph with 4 vertices and the next figure is an Hamiltonian graph with 8 vertices. The following is an example of a traceable graph which is non Hamiltonian.

***Example 2.6.4:***

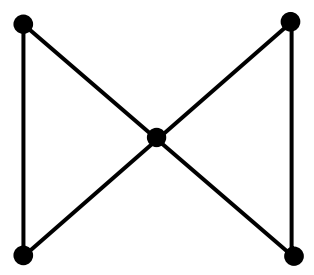

FIGURE: 2.6.5

a non Hamiltonian graph but which is a traceable graph.

Several interesting results in this direction are proved. We in this book just recall some of the results without proof.

***Result:*** If G is Hamiltonian then for every non empty proper subset S of V w (G – S ) $\leq |S|$.

***Result:*** Let G be a simple graph with n ( $\geq$ 3) vertices. If for every pair of non adjacent vertices u, v of G, d (u) + d (v) $\geq$ n, then G is Hamiltonian.

For proof please refer [47].

***Result:*** If G is a simple graph with n $\geq$ 3 and $\delta \geq n/2$ then G is Hamiltonian.

For proof please refer [47].

**THEOREM 2.6.4:** *Let G be a simple graph with n ( $\geq$ 3) vertices. If d(u) + d (v) $\geq$ n – 1 for every pair of non adjacent vertices u and v of G then G is traceable.*

***Example 2.6.5:*** If G = G (X, Y) is a bipartite Hamiltonian graph then we have |X| = |Y|.

**THEOREM 2.6.7:** *Let G be a simple graph of order n ( $\geq$ 3) vertices. Then G is Hamiltonian if and only if G + uv is Hamiltonian for every pair of non adjacent vertices u and v with d (u) + d (v) $\geq$ n.*

Now we proceed on to recall a nice and interesting property of a graph viz closure.

**DEFINITION 2.6.4:** *The closure of a graph G denoted by cl (G) is defined to be that super graph of G obtained from G by recursively joining pairs of nonadjacent vertices whose degrees sum is at least n until no such pair exists.*

The following example gives the closure of a graph.

***Example 2.6.6:***

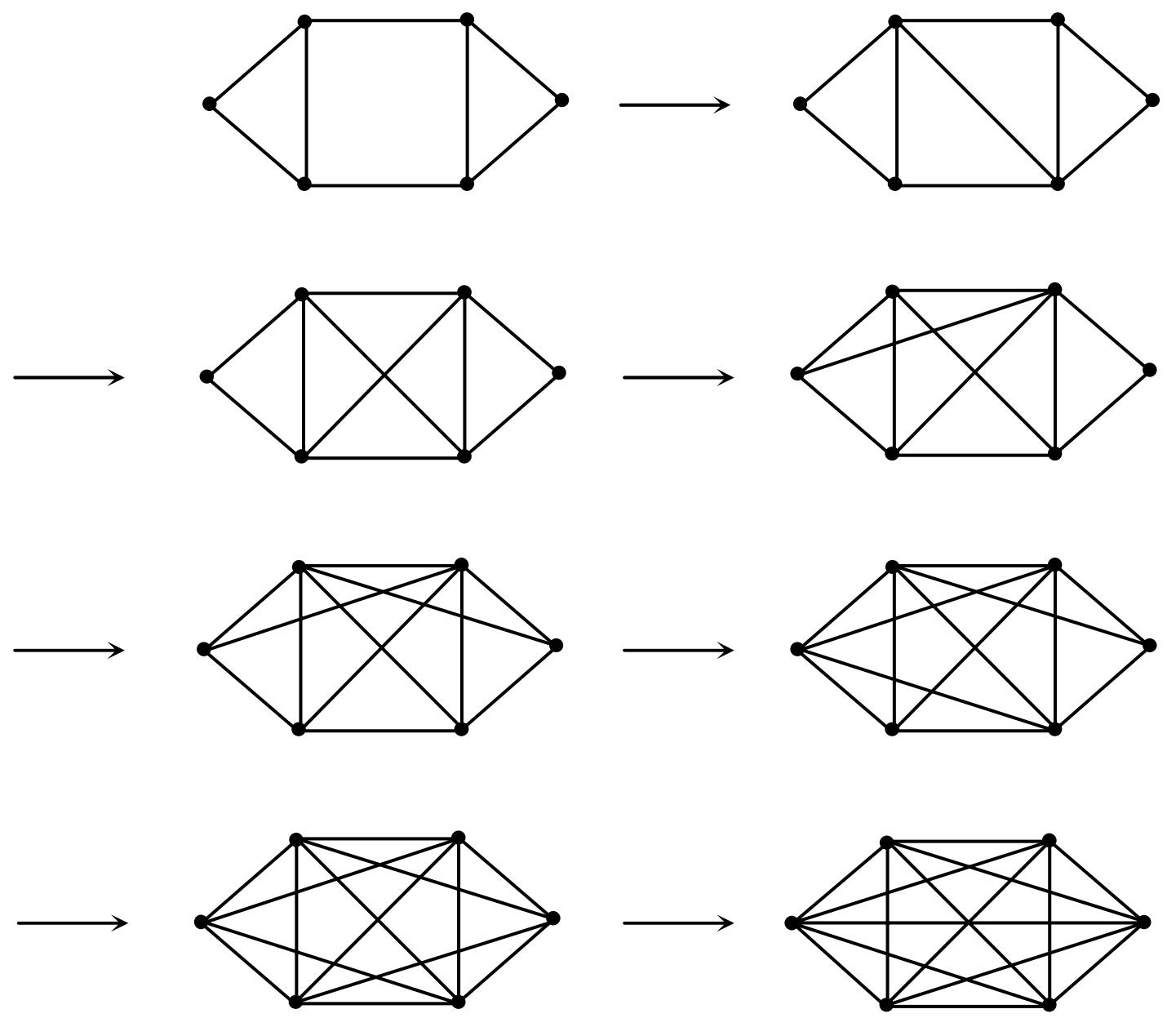

FIGURE: 2.6.6

The reader is expected to prove the following theorems:

**THEOREM 2.6.8:** *The closure cl (G) of a graph G is well defined.*

**THEOREM 2.6.9:** *If cl (G) is Hamiltonian then G is Hamiltonian.*

**THEOREM 2.6.10:** *If cl (G) is complete then G is Hamiltonian.*

The theorem on simple 2-connected graph given by [12] is an important and an interesting one.

**THEOREM [12]:** *If for a simple 2-connected graph G, $\alpha \leq K$ then G is Hamiltonian.*

For Proof please refer [12].

Finally we state yet another theorem for the reader to prove.

**THEOREM 2.6.11:** *If G is a simple graph with $n \geq 3$ vertices such that $d(u) + d(v) \geq n + 1$ for every pair of n on adjacent vertices of G, then G is Hamiltonian connected.*

Next we recall the definition of Pancyclic graph and their relation with the Hamiltonian graph.

**DEFINITION 2.6.5:** *A graph G of order n ($\geq 3$) is pancyclic if G contains cycles of all lengths from 3 to n. G is called vertex pancyclic if each vertex $\upsilon$ of G belongs to a cyclic of every length l, $3 \leq l \leq n$.*

The following graph is an example of a pancyclic graph which is not vertex pancyclic.

***Example:*** Pancyclic graph is given in figure 2.6.7.

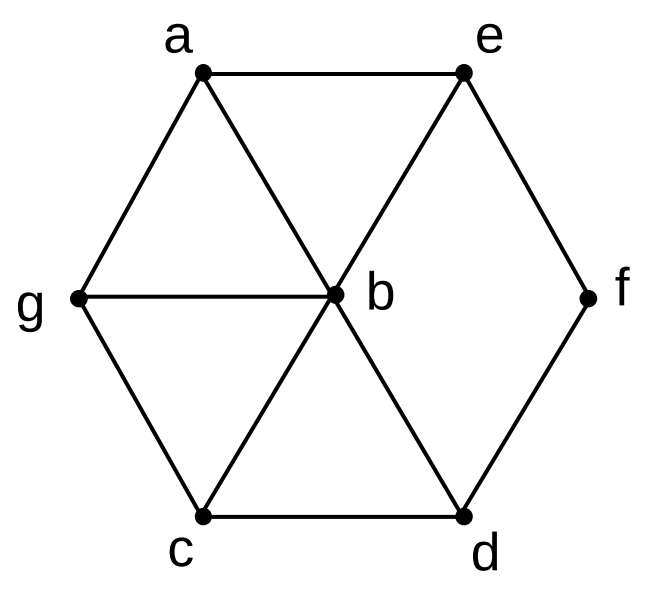

FIGURE: 2.6.7

**THEOREM 2.6.12:** *Let G be a simple Hamiltonian graph on n vertices with at least $\left[\frac{n^2}{2}\right]$ edges. Then G is either pancyclic or* *else is the complete bipartite graph* $K_{n/2, n/2}$. *In particular if G is Hamiltonian and*

$$m > \frac{n^2}{4},$$

*then G is pancyclic.*

The reader is expected to obtain the proof. For proof refer [12, 22-23].

**THEOREM 2.6.13:** *Let G* $\neq K_{B/2, n/2}$ *be a simple graph with n (≥ 3) vertices and let d (u) + d (v) ≥ n for every pair of nonadjacent vertices of G. Then G is pancyclic.*

*Proof:* Refer [22-3] .

Here we recall the condition for a Eulerian graph line graph to be Hamiltonian and Eulerian.

**THEOREM 2.6.14:** *If G is Eulerian then L (G), the line graph of G is both Hamiltonian and Eulerian.*

*Proof:* Refer [22-3].

**DEFINITION 2.6.6:** *A dominating trail of a graph G is a closed trail in G such that every edge of G not in T is incident with T.*

[22]characterized graphs that have Hamiltonian line graph.

**THEOREM [22]:** *The line graph of a graph G with at least three edges is Hamiltonian if and only if G has a dominating trial.*

*Proof:* Can be obtained from [22].

**THEOREM 2.6.15:** *The line graph of Hamiltonian graph is Hamiltonian.*

*Proof:* Let G be a Hamiltonian graph with Hamilton cycle C. Then C is a dominating trial of G. Hence L (G) is Hamiltonian.

**THEOREM [4]:** *Let G be any connected graph. If each edge of G belongs to a triangle in G then G has a spanning Eulerian subgraph.*

For proof refer [4].

The following results can be derived using the earlier results.

***Result:*** Let G be any connected graph. If each edge of G belongs to a triangle, then L(G) is Hamiltonian.

***Result:*** If G is connected and $\delta(G) \geq 3$ then $L^2$ (G) is Hamiltonian.

For proof please refer [11].

***Result:*** If G is a connected graph with at least three vertices then L ($G^2$) is Hamiltonian. For proof refer [11].

**THEOREM 2.6.16:** *Let G be a connected graph in which every edge belongs to a triangle. If $e_1$ and $e_2$ are edges of G such that G \ {$e_1$, $e_2$} is connected then there exists a spanning trail of G with $e_1$ and $e_2$ as its initial and terminal edges.*

For proof refer [11].

The following results are easy consequences of the above theorem.

***Result: [4]*:** Let G be any connected graph with $\delta$ (G) $\geq$ 4. Then $L^2$ (G) is Hamiltonian connected .

***Result [26]:*** The line graph of a 4-edge connected graph is Hamiltonian.

The theorem of [26] can be proved using the following lemma.

**LEMMA [26]:** *Let S be a set of vertices of a nontrivial tree T and let $|S| = 2K$, $K \geq 2$. Then there exists a set of K pairwise edge disjoint paths whose end vertices are all the vertices of S.*

Now we proceed on to recall the theorem on locally connected graphs.

**THEOREM [46]**: *A connected locally connected non trivial $K_{1,\,3}$ free graph is Hamiltonian.*

Proof is lengthy and the reader is expected to get from [46].

## 2.7 Graph Colorings

In this section we recall some of the basic results about graph colorings. We do not give any proof of the result we only recall the definition and state the results. It is up to the reader to prove the results or refer to get the results.

**DEFINITION 2.7.1:** *The chromatic number $\chi$ (G) of a graph G is the minimum number of independent subsets that partition the vertex set of G. Any such minimum partition is called the chromatic partition of V (G).*

A vertex coloring of G is a map f : V $\rightarrow$ S where S is the set of distinct colours it is proper if adjacent vertices of G receive distinct colours of S; that is if u v $\in$ E (G) then f (u) $\neq$ f (v). Thus $\chi$(G) is the minimum cardinality of S for which there exists a proper vertex coloring of G by colors of S. Clearly in any proper vertex coloring of G, the vertices that receive the same color are independent. The vertices that receive a particular color make up a color class. Thus in any chromatic partition of V (G), the parts of the partition constitute the color classes. This leads to an equivalent way of defining the chromatic number.

**DEFINITION 2.7.2:** *The chromatic number of a graph G is the minimum number of colors needed for a proper vertex coloring of G. G is K-chromatic if $\chi$(G) = K.*

**DEFINITION 2.7.3:** *A K-coloring of a graph G is a vertex coloring of G that uses K colors.*

**DEFINITION 2.7.4:** *A graph G is said to be K-colorable if G admits a proper vertex-coloring using K-colors.*

Clearly $\chi$ ($K_n$) = n. Further $\chi$(G) = 2 if and only if G is bipartite having at least one edge. In particular $\chi$(T) = 2 for any tree T with at least one edge (since any tree is bipartite).

Several interesting results about different types of graphs can be had from [2, 23, 27, 47].

**DEFINITION 2.7.5:** *A graph G is called critical if for every proper subgroup H of G; $\chi$(H) < $\chi$(G). Also G is K-critical if it is K-chromatic and critical.*

The reader is expected to prove the following:

i. Prove any critical graph is connected.
ii. Show that a graph is 3-critical if and only if it is an odd-cycle.
iii. If G is K – critical then $\delta \geq K - 2$.
iv. For any graph $\chi \leq 1 + \Delta$ .
v. Prove or disprove if G is K-chromatic then G contains $K_k$.
vi. In a critical graph G no vertex cut is a clique.
vii. Every critical graph is a block.
viii. If a connected graph G is neither an odd cycle nor a complete graph then $\chi(G) \leq \Delta(G)$.

Several results can be had in this direction from any text on graph theory. For proof of these results the reader if need be refer [5, 23].

**DEFINITION 2.7.6:** *The achromatic number $\phi$ (G) of a graph G is the maximum K for which G has a complete K – coloring.*

**DEFINITION 2.7.7:** *A graph G is triangle free if G contains $K_3$.*

**THEOREM [23]:** *For every positive integer K, there exists a triangle free graph with chromatic number K.*

For proof please refer [44]

Now we proceed on to analyze the edge colouring of the graphs.

**DEFINITION 2.7.8:** *An edge coloring of a loop less graph G is a function $\pi$: E (G) $\to$ S where S is a set of distinct colors: it is proper if no two adjacent edges receive the same color. Thus a proper edge coloring $\pi$ of G is a function $\pi$: E (G) $\to$ S such that $\pi(e) \neq \pi(e')$ whenever the edges e and e' are adjacent in G.*

**DEFINITION 2.7.9:** *The minimum K for which a loopless graph G has a proper K-edge coloring is called the edge chromatic number or chromatic index of G. It is denoted by $\chi'$ (G). G is K-edge chromatic number if $\chi'$ (G) = K.*

Further if an edge uv is colored by color C we say that C is represented at both u and v. If G has a proper K-edge coloring, E(G) is partitioned into K-edge disjoint matchings.

It is clear in a loopless graph $\chi'(G) \geq \Delta(G)$ since $\Delta(G)$ edge incident at a vertex u of maximum degree $\Delta(G)$ must all receive distinct colors. However for bipartite graphs equality holds.

The proof of the following theorems are left as a exercise.

**THEOREM [35]:** *If G is a loopless bipartite graph, $\chi'(G) = \Delta(G)$.*

For proof please refer [35].

**THEOREM 2.7.1:** $\chi'(K_n) = \begin{cases} n-1 & \text{if } n \text{ is even} \\ n & \text{if } n \text{ is odd} \end{cases}$

Prove the following results.

1. Show that a Hamiltonian cubic graph is 3-edge chromatic
2. Show that the Petersen graph is 4-edge chromatic
3. Describe a proper K edge coloring of a K-regular bipartite graph.

**THEOREM 2.7.2:** *For any simple graph G, $\Delta(G) \leq \chi'(G) \leq 1 + \Delta(G)$.*

For proof please refer.

A Graphs for which $\chi' = \Delta$ are called class I graphs and those for which $\chi' = 1 + \Delta$ are called 2 graphs. In view of the above theorem we have the following. If G is a cubic simple graph then $\chi'(G) = 3$ or 4. If G is a cubic simple graph with a cut-edge then $\chi'(G) = 4$.

**DEFINITION 2.7.10:** *A snark is a cyclically 4-edge connected cubic graph of girth at least 5 and has chromatic index 4.*

It is left for the reader to verify that no snark can be Hamiltonian.

**THEOREM 2.7.3:** *The Petersen graph P is the smallest snark and it is the unique snark with ten vertices.*

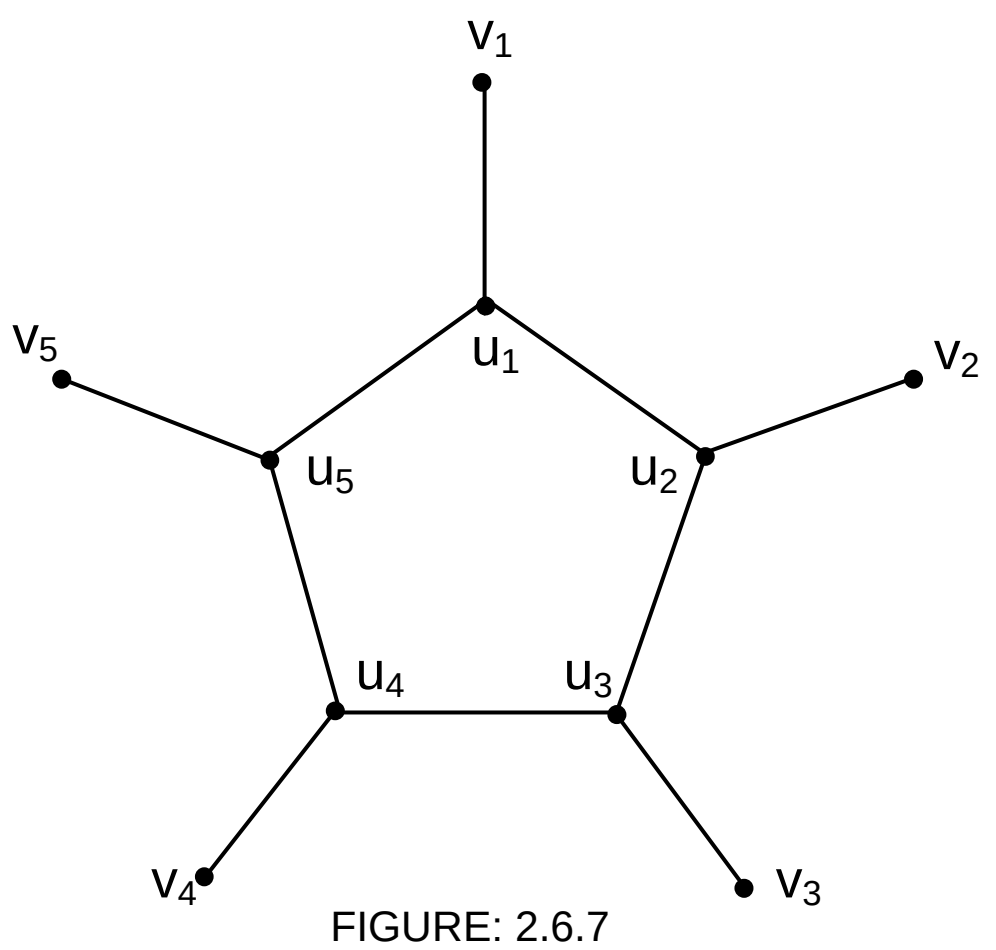

FIGURE: 2.6.7

For proof refer [23]. Several interesting results and applications can be had from any book on Graph Theory [5, 23, 25].

Let G be any graph, given a set of $\lambda$ colors. The function f(G: $\lambda$) is defined to be the number of ways of coloring G properly using $\lambda$ colors. Hence f (G; $\lambda$) = 0, when G has no proper $\lambda$-coloring. Clearly the minimum $\lambda$ for which f (G; $\lambda$) > 0 is the chromatic number $\chi$(G) of G. Clearly f ($K_n$, $\lambda$) = $\lambda$ ($\lambda$ – 1) … ($\lambda$ – n + 1) for $\lambda \geq n$.

$$f(K_3; \lambda) = \lambda(\lambda - 1)(\lambda - 2)$$
$$f(K_n^C; \lambda) = \lambda^n.$$

**THEOREM 2.7.4:** *If G is any graph. Then f (G; λ) = f (G – e, λ) – f(G e, λ) for any edge e of G.*

The proof is simple and the reader is advised to refer [3, 23, 25].

It is easily proved if G and H are adjoint then f (G $\cup$ H, $\lambda$) = f(G; $\lambda$); f (H; $\lambda$) f (G; $\lambda$) is called the chromatic polynomial of the graph G.

**THEOREM 2.7.5:** *For a simple graph G of order n and size m. f(G; λ) is a monic polynomial of degree n in λ with integer coefficients and constant term zero. In addition its coefficients alternate in sign and the coefficient of* $\lambda^{n-1}$ *is – m.*

For proof refer [23].

The following results are given for the reader to prove.

***Result:*** A simple graph G on n vertices in a tree if and only if f(G; $\lambda$) = $\lambda$ $(\lambda - 1)^{n-1}$.

***Result:*** If G has w components i then show that $\lambda^w$ is a factor of f(G ; $\lambda$).

For more about graph colouring. Refer Jensen and Toft [27] Fiorini and Wilson [17].

## 2.8 Application of Graphs to Fuzzy Models

Graphs have been the basis for several of the fuzzy models like binary fuzzy relations, sagittal diagrams, fuzzy compatibility relations fuzzy partial ordering relations, fuzzy morphisms fuzzy cognitive maps and fuzzy relational map models. Most of these models also basically rely on the under lying matrix. Both square and rectangular matrices are used. For more about these please refer [32-4, 64-6].

Just for the sake of completeness we at each stage illustrate each of these by giving a brief definition and by example.

**DEFINITION 2.8.1:** *Let X and Y be two finite sets contrary to functions from X to Y, binary relations R (X, Y) may assign to each element X two or more elements of Y some basic operations on functions such as composition or inverse may also be applicable to binary fuzzy relations.*

*The fuzzy relation R (X, Y) from the domain set X to the range set Y depicted by bi partite graph with edge weights.*

***Example 2.8.1:*** Let X = {$x_1$, $x_2$, $x_3$ $x_4$ $x_5$} and Y = {$y_1$, $y_2$, $y_3$, $y_4$, $y_5$, $y_6$} suppose R (XY) fuzzy membership relation.

The related bipartite graph with associated edge weights.

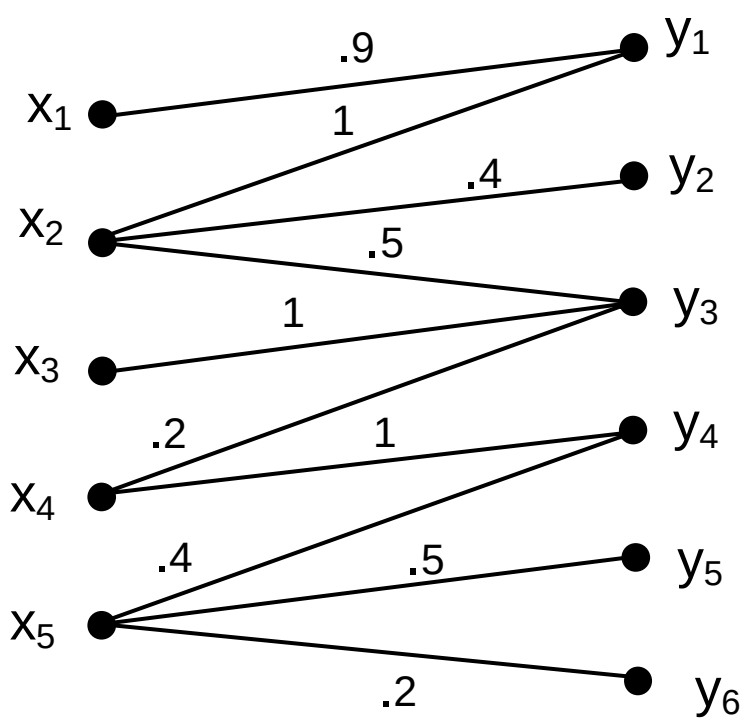

FIGURE: 2.8.1

The relational membership matrix

$$R(X, Y) = R = \begin{array}{c} \\ x_1 \\ x_2 \\ x_3 \\ x_4 \\ x_5 \end{array} \begin{array}{c} \begin{array}{cccccc} y_1 & y_2 & y_3 & y_4 & y_5 & y_6 \end{array} \\ \left[\begin{array}{cccccc} .9 & 0 & 0 & 0 & 0 & 0 \\ 1 & .4 & .5 & 0 & 0 & 0 \\ 0 & 0 & 1 & 0 & 0 & 0 \\ 0 & 0 & .2 & 1 & 0 & 0 \\ 0 & 0 & 0 & .4 & .5 & .2 \end{array}\right] \end{array}$$

The extension of two binary relation is a ternary relation denoted by R (X, Y, Z) where P (X, Y) and Q (Y, Z) are fuzzy relations from X to Y and Y to Z respectively.

Suppose X = $\{x_1\ x_2\ x_3\ x_4\}$ Y = $\{y_1\ y_2\ y_3\}$ and Z = $\{Z_1\ Z_2\}$ The representation of the composition of the two binary relations given by R (X, Y, Z) is given by the following graph

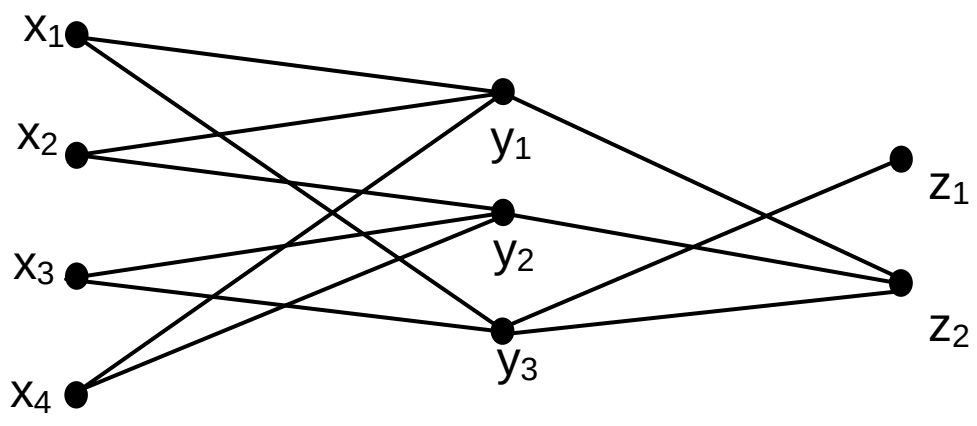

FIGURE: 2.8.2

with the desired edge weights

***Example 2.8.2***: The Binary relation on a single set. The related weighted graph or the sagittal diagram is follows. Let X = {$x_1$, $x_2$ $x_3$ $x_4$ $x_5$}.

Membership matrix

$$\begin{array}{c} \\ x_1 \\ x_2 \\ x_3 \\ x_4 \\ x_5 \end{array} \begin{array}{c} \begin{array}{ccccc} x_1 & x_2 & x_3 & x_4 & x_5 \end{array} \\ \begin{bmatrix} .8 & .3 & 0 & 0 & 0 \\ 0 & .7 & 0 & .8 & 0 \\ 0 & 0 & 0 & 0 & .3 \\ 0 & 0 & 0 & .4 & .3 \\ 0 & .1 & .6 & 0 & 0 \end{bmatrix} \end{array}$$

The related bipartite graph

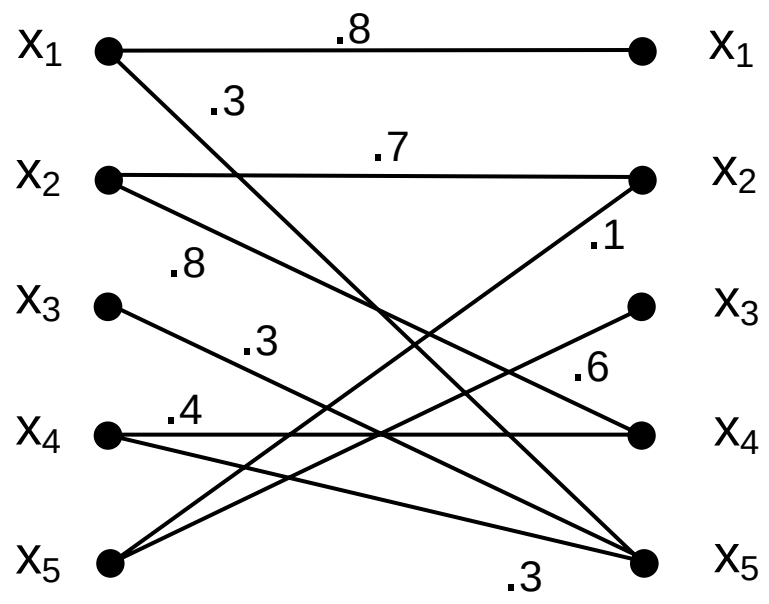

FIGURE: 2.8.3

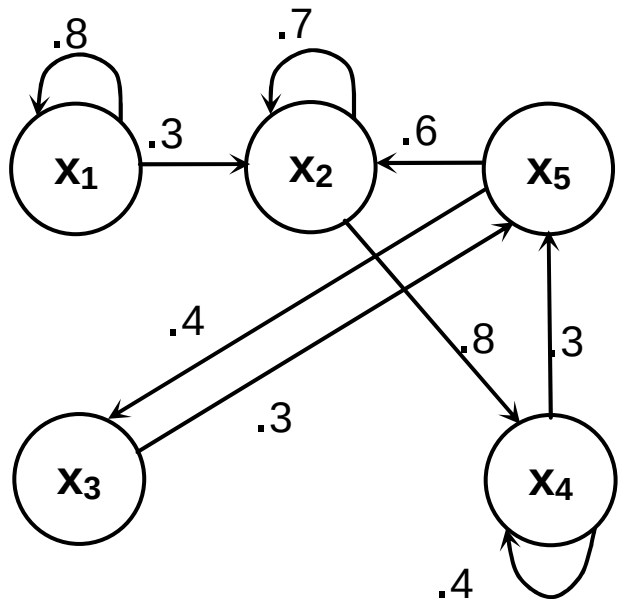

FIGURE: 2.8.4

The graph representation by a single diagram.

Study in this direction using this graph will be carried out for any appropriate models which is under investigation.

**DEFINITION 2.8.2:** *A binary relation R (X, X) that is reflexive and symmetric is usually called a compatibility relation or tolerance relation.*

***Example 2.8.3:*** Consider a fuzzy relation R (X, X) defined on X = {$x_1$, $x_2$, …, $x_7$} by the following membership matrix.

$$
\begin{array}{c} \\ x_1 \\ x_2 \\ x_3 \\ x_4 \\ x_5 \\ x_6 \\ x_7 \end{array}
\begin{array}{c}
\begin{array}{ccccccc} x_1 & x_2 & x_3 & x_4 & x_5 & x_6 & x_7 \end{array} \\
\left[\begin{array}{ccccccc}
1 & .3 & 0 & 0 & 0 & 0 & .2 \\
.3 & 1 & 0 & 0 & 0 & 0 & 0 \\
0 & 0 & 1 & 1 & 0 & 0.6 & 0 \\
0 & 0 & 1 & 1 & 0 & .1 & 0 \\
0 & 0 & 0 & 0 & 1 & .4 & 0 \\
0 & 0 & .6 & .1 & .4 & 1 & 0 \\
.2 & 0 & .0 & 0 & 0 & 0 & 1
\end{array}\right]
\end{array}
$$

Graph of the compatibility relation or the compatibility relation graph.

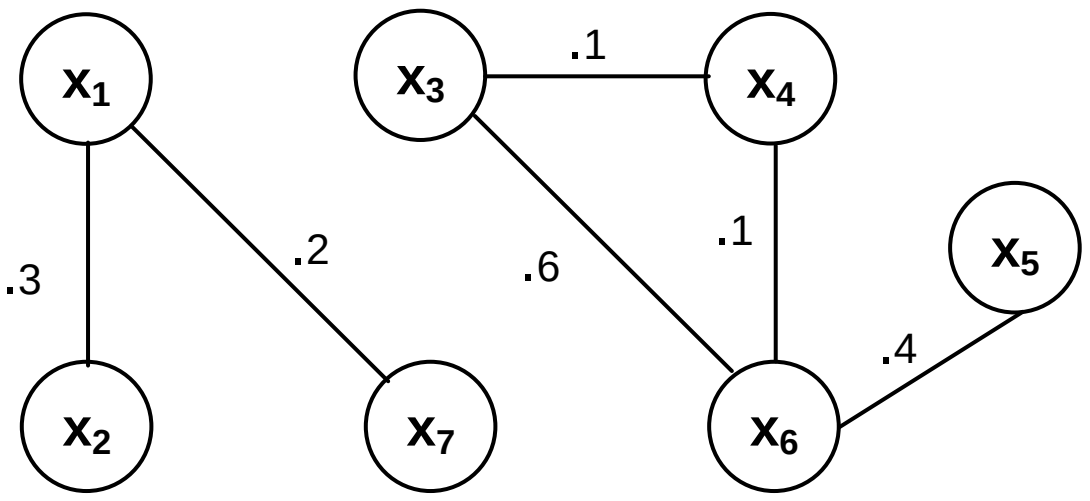

FIGURE: 2.8.5

Now we proceed on to define fuzzy partial ordering and its graphical representation.

**DEFINITION 2.8.3**: *A fuzzy binary relation R on a set X is a fuzzy partial ordering if and only if it is reflexive antisymmetric and transitive under some form of fuzzy transitivity.*

We represent this by the following example.

***Example 2.8.4:*** Let X = {a, b, c, d, e, f} P Q and Q denote the crisp partial ordering on the set X which are defined by their membership matrices and their graphical representation.

$$P = \begin{array}{c} \\ a \\ b \\ c \\ d \\ e \\ f \end{array} \begin{array}{c} \begin{array}{cccccc} a & b & c & d & e & f \end{array} \\ \begin{bmatrix} 1 & 0 & 0 & 0 & 0 & 0 \\ \underline{1} & 1 & 0 & 0 & 0 & 0 \\ 1 & \underline{1} & 1 & 0 & 0 & 0 \\ 1 & 1 & \underline{1} & 1 & 0 & 0 \\ 1 & 1 & 1 & \underline{1} & 1 & 0 \\ 1 & 1 & 1 & 1 & \underline{1} & 1 \end{bmatrix} \end{array}$$

The graph

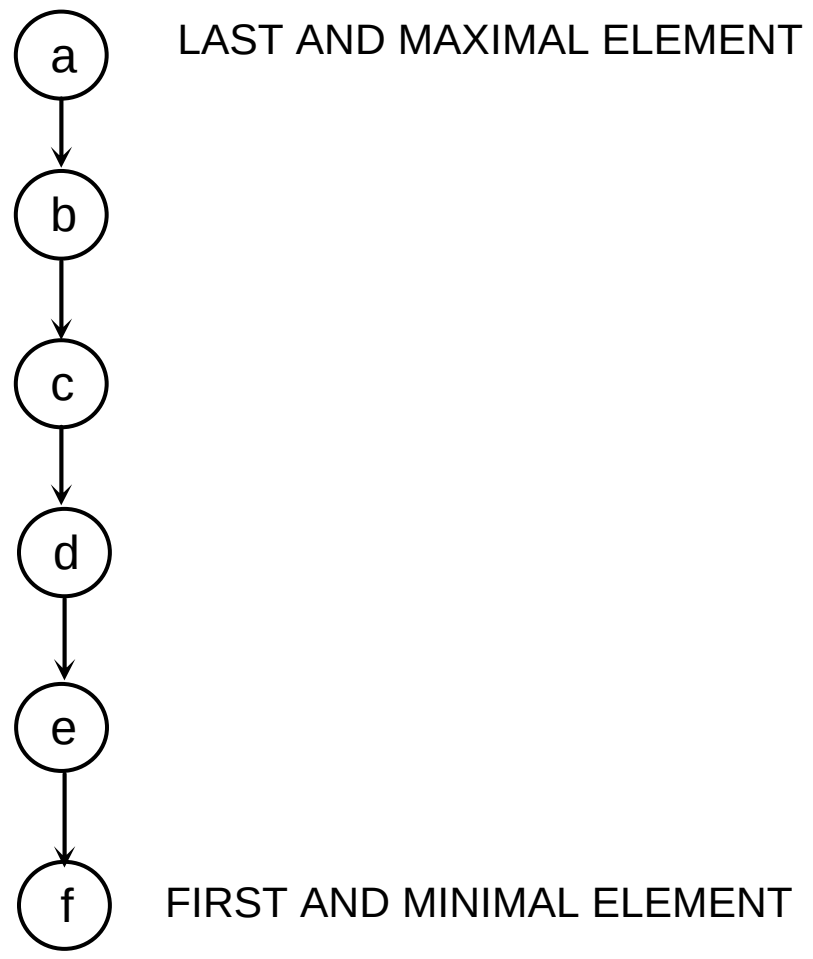

FIGURE: 2.8.6

The related matrix Q and its graphical representation

$$Q = \begin{array}{c} \\ a \\ b \\ c \\ d \\ r \\ f \end{array} \begin{array}{c} \begin{array}{cccccc} a & b & c & d & e & f \end{array} \\ \begin{bmatrix} 1 & 0 & 0 & 0 & 0 & 0 \\ \underline{1} & 1 & 0 & 0 & 0 & 0 \\ 1 & \underline{1} & 1 & 0 & 0 & 0 \\ 1 & \underline{1} & 0 & 1 & 0 & 0 \\ 1 & 1 & \underline{1} & \underline{1} & 1 & 0 \\ 1 & 1 & 1 & 1 & \underline{1} & 1 \end{bmatrix} \end{array}$$

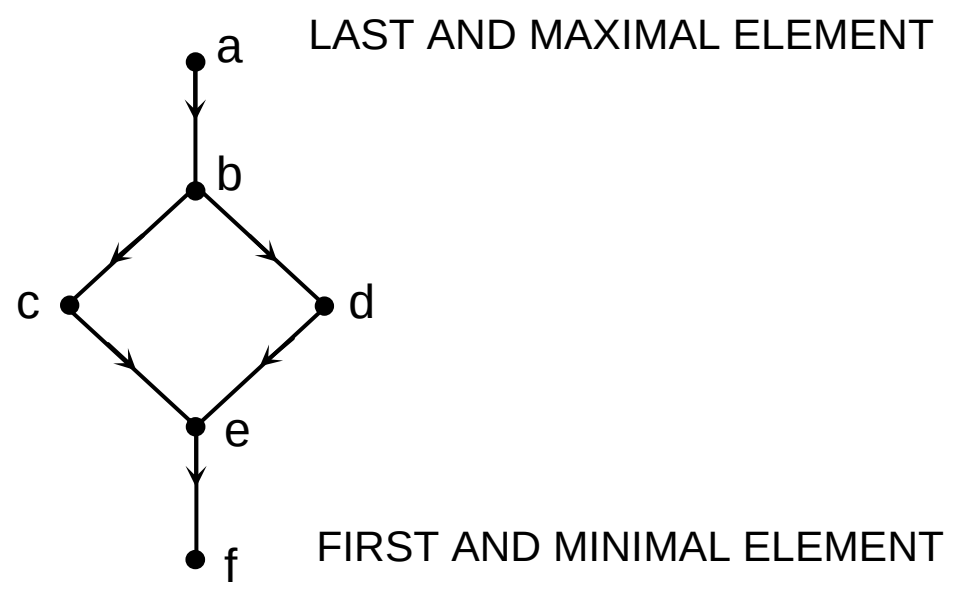

FIGURE: 2.8.7

The membership matrix of R

$$R = \begin{array}{c} \\ a \\ b \\ c \\ d \\ r \\ f \end{array} \begin{array}{c} \begin{array}{cccccc} a & b & c & d & e & f \end{array} \\ \left[\begin{array}{cccccc} 1 & 0 & 0 & 0 & 0 & 0 \\ \underline{1} & 1 & 0 & 0 & 0 & 0 \\ 1 & \underline{1} & 1 & 0 & 0 & 0 \\ \underline{1} & 0 & 0 & 1 & 0 & 0 \\ 1 & 1 & \underline{1} & \underline{1} & 1 & 0 \\ 1 & 1 & 1 & 1 & \underline{1} & 1 \end{array}\right] \end{array}$$

The related graph

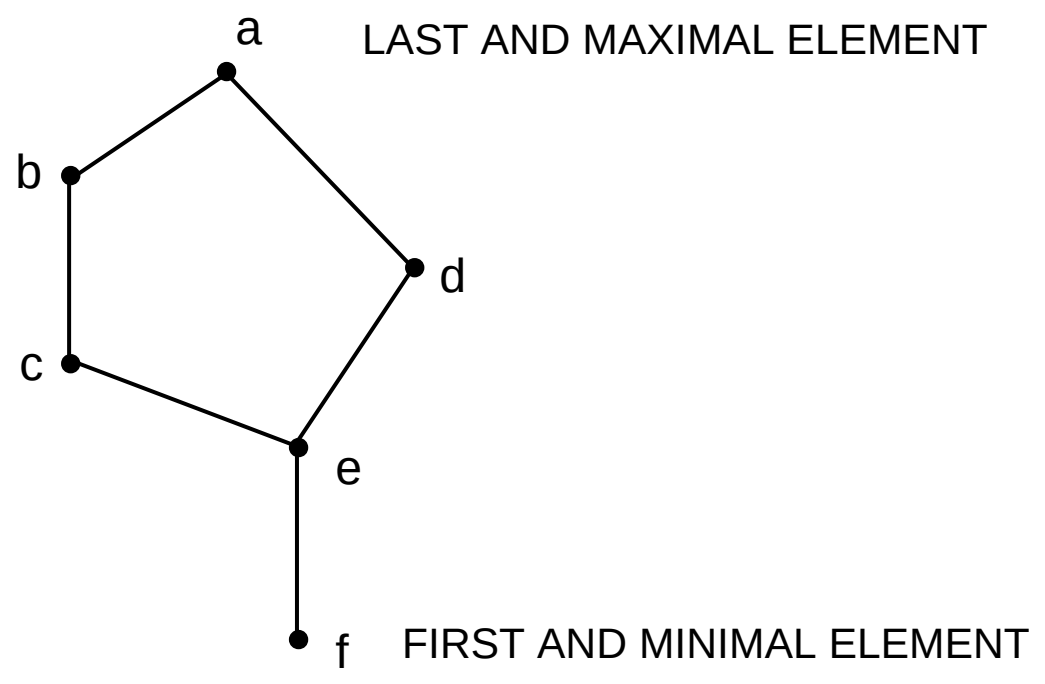

FIGURE: 2.8.8

Next we recall the definition of morphism and there graphical representation.

**DEFINITION 2.8.4**: *Let R (X, X) and Q (Y, Y) be fuzzy binary relations on the sets X and Y respectively A function h : X → Y is said to be a fuzzy homomorphism from (X, R) to (Y Q) if R ($x_1$, $x_2$)≤ Q (h ($x_1$), h ($x_2$) for all $x_1$, $x_2$ ∈ X and their images h ($x_1$), h ($x_2$)∈ Y.*

*Thus, the strength of relation between two elements under R is equated or excepted by the strength of relation between their homomorphic images under Q.*

It is possible for a relation to exist under Q between the homomorphic images of two elements that are themselves unrelated under R.

When this is never the case under a homomorphic function h, the function is called a strong homomorphism. It satisfies the two implications ($x_1$, $x_2$) ∈ R implies (h ($x_1$), h ($x_2$)) ∈ Q for all $x_1$, $x_2$ ∈ X and ($y_1$, $y_2$) ∈ Q implies ($x_1$, $x_2$) ∈ R for all $y_1$ ,$y_2$ ∈ Y where $x_1 \in h^{-1}$ ($y_1$) and $x_2 \in h^{-1}$ ($y_2$).

***Example 2.8.5:*** Let X = {a, b, c} y = {α, β, γ, δ} be sets with the following membership matrices which represent the fuzzy relations R (X, X) and Q (Y, Y)

$$R = \begin{matrix} \\ a \\ b \\ c \end{matrix} \overset{\begin{matrix} \alpha & \beta & \gamma & \delta \end{matrix}}{\begin{bmatrix} 0 & .5 & 0 & 0 \\ 0 & 0 & .9 & 0 \\ 1 & 0 & 0 & .5 \\ 0 & .6 & 0 & 0 \end{bmatrix}}$$

and

$$Q = \begin{matrix} \alpha \\ \beta \\ \gamma \\ \delta \end{matrix} \overset{\begin{matrix} \alpha & \beta & \gamma & \delta \end{matrix}}{\begin{bmatrix} 0 & .5 & 0 & 0 \\ 0 & 0 & .9 & 0 \\ 1 & 0 & 0 & .5 \\ 0 & .6 & 0 & 0 \end{bmatrix}}$$

The graph of the ordinary fuzzy homomorphism.

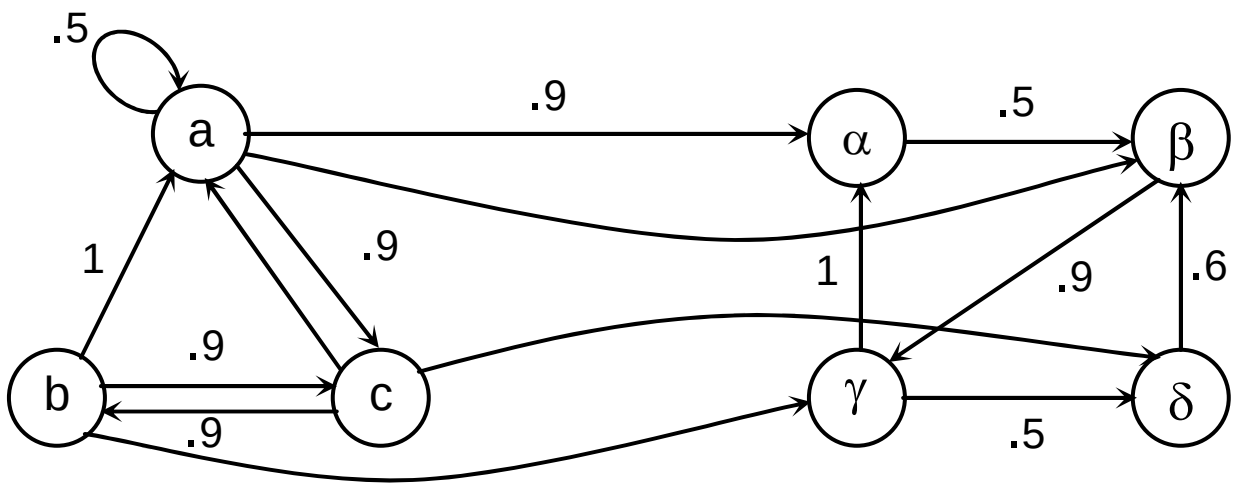

FIGURE: 2.8.9

Now we leave it for the reader to construct a strong fuzzy homomorphism and illustrate it by a graph. Next we just recall the definition of Fuzzy Cognitive maps (FCMs) and illustrate the directed graph of a FCM by an example.

**DEFINITION 2.8.5:***An FCM is a directed graph with concepts like policies, events etc as nodes and causalities as edges. It represents causal relationship between concepts.*

The problem studied in this case is for a fixed source S, a fixed destination D and a unique route from the source to the destination, with the assumption that all the passengers travel in the same route, we identify the preferences in the regular services at the peak hour of a day.

We have considered only the peak-hour since the passenger demand is very high only during this time period, where the transport sector caters to the demands of the different groups of people like the school children, the office goers, the vendors etc.

We have taken a total of eight characteristic of the transit system, which includes the level of service and the convenience factors.

We have the following elements, Frequency of the service, in-vehicle travel time, the travel fare along the route, the speed of the vehicle, the number of intermediate points, the waiting time, the number of transfers and the crowd in the bus or equivalently the congestion in the service.

Before defining the cognitive structure of the relationship, we give notations to the concepts involved in the analysis as below.

| | | |
|---|---|---|
| $C_1$ | - | Frequency of the vehicles along the route |
| $C_2$ | - | In-vehicle travel time along the route |
| $C_3$ | - | Travel fare along the route |
| $C_4$ | - | Speed of the vehicles along the route |
| $C_5$ | - | Number of intermediate points in the route |
| $C_6$ | - | Waiting time |
| $C_7$ | - | Number of transfers in the route |
| $C_8$ | - | Congestion in the vehicle. |

The graphical representation of the inter-relationship between the nodes is given in the form of directed graph given in Figure: 2.8.10.

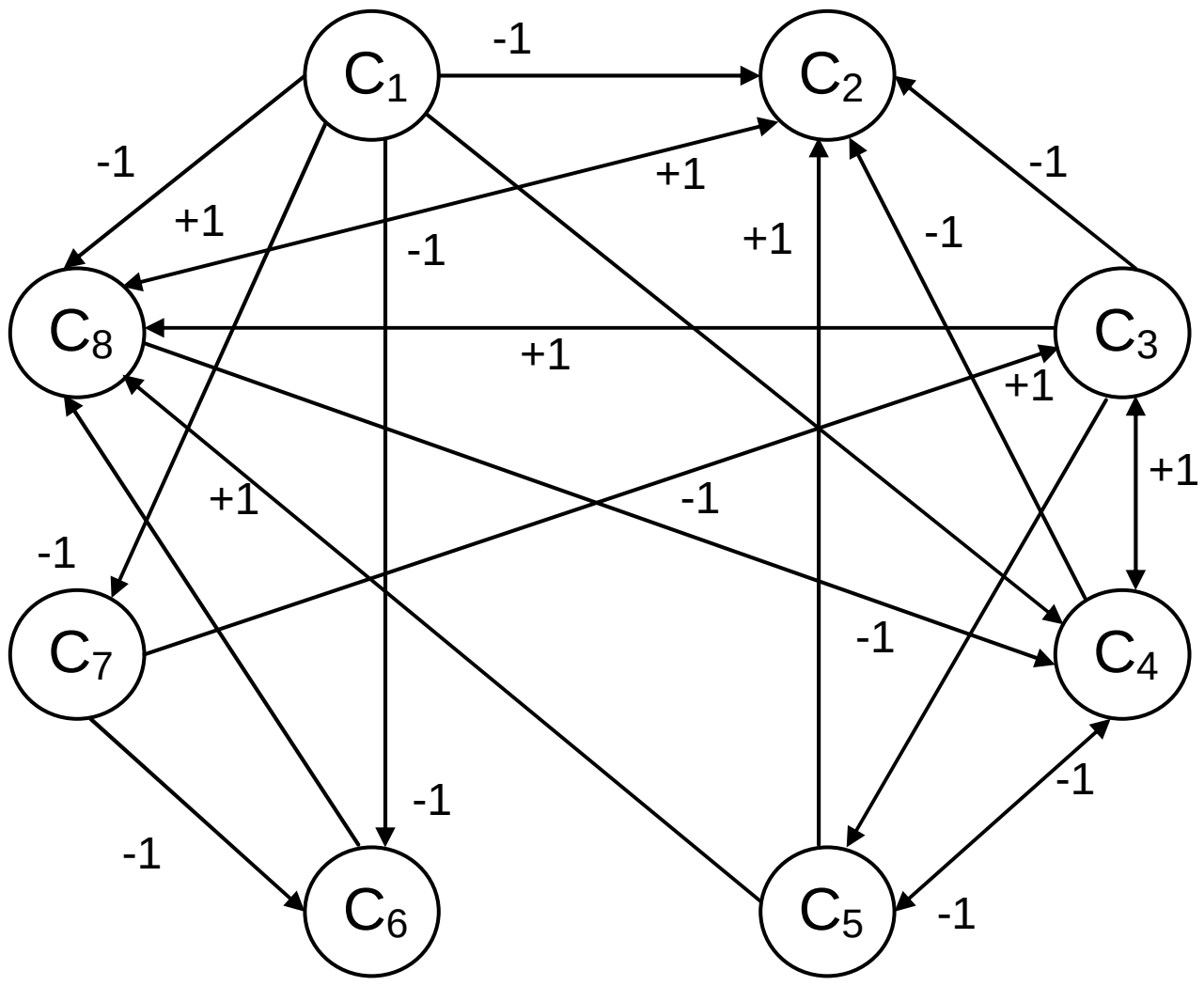

FIGURE: 2.8.10

From the above signed directed graph, we obtain a connection matrix E, since the number of concepts used here are eight, the connection matrix is a $8 \times 8$ matrix.

Thus we have $E = [A_y]_{8\times8}$

$$E = \begin{bmatrix} 0 & -1 & 0 & 1 & 0 & -1 & -1 & -1 \\ 0 & 0 & 0 & 0 & 0 & 0 & 0 & 1 \\ 0 & -1 & 0 & 1 & -1 & 0 & 0 & -1 \\ 0 & -1 & 1 & 0 & -1 & 0 & 0 & 0 \\ 0 & 1 & 0 & -1 & 0 & 0 & 0 & 1 \\ 0 & 0 & 0 & 0 & 0 & 0 & 0 & 1 \\ 0 & 0 & 1 & 0 & 0 & -1 & 0 & 0 \\ 0 & 1 & 0 & -1 & 0 & 0 & 0 & 0 \end{bmatrix}.$$

Now we just recall the definition of Fuzzy Relational Maps. (FRM).

**DEFINITION 2.8.6:** *A FRM is a directed graph or a map from D to R with concepts like policies or events etc as nodes and causalities as edges. It represents causal relation between spaces D and R.*

The employee-employer relationship is an intricate one. For, the employers expect to achieve performance in quality and production in order to earn profit, on the other hand employees need good pay with all possible allowances. Here we have taken three experts opinion in the study of Employee and Employer model.

The three experts whose opinions are taken are the Industry Owner, Employees' Association Union Leader and an Employee.

The data and the opinion are taken only from one industry. Using the opinion we obtain the hidden patterns. The following concepts are taken as the nodes relative to the employee.

We can have several more nodes and also several experts' opinions for it a clearly evident theory which professes that more the number of experts the better is the result.

We have taken as the concepts / nodes of domain only 8 notions which pertain to the employee.

$D_1$ – Pay with allowances and bonus to the employee
$D_2$ – Only pay to the employee
$D_3$ – Pay with allowances (or bonus) to the employee
$D_4$ – Best performance by the employee
$D_5$ – Average performance by the employee
$D_6$ – Poor performance by the employee
$D_7$ – Employee works for more number for hours
$D_8$ – Employee works for less number of hours.

$D_1, D_2, \ldots, D_8$ are elements related to the employee space which is taken as the domain space.

We have taken only 5 nodes / concepts related to the employer in this study.

These concepts form the range space which is listed below.

| | | |
|---|---|---|
| $R_1$ | – | Maximum profit to the employer |
| $R_2$ | – | Only profit to the employer |
| $R_3$ | – | Neither profit nor loss to the employer |
| $R_4$ | – | Loss to the employer |
| $R_5$ | – | Heavy loss to the employer |

The directed graph as given by the employer is given in Figure 1.6.1.

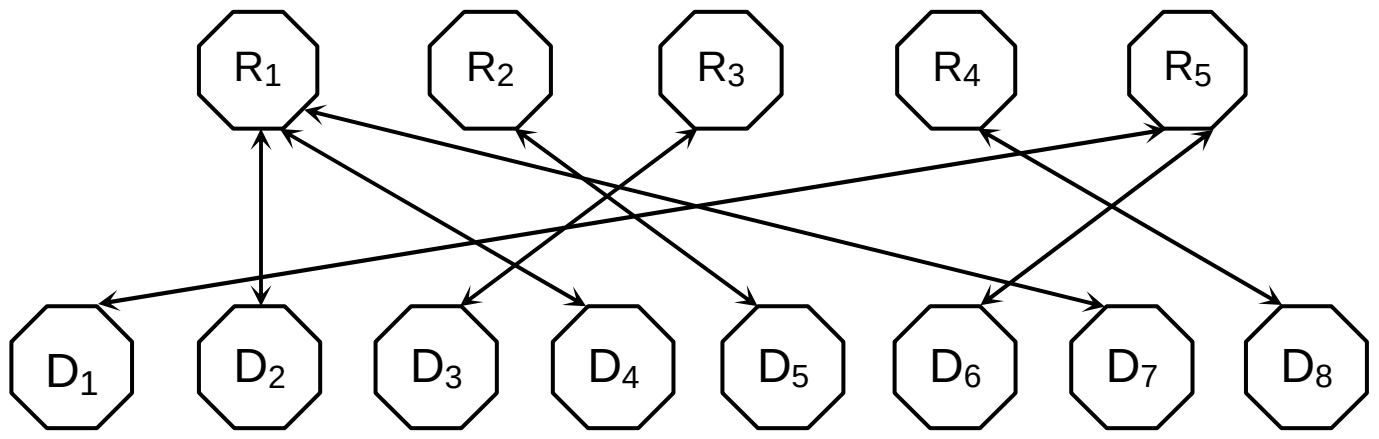

FIGURE:2.8.11

The associated relational matrix $E_1$ of the employer as given by following.

$$E_1 = \begin{bmatrix} 0 & 0 & 0 & 0 & 1 \\ 1 & 0 & 0 & 0 & 0 \\ 0 & 0 & 1 & 0 & 0 \\ 1 & 0 & 0 & 0 & 0 \\ 0 & 1 & 0 & 0 & 0 \\ 0 & 0 & 0 & 0 & 1 \\ 1 & 0 & 0 & 0 & 0 \\ 0 & 0 & 0 & 1 & 0 \end{bmatrix}$$

Thus we several of the fuzzy models make use of the graphs.

Chapter Three

# NEUTROSOPHIC GRAPHS AND THEIR APPLICATIONS TO NEUTROSOPHIC MODELS

In this chapter we introduce the notion of neutrosophic graphs and bring in several of its analogous properties and its applications. The study of neutrosophy is giving proper representation to the concept of indeterminacy present in all problems be it real or otherwise.

Neutrosophy was studied by [39, 53-4]. We have studied neutrosophic models like neutrosophic morphism, neutrosophic relational maps and neutrosophic cognitive maps. All these models exploit the notion the neutrosophic graphs [64-5]. So to study neutrosophic graphs one needs the concept of neutrosophic matrices and hence the neutrosophic field. Thus a few of the neutrosophic algebraic structures were introduced in the chapter I of this book. Also we use the notion of neutrosophic directed graphs and neutrosophic bipartite graphs. So we in this chapter give several of the neutrosophic analogous of basic graph theory results.

This chapter has eight sections. The first two sections gives the definition of neutrosophic point / vertex graphs and a few of its properties. In section three we define the notion of neutrosophic regular graphs and show a neutrosophic regular graph in general is not a neutrosophic strongly regular graph. Section four is devoted to the introduction of neutrosophic trees and neutrosophic Eulerian graphs. Just we introduce the notion of neutrosophic graph colourings in section five.

The famous neutrosophic Petersen graphs are introduced in section six. Unlike in graphs in which we have only one Petersen graph we have several neutrosophic Petersen graphs. The seventh section deals with the application of neutrosophic graphs to neutrosophic models. The neutrosophic graphs of the these models are fully explained and illustrated. The main reason for it is that in case of fuzzy models, we have been studying and have used fuzzy models for at least 2 decades but these neutrosophic models are very few for they have been built, defined and put to use only in 2003 and 2004. The final section compares both fuzzy and neutrosophic notions

To the best of our knowledge we do not have any books or papers on these models except [64-6]. So we felt it necessary to explain them fully. Thus this section completely gives the methods and implementation with a motive to make this book a self contained one. Thus a complete section is devoted to show the application of neutrosophic graphs in these neutrosophic models.

## 3.1. Neutrosophic point / vertex Graphs and its properties

The term neutrosophic vertex means a node or a vertex, which is an indeterminable. Thus when we say neutrosophy we mean the concept / attribute / the node of it is not determinable may be at that time or for some specified interval of time or under some specified conditions, varying with the time and the related conditions.

Thus the indeterminable node may after a period of time and with additional circumstances may become partially determinable or determinable depending mainly on the problem under investigation. Thus this new notion of neutrosophic vertex graphs will finds its applications in real world problems, like NCMs and NRMs.

**DEFINITION 3.1.1:** *A neutrosophic point graph $G_N$ is a graph G with finite non empty set $V_N = V_N$ (G) of p-points where at least one of the point in $V_N$ (G) are indeterminate node, element, point or vertex.*

Note here $V_N(G) = V(G) + N$ where V(G) are points or vertices of the graph G and N the non empty set of points which are indeterminate node.

***Example 3.1.1:*** $V_3$ (G) = V(G) ∪ N(N = {$N_1$, $N_2$, $N_3$}) $G_N$ is a neutrosophic point graph with 3 nodes which are indeterminates and V(G) is the graph, given in figure 3.1.1.

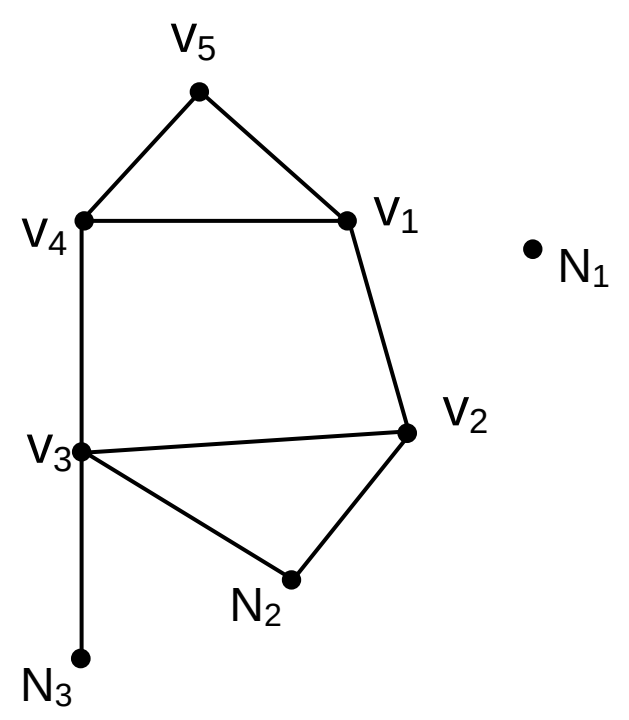

FIGURE: 3.1.1

**DEFINITION 3.1.2:** *Elements of $V_N(G)$, where $V_N(G)$ is a neutrosophic point graph are called the neutrosophic vertices of G. The number of elements in $V_N(G)$ is n(G) + N where n(G) is called the order of G and N is the number of indeterminate nodes used in $V_N(G)$.*

***Example 3.1.2:*** Thus in $V_5(G)$ where G is the graph with four elements

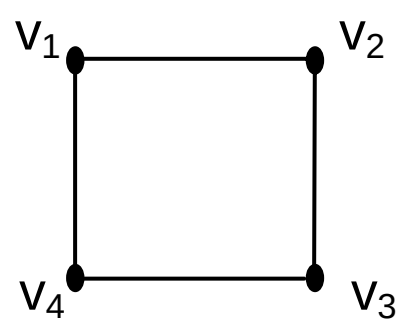

FIGURE: 3.1.2

then order of $V_5(G)$ is 4 + 5 i.e. 9.

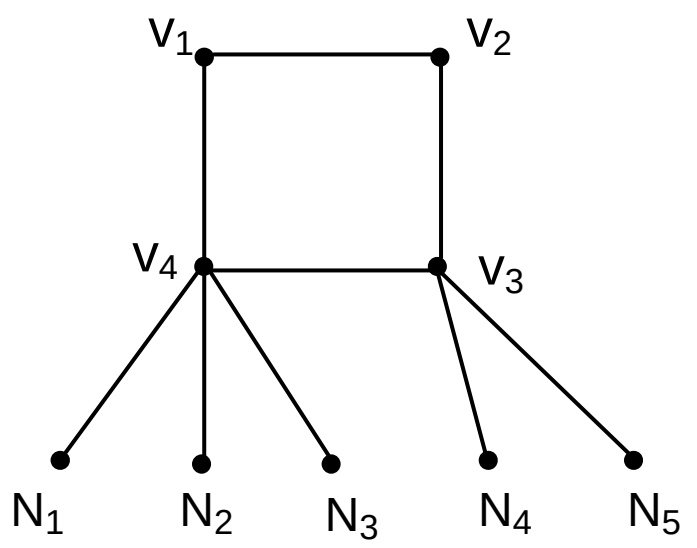

FIGURE: 3.1.3

$E(V_N(G))$ denotes the edge set of $V_N(G)$.

***Example 3.1.3:*** $V_3(G)$ is given by the following figure:

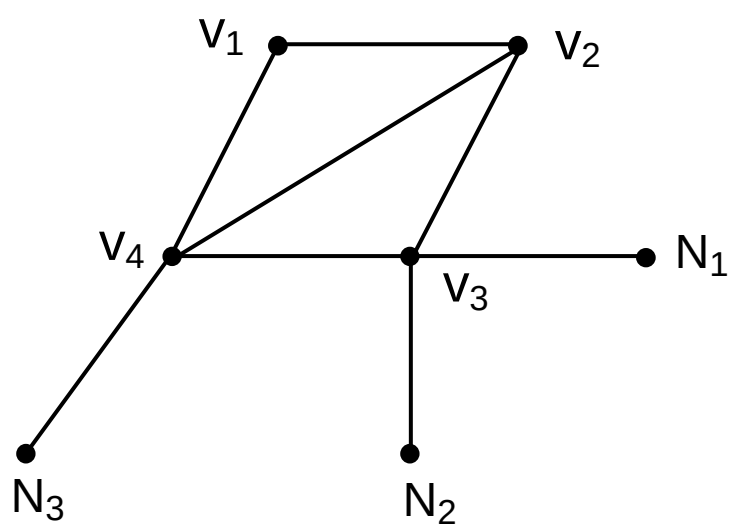

FIGURE: 3.1.4

The practical application of the neutrosophic point graph will be discussed in the last section of this chapter.

Now we proceed onto define the notion of neutrosophic edge graph. In a neutrosophic edge graph, we will have all the nodes to be real or none of the nodes will be indeterminate ones, only some of the edges relating the nodes will be indeterminates.

**DEFINITION 3.1.3:** *Let V(G) be the set of all vertices of the graph G. If the edge set E(G) where at least one of the edges of G is an indeterminate one. Then we call such graphs as neutrosophic edge graphs.*

***Example 3.1.4:*** The graph given in figure 3.1.5 has $(v_1, \ldots, v_5)$ as vertex set is a neutrosophic edge graph. The edges $\{v_2, v_3\}$ and $\{v_4, v_5\}$ are indeterminate edges so we see the graph is a neutrosophic edge graph.

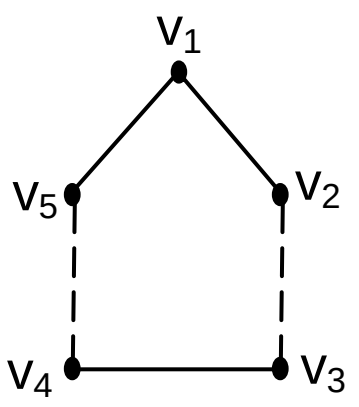

FIGURE: 3.1.5

Thus we see the neutrosophic vertex graph is distinctly different from the neutrosophic edge graph. They differ from each other on the edge set and the vertex set. The edge set of a neutrosophic vertex graph are all usual edges where as only the vertex sets are indeterminates, on the contrary the vertex set of the neutrosophic edge graph have the vertex set to be the usual set, the difference lies only in the edges set where some of the edges are indeterminates.

The neutrosophic edges graph has several applications like in the neutrosophic cognitive maps, neutrosophic relational maps and so on.

Thus from now on wards we make some compromise and call the neutrosophic edge graphs as just neutrosophic graphs. Thus neutrosophic vertex graphs will continue to be so. Hence the reformulated or to be more precise is to be restated.

The following are examples of neutrosophic graphs:

***Example 3.1.5:***

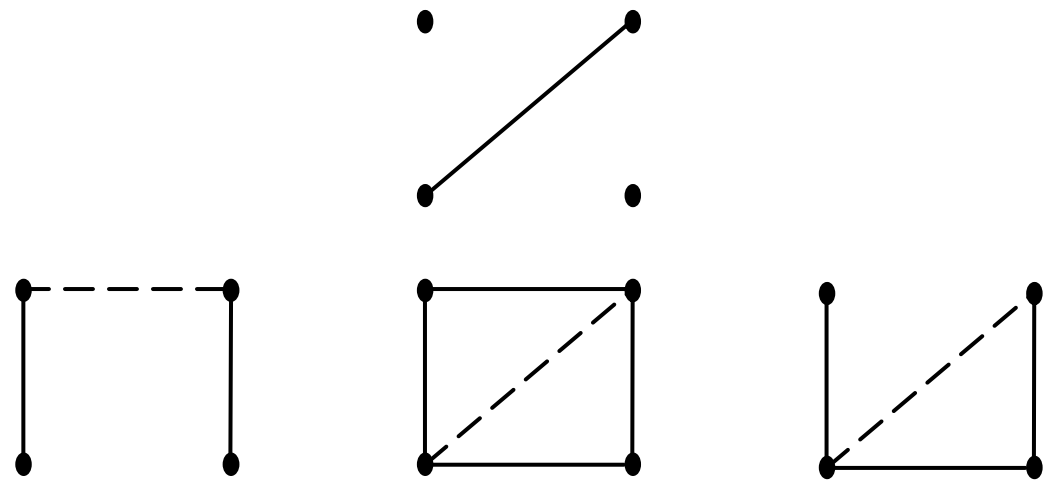

FIGURE: 3.1.6

All graphs in general are not neutrosophic graphs.

***Example 3.1.6:*** The following graphs are not neutrosophic graphs.

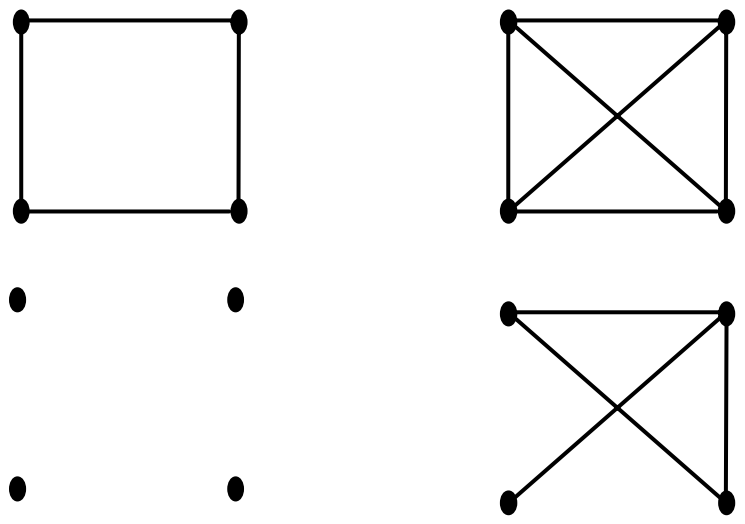

FIGURE: 3.1.7

## 3.2 Properties of Neutrosophic graphs

Now we proceed on to define the neutrosophic graph.

**DEFINITION 3.2.1:** *A neutrosophic graph is a graph in which at least one edge is an indeterminacy denoted by dotted lines.*

**NOTATION***:* The indeterminacy of an edge between two vertices will always be denoted by dotted lines.

***Example 3.2.1:*** The following are neutrosophic graphs:

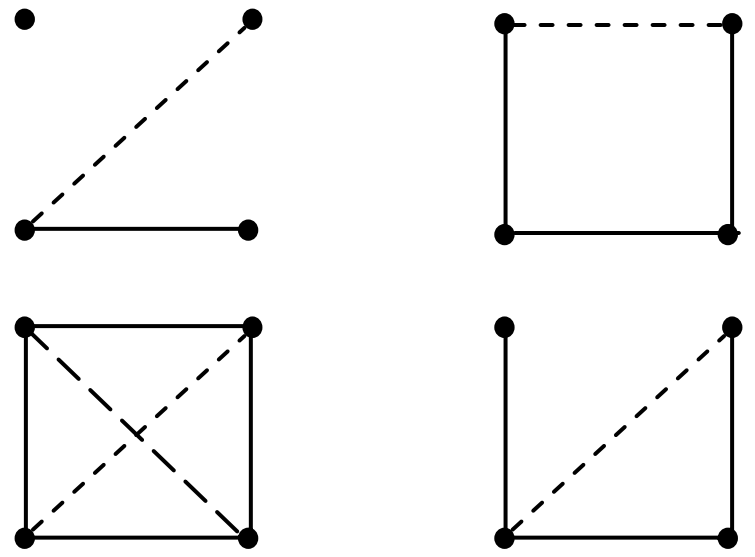

FIGURE: 3.2.1

All graphs in general are not neutrosophic graphs.

***Example 3.2.2:*** The following graphs are not neutrosophic graphs given in Figure 3.2.2:

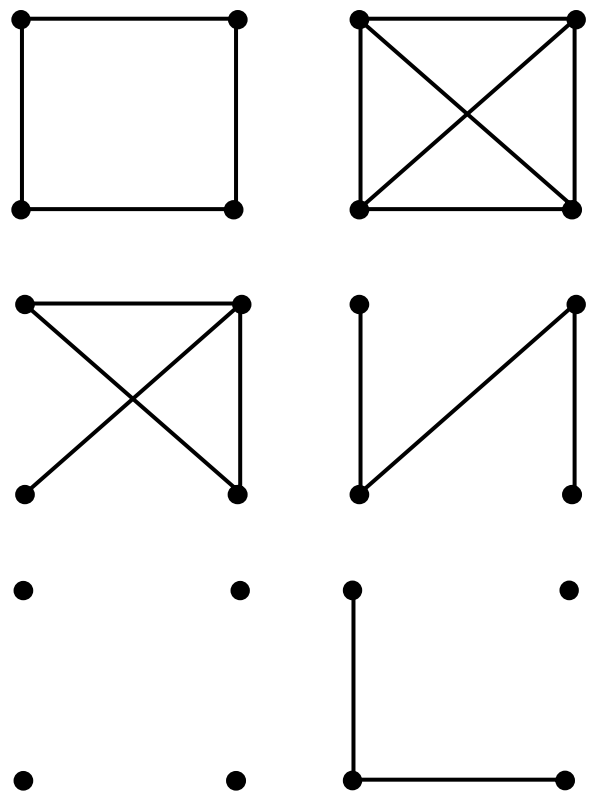

FIGURE: 3.2.2

**DEFINITION 3.2.2:** *A neutrosophic directed graph is a directed graph which has at least one edge to be an indeterminacy.*

**DEFINITION 3.2.3:** *A neutrosophic oriented graph is a neutrosophic directed graph having no symmetric pair of directed indeterminacy lines.*

**DEFINITION 3.2.4:** *A neutrosophic subgraph H of a neutrosophic graph G is a subgraph H which is itself a neutrosophic graph.*

**THEOREM 3.2.1:** *Let G be a neutrosophic graph. All subgraphs of G are not neutrosophic subgraphs of G.*

*Proof:* By an example. Consider the neutrosophic graph given in Figure 3.2.3.

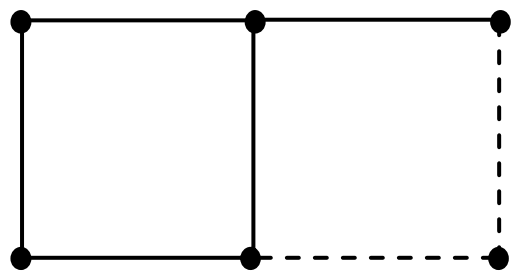

FIGURE: 3.2.3

This has a subgraph given by Figure 3.2.4.

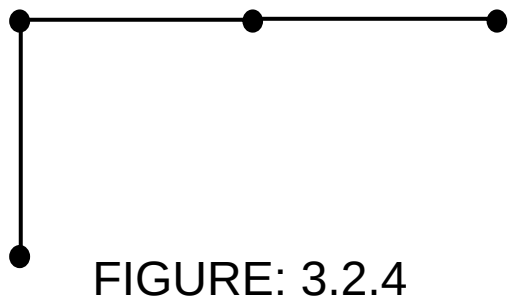

FIGURE: 3.2.4

which is not a neutrosophic subgraph of G.

**THEOREM 3.2.2:** *Let G be a neutrosophic graph. In general the removal of a point from G need not be a neutrosophic subgraph.*

*Proof:* Consider the graph G given in Figure 3.2.5.

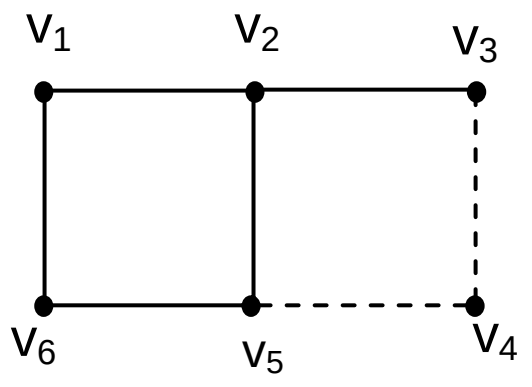

FIGURE: 3.2.5

G \ $v_4$ is only a subgraph of G but is not a neutrosophic subgraph of G.

Thus it is interesting to note that this is a main feature by which a graph differs from a neutrosophic graph.

**DEFINITION 3.2.5:** *Two graphs G and H are neutrosophically isomorphic if*

*i. They are isomorphic.*
*ii. If there exists a one to one correspondence between their point sets which preserve indeterminacy adjacency.*

**DEFINITION 3.2.6:** *A neutrosophic walk of a neutrosophic graph G is a walk of the graph G in which at least one of the lines is an indeterminacy line. The neutrosophic walk is neutrosophic closed if* $v_0 = v_n$ *and is neutrosophic open otherwise.*

*It is a neutrosophic trial if all the lines are distinct and at least one of the lines is a indeterminacy line and a path, if all points are distinct (i.e. this necessarily means all lines are distinct and at least one line is a line of indeterminacy). If the neutrosophic walk is neutrosophic closed then it is a neutrosophic cycle provided its n points are distinct and* $n \geq 3$.

*A neutrosophic graph is neutrosophic connected if it is connected and at least a pair of points are joined by a path. A neutrosophic maximal connected neutrosophic subgraph of G is called a neutrosophic connected component or simple neutrosophic component of G.*

*Thus a neutrosophic graph has at least two neutrosophic components then it is neutrosophic disconnected. Even if one is a component and another is a neutrosophic component still we do not say the graph is neutrosophic disconnected.*

***Example 3.2.3:*** Neutrosophic disconnected graphs are given in Figure 3.2.6.

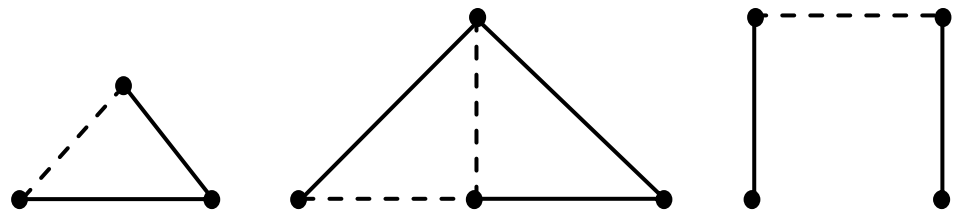

FIGURE: 3.2.6

***Example 3.2.4:*** Graph which is not neutrosophic disconnected is given by Figure 3.2.7.

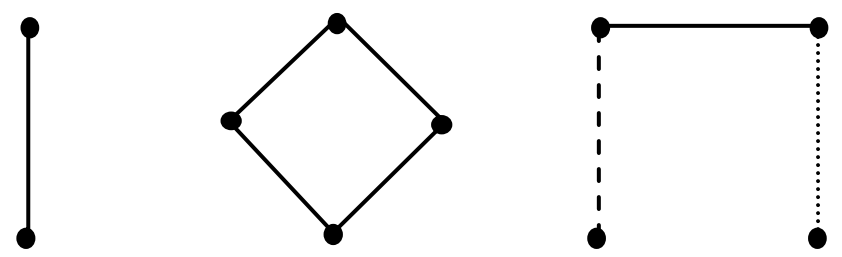

FIGURE: 3.2.7

Several results in this direction can be defined and analyzed.

**DEFINITION 3.2.7:** *A neutrosophic bigraph, G is a bigraph, G whose point set V can be partitioned into two subsets $V_1$ and $V_2$ such that at least a line of G which joins $V_1$ with $V_2$ is a line of indeterminacy.*

This neutrosophic bigraphs will certainly play a role in the study of FRMs and in fact we give a method of conversion of data from FRMs to FCMs.

As both the models FRMs and FCMs work on the adjacency or the connection matrix we just define the neutrosophic adjacency matrix related to a neutrosophic graph G given by Figure 3.2.8.

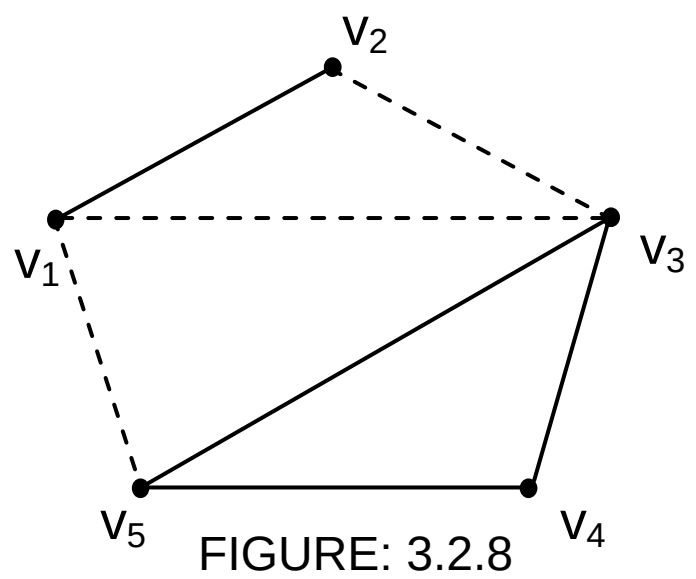

FIGURE: 3.2.8

The neutrosophic adjacency matrix is N(A)

$$N(A) = \begin{bmatrix} 0 & 1 & I & 0 & I \\ 1 & 0 & I & 0 & 0 \\ I & I & 0 & 1 & 1 \\ 0 & 0 & 1 & 0 & 1 \\ I & 0 & 1 & 1 & 0 \end{bmatrix}.$$

Its entries will not only be 0 and 1 but also the indeterminacy *I*.

**DEFINITION 3.2.8:** *Let G be a neutrosophic graph. The adjacency matrix of G with entries from the set (I, 0,* 1*) is called the neutrosophic adjacency matrix of the graph.*

Thus one finds a very interesting application of neutrosophy graphs in Neutrosophic Cognitive Maps and Neutrosophic Relational Maps. Application of these concepts will be dealt in the last section of this chapter.

## 3.3 Neutrosophic regular graphs and its properties

Now we proceed on to define the notion doubly neutrosophic graphs.

**DEFINITION 3.3.1:** *A graph G is said to be a doubly or strongly neutrosophic graph if the graph has both indeterminate vertices and indeterminate edges. The indeterminate edges as usual will be denoted by dotted lines where as the indeterminate vertices will be denoted by $N_1, \ldots, N_k$.*

***Example 3.3.1:***

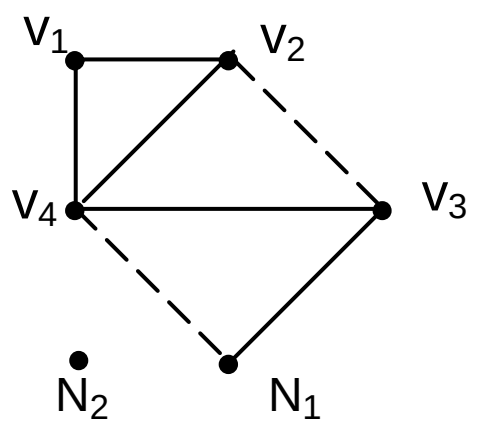

FIGURE: 3.3.1

The graph given in the figure is a strongly neutrosophic graph.

**NOTATION:** $V_N(G)$ will denote the vertices of the neutrosophic graph with N vertices, which are indeterminates and G the number of vertices that are not indeterminates. Thus the number of vertices in the graph is n(G) +N. The edge set $E_N$ (G) will include the edges as well as the dotted lines that is indeterminate edge or the neutrosophic edge. It is important to note that if $V_N$ (G) is a neutrosophic vertex graph then from removing the N indeterminate vertices $V_N(G)$ becomes a usual graph or neutrosophic graph which is not a neutrosophic strong graph.

In view of this we have the following theorem:

**THEOREM 3.3.1:** *Let G be a neutrosophic vertex graph, $V_N$ (G) its vertices. All subgraphs of $V_N(G)$ need not in general be neutrosophic vertex graphs.*

*Proof:* Let H be a subgraph of G such that V(H) has no neutrosophic vertex. Then the graph H is not a neutrosophic vertex graph.

We illustrate this by the following example:

***Example 3.3.2:***

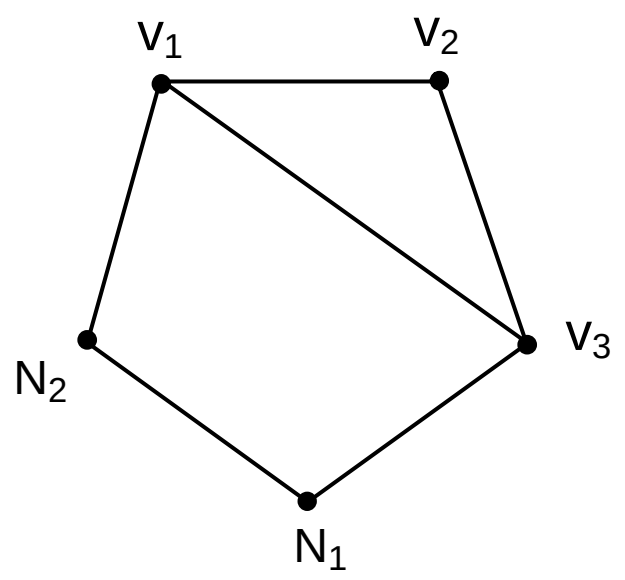

FIGURE: 3.3.2

$V_2$(G) is graph, which is a neutrosophic vertex graph. The subgraph H = {$v_1$, $v_2$, $v_3$} is a subgraph which is not a neutrosophic vertex graph.

Now we have to define the isomorphism of neutrosophic vertex graphs.

**DEFINITION 3.3.2:** *A neutrosophic vertex graph $G_N$ is said to be neutrosophic simple if the graph has no loops or multiple edges connecting indeterminate vertex or two indeterminate vertices.*

Thus a neutrosophic vertex simple graph $G_N$ need not in general be neutrosophic vertex simple.

***Example 3.3.3:***

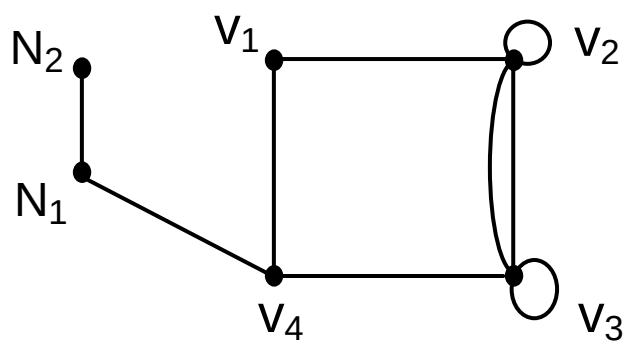

FIGURE: 3.3.3

The graph shown in the figure is neutrosophic vertex simple. Clearly the graph is not simple.

***Example 3.3.4:***

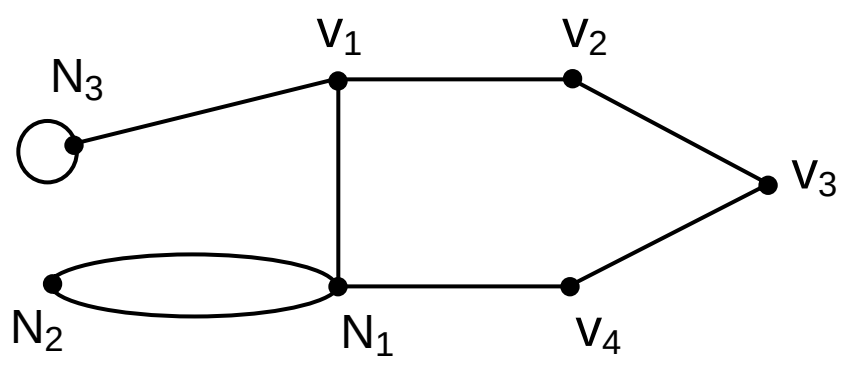

FIGURE: 3.3.4

The graph shown in this figure is not neutrosophic vertex simple. But it can be said to be a simple graph. So while defining properties of neutrosophic graphs we are more interested only in the study of its indeterminate vertices.

**DEFINITION 3.3.3:** *A neutrosophic vertex graph isomorphism from the neutrosophic vertex graphs $G_N$ to $H_N$ is as follows:*

*Let*

$$G_{N_1} = \{V(G_{N_1}), E(G_{N_1}), I_{GN_1}\} \text{ and}$$
$$H_{N_2} = (V(H_{N_2}), E(H_{N_2}), I_{HN_2})$$

*be neutrosophic vertex graphs. A neutrosophic vertex graph isomorphism from $G_{N_1}$ to $H_{N_2}$ is a pair ($\phi_N$, $\theta_N$) where*

$$\phi_N: V(G_{N_1}) \to V(H_{N_2})$$

*and such that $\phi$: $G \to H$ where*

$$G = (G_{N_1} \setminus N_1) \text{ and } (H = H_{N_2} \setminus N_2)$$

*is a graph isomorphism and $\phi_N$ maps elements of $N_1$ to $N_2$ (i.e. indeterminate elements are mapped onto indeterminate elements only) from G to H and no intermingling of indeterminate vertex and vertex of G takes place. $\theta_N : E\ (G_{N_1}) \to E\ (H_{N_2})$.*

**DEFINITION 3.3.4:** *Let $G_N$ be a neutrosophic vertex graph. Let $N_1 \in G_N$ i.e. $N_1$ is an indeterminate vertex of $G_N$. The number of* *edges incident of $N_1 \in G_N$ is called the neutrosophic degree (or valiancy) of the neutrosophic vertex $N_1$ in G and is denoted by $d_{G_N}$ ($N_1$) or simply ($N_1$).*

A loop at any neutrosophic vertex of a graph $G_N$ as in case of graphs is counted twice.

The minimum (or respectively the maximum) of neutrosophic degrees of the indeterminate vertices of a graph $G_N$ is denoted by $\delta_N$ ($G_N$) or $\delta_N$ (respectively $\Delta_N$ (G) or $\Delta_N$). A neutrosophic vertex graph $G_N$ is K-neutrosophic regular if every indeterminate vertex of $G_N$ has degree K.

A neutrosophic vertex graph $G_N$ is neutrosophic regular if it is K-neutrosophic regular. A neutrosophic vertex graph is neutrosophic strongly regular if every vertex of $G_N \setminus N$ is K-regular and $G_N$ is K-neutrosophic regular. Thus we have the following theorem:

**THEOREM 3.3.2:** *Let $G_N$ be a neutrosophic strongly regular graph then $G_N$ is neutrosophic regular. However a neutrosophic regular graph need not be neutrosophic strongly regular.*

*Proof:* Follows from the very definition.

However to prove the neutrosophic regular graph in general need not be neutrosophic strongly regular we give the following example. Consider the graph $G_3$ given by the following example:

***Example 3.3.5:*** $G_3$ is a neutrosophic regular graph, but $G_3$ is not a neutrosophic strongly regular graph.

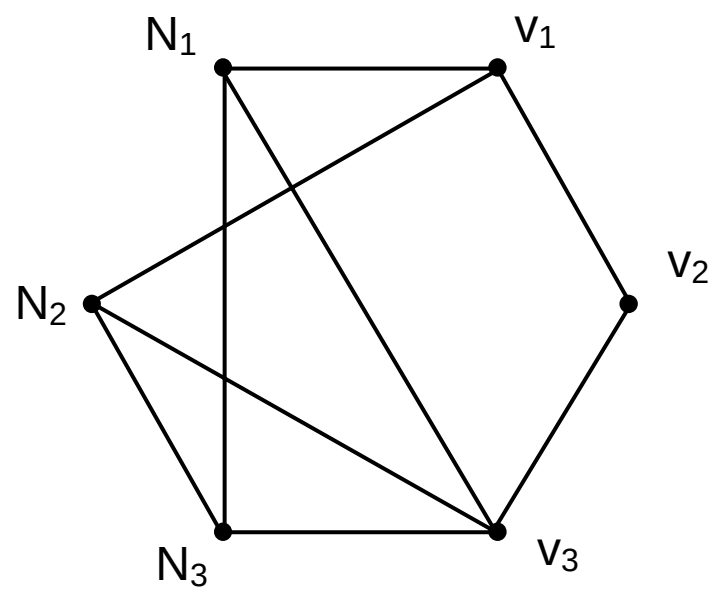

FIGURE: 3.3.5

***Example 3.3.6:*** The following example gives a neutrosophic strongly regular graph $G_3$:

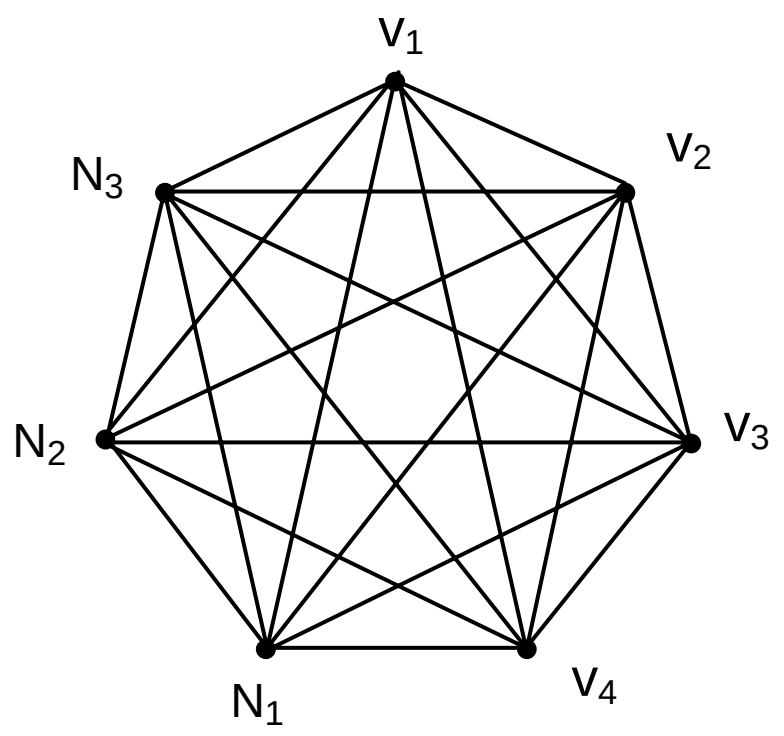

FIGURE: 3.3.6

This is in fact a 6-regular graph.

**DEFINITION 3.3.5:** *A indeterminate vertex of degree zero in a neutrosophic vertex graph is an isolated neutrosophic vertex of* $G_N$.

*A pendent neutrosophic vertex of* $G_N$ *is an indeterminate vertex of degree 1 and the unique edge of* $G_N$ *incident to such a vertex of* $G_N$ *is called the neutrosophic pendent edge of* $G_N$.

***Example 3.3.7:*** $G_5$ given by the following diagram has neutrosophic isolated vertex:

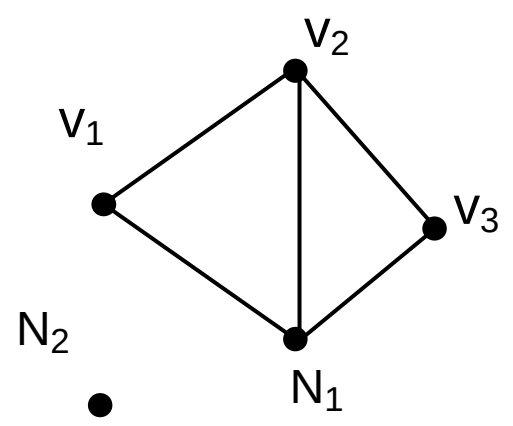

FIGURE: 3.3.7

$N_2$ is the neutrosophic-isolated vertex of $G_5$.

***Example 3.3.8:***

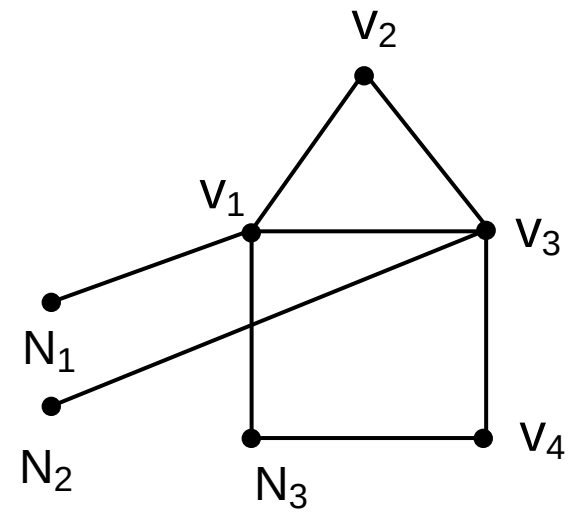

FIGURE: 3.3.8

The graph $G_3$ has neutrosophic pendent vertices at $N_1$ and $N_2$ and hence they have neutrosophic pendent edges. However this graph has no pendent vertex.

Euler theorem can be extended in case of neutrosophic graphs also.

A neutrosophic subgraph $H_N$ of $G_N$ is a neutrosophic spanning subgraph of $G_N$ if $V(H_N) = V(G_N)$. A neutrosophic vertex subgraph $H_{N_1}$ of $G_N$ is said to be an neutrosophic-induced subgraph of $G_N$ if each edge of $G_N$ having its ends in V ( $H_{N_1}$ ) is also an edge of $H_N$.

A neutrosophic vertex subgraph $H_{N_1}$ of $G_N$ is a neutrosophic spanning subgraph of G if V ( $H_{N_1}$ ) = $V(G_N)$.

Later we will be defining their concepts in case of neutrosophic vertex graphs or neutrosophic graphs and in case of doubly or strongly neutrosophic graphs.

**DEFINITION 3.3.6:** *A neutrosophic walk of a neutrosophic graph is a walk of the graph with an alternating sequence of points and lines where at least one of the point must necessarily be indeterminate vertex and in the lines or the edges at least one of the edge is an indeterminate one. If in the neutrosophic walk of a neutrosophic graph in which both $v_o$ and $v_n$ are indeterminate or in which both $v_o$ and $v_n$ are real vertices is said to be neutrosophic closed if $v_o = v_n$. ($v_o$ ,$v_n \in N$) and neutrosophic open otherwise.*

*It is a neutrosophic trial if lines or edges are distinct and at least one edge is an indeterminate and a neutrosophic path if all the points or vertices are distinct with at least one point or vertex to be an indeterminate. If the neutrosophic walk is closed then it is a neutrosophic cycle provided its n points are distinct n ≥ 3 of which at least one point or vertex in n is an indeterminate point or vertex.*

***Example 3.3.9:*** The labeled graph $G_3$ is a neutrosophic walk which is not a neutrosophic trial i.e. i.e. $v_1v_2$ $N_1$ $v_2$ $v_3$ is a neutrosophic walk and is not a neutrosophic trial

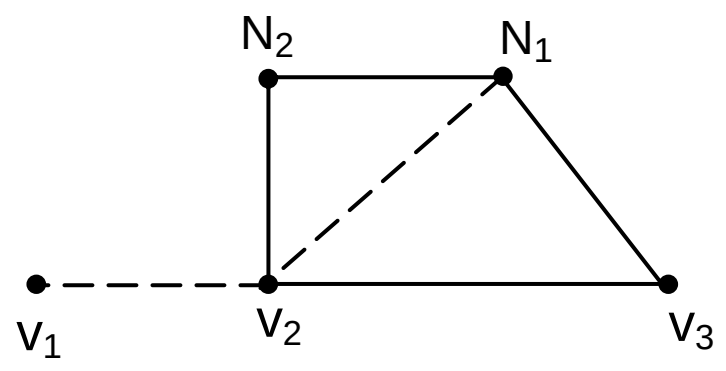

FIGURE: 3.3.9

$v_1v_2$ $N_1$ $N_2$ $v_2$ $v_3$ is a neutrosophic trial which is not a neutrosophic path.

$v_1v_2$ $N_1$ $N_2$ is a neutrosophic path and $v_2$ $N_2$ $N_1$ $v_2$ is a cycle and not a neutrosophic cycle. We denote by $NC_N$ the neutrosophic graph consisting of neutrosophic cycle with N points and by $NP_N$ a neutrosophic path with N points. $NC_3$ is often called a neutrosophic triangle.

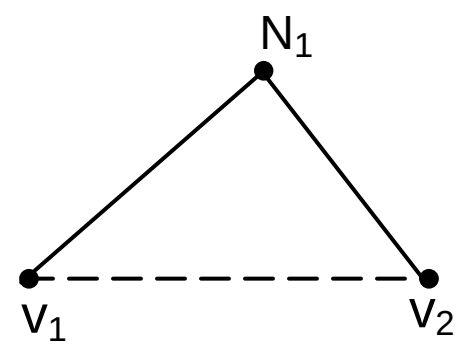

FIGURE: 3.3.10

A neutrosophic graph is connected if every pair of points is joined by paths (Here G has both indeterminate vertex and indeterminate edge). A maximal connected neutrosophic subgraph of G is called a neutrosophic connected component or simply neutrosophic component of G.

Thus the component of a graph G is not a neutrosophic component of G. It is interesting question to find out when does all neutrosophic components contain completely the components of the neutrosophic graphs. A neutrosophic disconnected graph has at least two neutrosophic components.

***Example 3.3.10:*** This is not a neutrosophic disconnected graph.

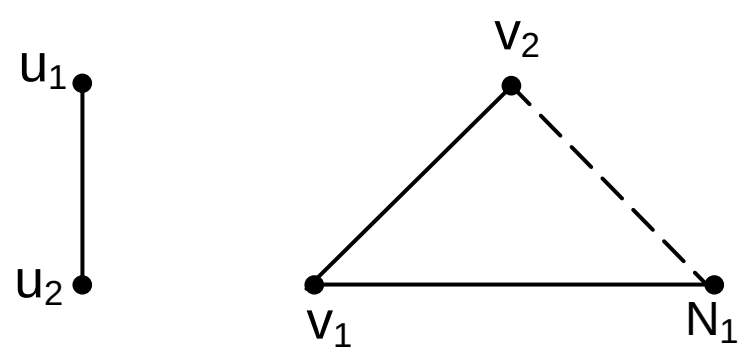

FIGURE: 3.3.11

***Example 3.3.11:*** This example is a neutrosophic disconnected graph.

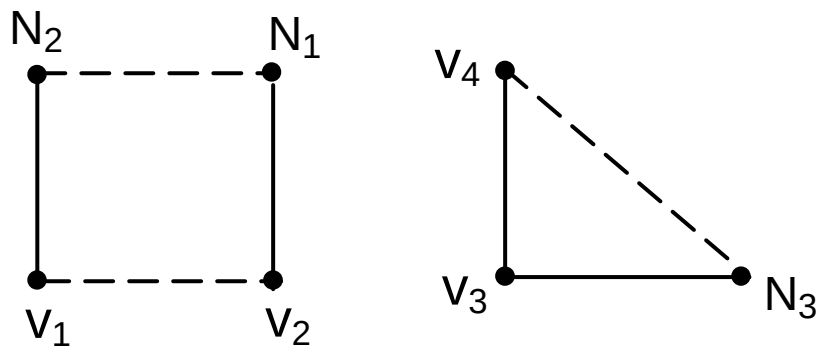

FIGURE: 3.3.12

***Example 3.3.12:***

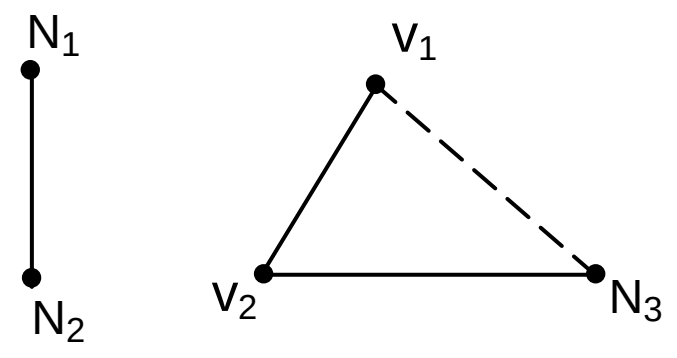

FIGURE: 3.3.13

The graph given in the above example in figure 3.3.13 is also not a neutrosophic disconnected graph.

The study of neutrosophic vertex graph, neutrosophic edge graph and strongly or doubly neutrosophic graphs was introduced and analyzed in the earlier sections.

In this section we introduce the notion of neutrosophic graphs and study them. As trees form an important class of graphs neutrosophic trees will certainly form an important class of neutrosophic graphs. Thus the study of neutrosophic trees is very important.

**DEFINITION 3.4.1:** *A neutrosophic tree is a neutrosophic graph which is neutrosophic connected without cycles. A neutrosophic graph without cycles is called a neutrosophic acyclic graph or a neutrosophic forest. Hence each component of a neutrosophic forest is a neutrosophic tree.*

We can equivalently give a definition for a simple neutrosophic graph to be a tree.

**DEFINITION 3.4.2:** *A simple neutrosophic graph G is a neutrosophic tree (G is a simple graph) if and only if any two distinct indeterminant vertices are connected by a unique neutrosophic path.*

A spanning neutrosophic subgraph of a neutrosophic graph, which is also, a neutrosophic tree called a spanning neutrosophic tree of the neutrosophic graph.

***Example 3.4.1:***

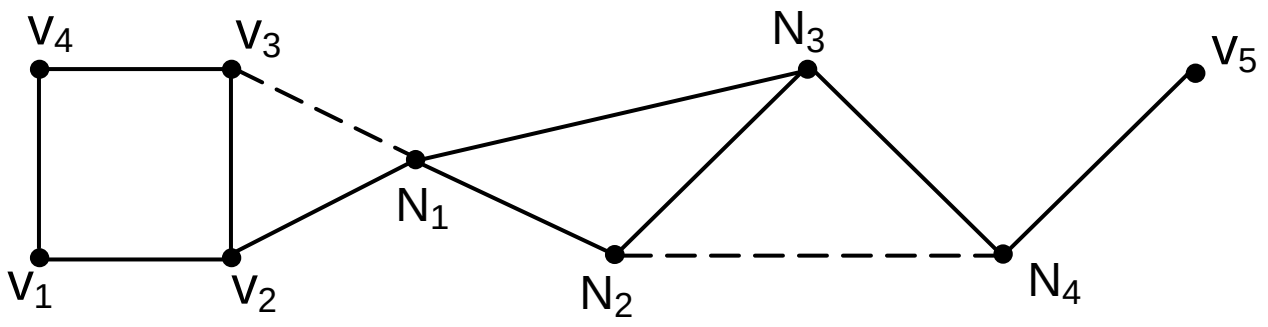

FIGURE: 3.4.1

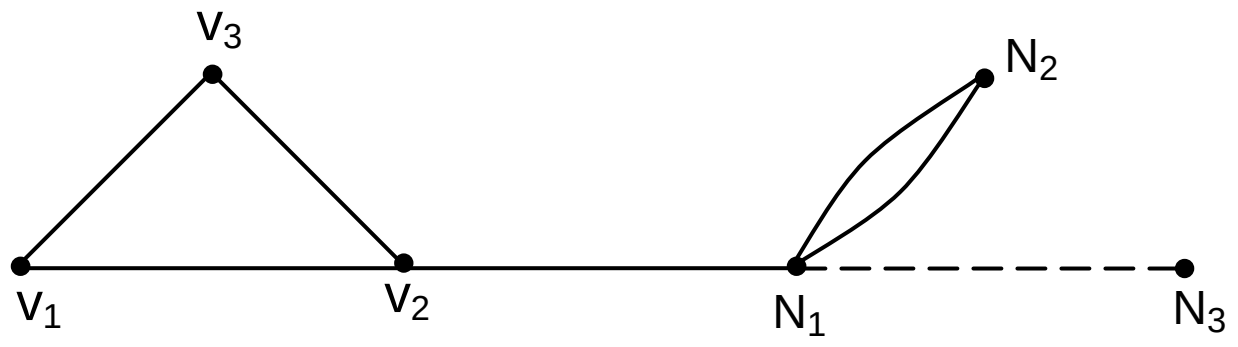

FIGURE: 3.4.2

The above two graphs shows the neutrosophic trees.

**DEFINITION 3.4.3:** *Let $G_N$ be a neutrosophic connected graph the neutrosophic diameter of $G_N$ is defined as max $\{d(n_1, n_2) \mid n_1, n_2 \in N\}$ and is denoted by diam ($G_N$).*

*If n is an indeterminate vertex of $G_N$ its neutrosophic eccentricity e (n) is defined by*

$$\max \{d(n, m) \mid m \in N\}.$$

*The neutrosophic radius of $G_N$, $r(G_N)$ is the minimum neutrosophic eccentricity of $G_N$, that is $r(G_N) = \min \{e(n) \mid n \in N\}$. Note that diam*

$$(G_N) = \max \{e(n) \mid n \in N\}.$$

*A indeterminate vertex n of $G_N$ is called a neutrosophic central vertex if $e(n) = r(G_N)$.*

*The set of all neutrosophic central vertices of $G_N$ is called the neutrosophic center of $G_N$.*

**DEFINITION 3.4.4:** *A neutrosophic Eulerian graph G is one which is a Eulerian graph such that at least one edge or one vertex is an indeterminate. 'or' not used in the mutually exclusive sense.*

*If the neutrosophic Eulerian graph G has both, at least one vertex and one edge to be an indeterminate then G is called as the strong or double neutrosophic Eulerian graph.*

First we give an example of a neutrosophic Eulerian graph.

***Example 3.4.2:***

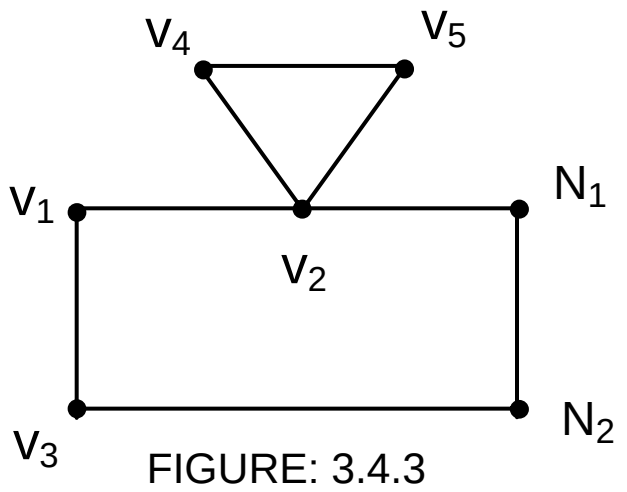

FIGURE: 3.4.3

were $N_1$ $N_2$ are indeterminate vertices and $v_1 \ldots v_5$ are vertices

***Example 3.4.3:***

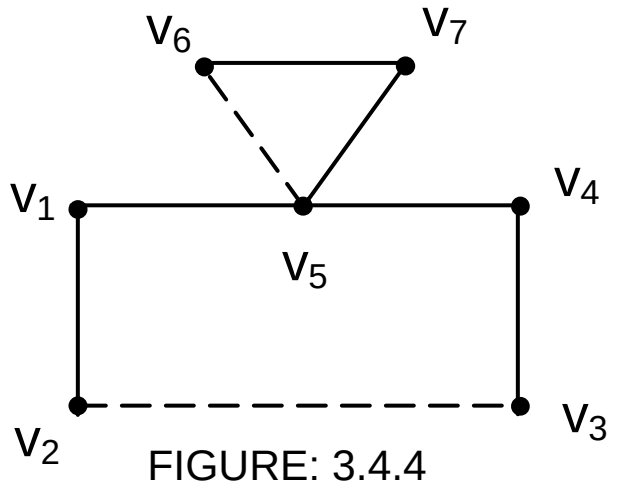

FIGURE: 3.4.4

The following is also a neutrosophic Eulerian graph with edges:

***Example 3.4.4:***

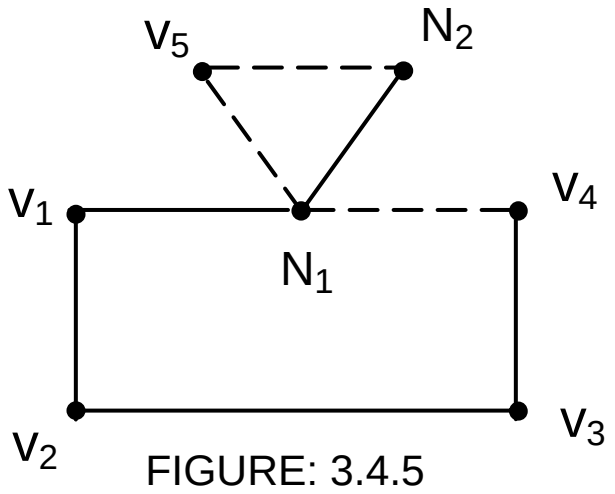

FIGURE: 3.4.5

The following is an example of a double of strong neutrosophic Eulerian graph:

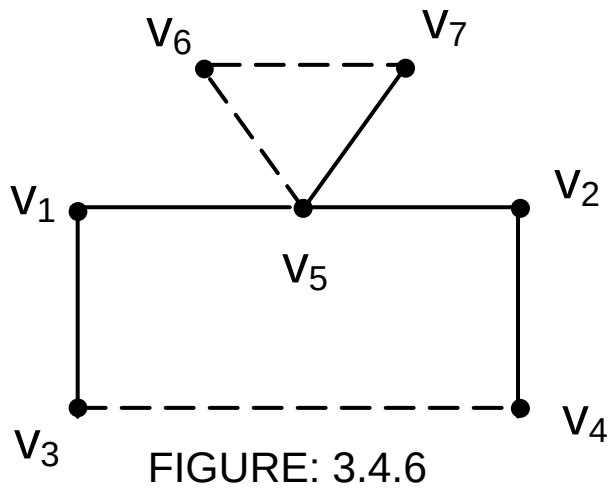

FIGURE: 3.4.6

where we have 3 indeterminate edges and no vertex is an indeterminate. Now we proceed on to give an example of a strong Eulerian graph or the double neutrosophic graph.

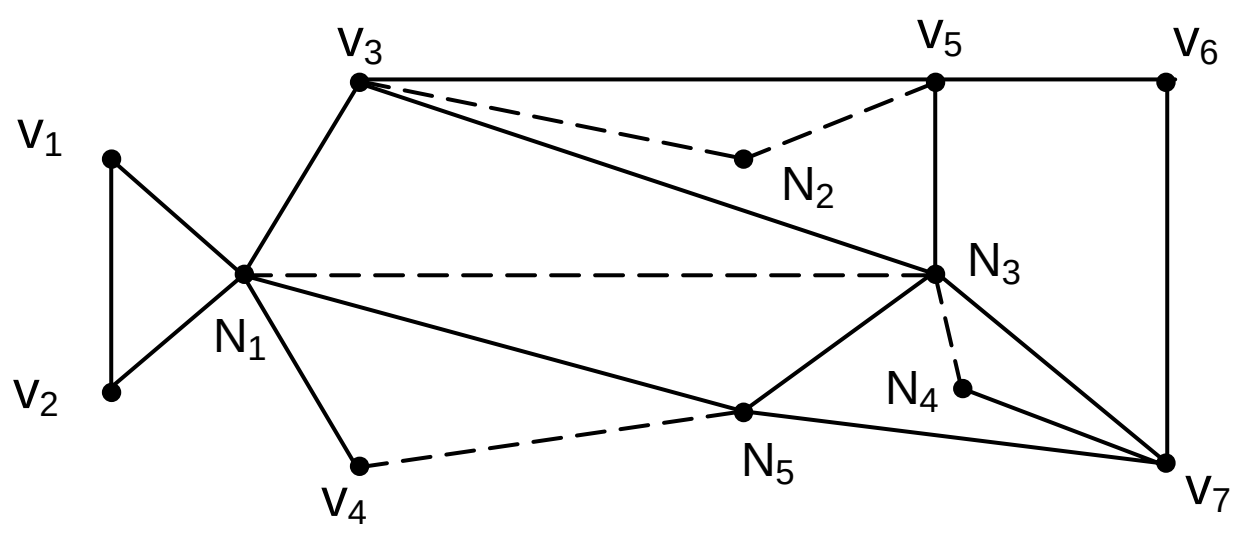

FIGURE: 3.4.7

This is a strong neutrosophic Eulerian graph with 12 vertices of which 7 vertices are real and 5 of the vertices are indeterminates.

We propose several interesting problems in this direction in the final chapter of this book. Now we proceed onto define the notion of neutrosophic Eulerian digraph.

**DEFINITION 3.4.5:** *A neutrosophic Eulerian trial in a neutrosophic digraph D is a closed spanning walk in which each arc of D occurs exactly once.*

A digraph is neutrosophic Eulerian if it has such a trial.

A result about indegree and outdegree is proposed as a problem in the last chapter.

## 3.5 Neutrosophic Graph Colourings

In this section we introduce and define the graph colourings of a neutrosophic graph. Throughout this section we may have neutrosophic graphs with indeterminate edges or indeterminate vertices.

**DEFINITION 3.5.1:** *The neutrosophic chromatic number $\chi_N$ (G) of a neutrosophic graph G is the minimum number of independent subsets that partition vertex set V of G (the real vertex set is N + V, N-number of neutrosophic vertices V-usual vertices).*

*Any such minimum partition is called the neutrosophic chromatic partition of $V_N$ (G) If N = 0 then the chromatic partition of G coincides with the neutrosophic chromatic partition.*

Let G be a neutrosophic graph; a vertex N-colouring of G is a map $f : V \rightarrow S$ where S is the set of distinct colours, it is proper if the adjacent real vertices of G receive distinct colours.

If an adjacent vertex is an indeterminant vertex then we do not demand the adjacent vertices to receive distinct colours. If the graph has no indeterminate vertices then the vertex N-colouring is the same as vertex colouring if $u, v \in E(G)$ then $f(u) \neq f(v)$ provided u and v are not indeterminant vertices otherwise u can be equal to v.

This $\chi_N(G)$ is the minimum cardinality of S for which there exists a proper vertex colouring of G by colours of S. (Here proper vertex colouring does not mean indeterminate vertex adjacent to any vertex, need to get distinct colours).

Clearly in any proper vertex coloring of G the vertices (the indeterminant vertices not included) that receive the same colour are independent. The vertices that receive a particular colour make up a color class (Color class will include indeterminate vertex also).

Thus in any chromatic partition of V (G) the parts of the partition constitute colour classes of only the real vertices (The indeterminate vertices are not included). Thus the neutrosophic chromatic number of a neutrosophic graph G is the minimum number of colors needed for a proper real vertex colouring of G, i.e., the indeterminate vertex can take colors very arbitrarily.

G is K-Chromatic if $X_N(G) = K$. Clearly $\chi_N(G) = \chi(G)$ for the colours of indeterminate vertex is not taken into account i.e., is calculated by jumping over the indeterminate vertex. Thus a neutrosophic graph G with N neutrosophic indeterminate vertex and V real vertices is said to be K colourable if the real vertices of G admits a proper vertex coloring using K-colors.

Now we proceed on to define the edge coloring of the neutrosophic graph. Here the graphs will have at least one indeterminate edge.

Suppose part of an indeterminate edge is one colour other part another colour one can start to define even the notion of indeterminate colourings or fuzzy colourings.

**DEFINITION 3.5.2:** *The edge coloring of a loop less neutrosophic graph (which has at least one indeterminate edge) is a function* $\pi : E(G) \rightarrow S$ *where S is the set of distinct colours, it is proper if no two adjacent real edges receive the same colour (So the following can happen, two adjacent indeterminate edges can* *receive same colour or two adjacent edges in which one is a real edge another an indeterminate edge can also receive same colour).*

The minimum K for which loop less neutrosophic graph has a proper K-edge coloring is called the neutrosophic edge chromatic number or N-chromatic index of G denoted by $\chi_N$' (G).

### 3.6 Petersen Graphs

Unlike Petersen graph, which is only one, we can have several neutrosophic Petersen graphs. We can define several neutrosophic Petersen graphs with in determinate nodes alones, with indeterminate edges alone with both indeterminate edges and nodes.

The study of finding the number of neutrosophic Petersen graphs is an interesting problem.

***Example 3.6.1:*** $\{N_1 – N_5\}$ are indeterminate vertices / nodes of G and $v_1$ … $v_5$ are real vertices of G. G in Figure 3.6.1 is a neutrosophic vertex Petersen graph.

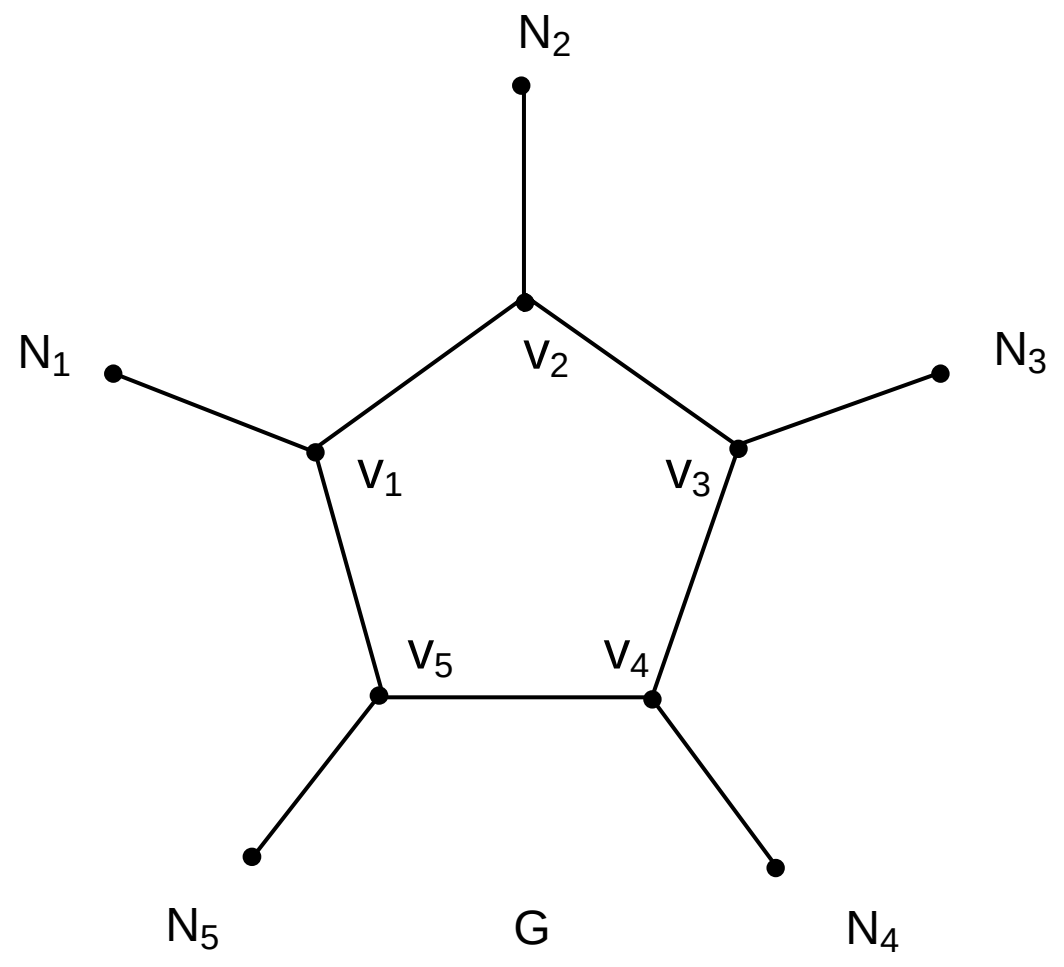

FIGURE: 3.6.1

***Example 3.6.2:***

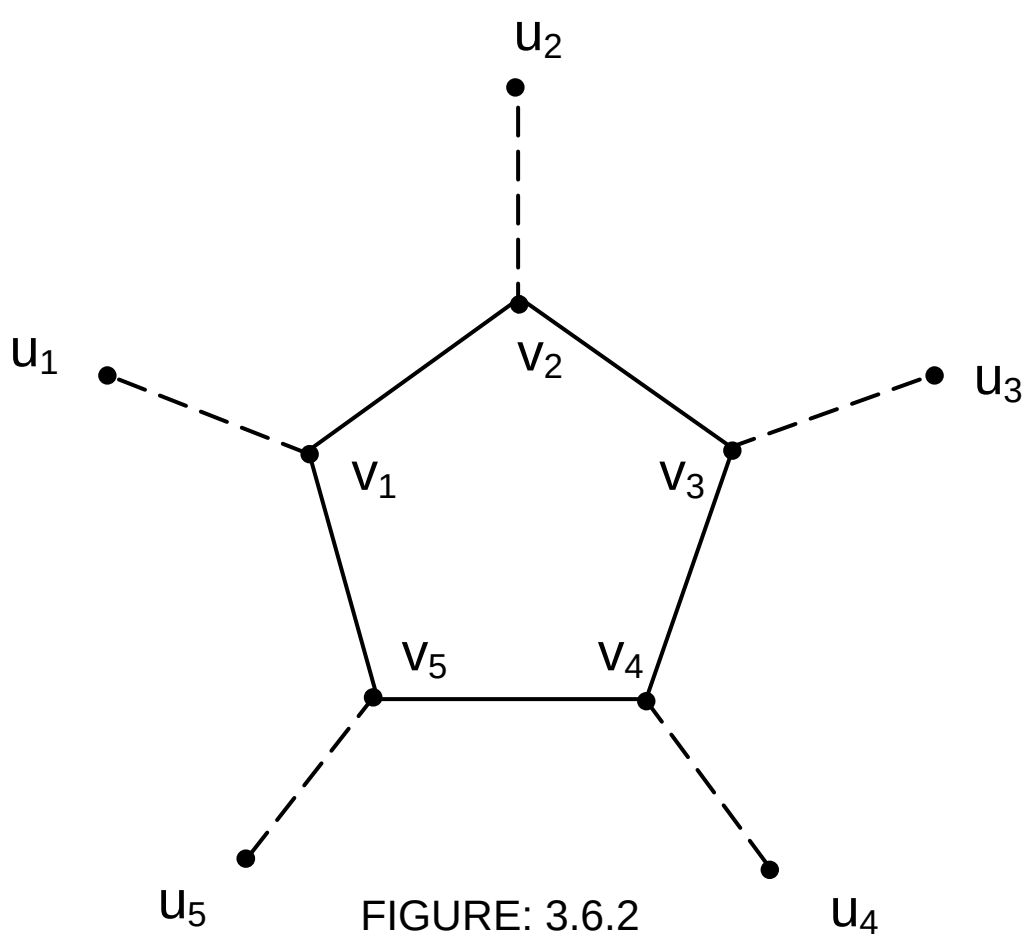

FIGURE: 3.6.2

This neutrosophic Petersen graph has 5 indeterminate edges all the vertices are real.

***Example 3.6.3:***

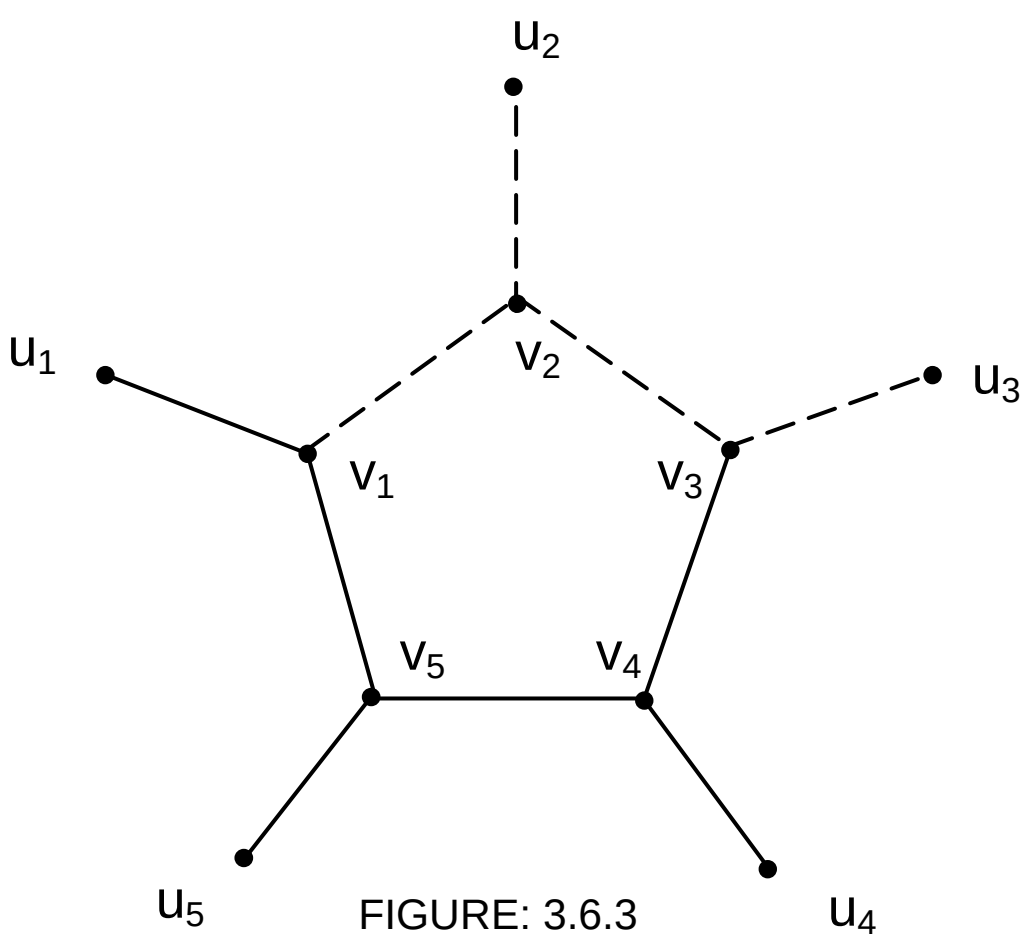

FIGURE: 3.6.3

This neutrosophic Petersen graph has 4 indeterminate edges all the 10 vertices are real

## 3.7 Application of neutrosophic graphs to neutrosophic models

The study of neutrosophic models like neutrosophic relational equations, neutrosophic cognitive maps, binary neutrosophic relations, composition of neutrosophic binary relations, neutrosophic sagittal diagram are very recent in 2004 [66]. But all these models make use of the neutrosophic graphs. For these models neutrosophic graphs describe explicitly the relations. Now we recall the applications of these neutrosophic graphs in the neutrosophic models.

### Binary neutrosophic Relation and their properties

In this section we introduce the notion of neutrosophic relational equations and fuzzy neutrosophic relational equations and analyze and apply them to real-world problems, which are abundant with the concept of indeterminacy. We also mention that most of the unsupervised data also involve at least to certain degrees the notion of indeterminacy.

Throughout this section by a neutrosophic matrix we mean a matrix whose entries are from the set N = [0, 1] $\cup$ *I* and by a fuzzy neutrosophic matrix we mean a matrix whose entries are from N' = [0, 1] $\cup$ {n*I* / n $\in$ (0,1]}.

Now we proceed on to define binary neutrosophic relations and binary neutrosophic fuzzy relation.

A binary neutrosophic relation $R_N(X, Y)$ may assign to each element of X two or more elements of Y or the indeterminate *I*. Some basic operations on functions such as the inverse and composition are applicable to binary relations as well. Given a neutrosophic relation $R_N(X, Y)$ its domain is a neutrosophic set on X $\cup$ *I* domain R whose membership function is defined by dom $R(x) = \max_{y \in X \cup I} R_N(x, y)$ for each x $\in$ X $\cup$ *I*.

That is each element of set X $\cup$ *I* belongs to the domain of R to the degree equal to the strength of its strongest relation to any member of set Y $\cup$ *I*. The degree may be an indeterminate *I* also. Thus this is one of the marked difference between the binary fuzzy relation and the binary neutrosophic relation. The range of $R_N(X,Y)$ is a neutrosophic relation on Y, ran R whose membership is defined by ran $R(y) = \max_{x \in X} R_N(x, y)$ for each y $\in$ Y, that is the strength of the strongest relation that each element of Y has to an element of X is equal to the degree of that element's membership in the range of R or it can be an indeterminate *I*.

The height of a neutrosophic relation $R_N(x, y)$ is a number h(R) or an indeterminate *I* defined by $h_N(R) = \max_{y \in Y \cup I} \max_{x \in X \cup I} R_N(x, y)$. That is $h_N(R)$ is the largest membership grade attained by any pair (x, y) in R or the indeterminate *I*.

A convenient representation of the neutrosophic binary relation $R_N(X, Y)$ are membership matrices $R = [\gamma_{xy}]$ where $\gamma_{xy} \in R_N(x, y)$. Another useful representation of a binary neutrosophic relation is a neutrosophic sagittal diagram. Each of the sets X, Y represented by a set of nodes in the diagram, nodes corresponding to one set are clearly distinguished from nodes representing the other set. Elements of X' × Y' with non-zero membership grades in $R_N(X, Y)$ are represented in the diagram by lines connecting the respective nodes. These lines are labeled with the values of the membership grades.

An example of the neutrosophic sagittal diagram is a binary neutrosophic relation $R_N(X, Y)$ together with the membership neutrosophic matrix which is given below:

|  | $y_1$ | $y_2$ | $y_3$ | $y_4$ |
|---|---|---|---|---|
| $x_1$ | *I* | 0 | 0 | 0.5 |
| $x_2$ | 0.3 | 0 | 0.4 | 0 |
| $x_3$ | 1 | 0 | 0 | 0.2 |
| $x_4$ | 0 | *I* | 0 | 0 |
| $x_5$ | 0 | 0 | 0.5 | 0.7 |

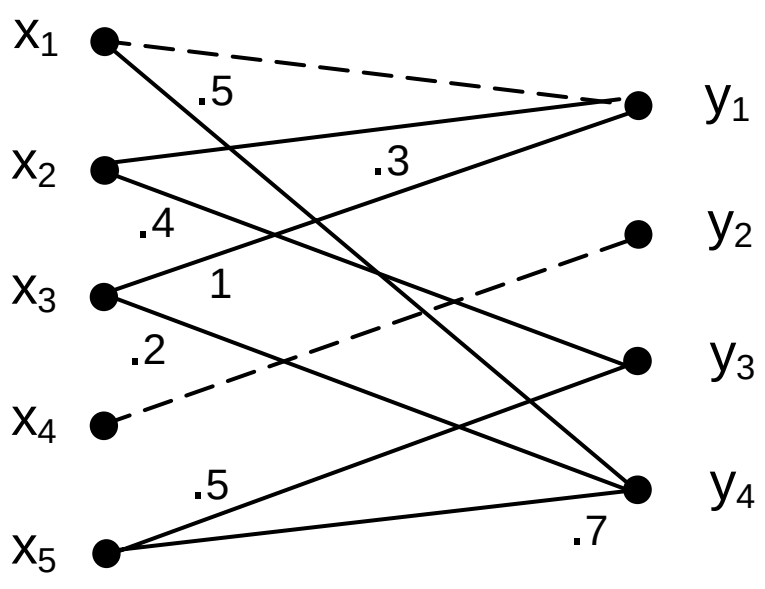

FIGURE: 3.7.1

The inverse of a neutrosophic relation $R_N(X, Y)$ which is denoted by $R^{-1}$ (Y, X) is a relation on $Y \times X$ defined by $R^{-1}(y, x) = R(x, y)$ for all $x \in X$ and all $y \in Y$. A neutrosophic membership matrix $R^{-1} = [r^{-1}_{yx}]$ representing $R^{-1}_N(Y, X)$ is the transpose of the matrix R for $R_N(X, Y)$ which means that the rows of $R^{-1}$ equal the columns of R and the columns of $R^{-1}$ equal rows of R. Clearly $(R^{-1})^{-1} = R$ for any binary neutrosophic relation.

Consider any two binary neutrosophic relation $P_N(X, Y)$ and $Q_N(Y, Z)$ with a common set Y. The standard composition of these relations which is denoted by $P_N(X, Y) \bullet Q_N(Y, Z)$ produces a binary neutrosophic relation $R_N(X, Z)$ on $X \times Z$ defined by $R_N(x, z) = [P \bullet Q]_N(x, z) = \max_{y \in Y} \min[P_N(x, y), Q_N(x, y)]$ for all $x \in X$ and all $z \in Z$.

This composition which is based on the standard $t_N$-norm and $t_N$-co-norm, is often referred to as the max-min composition. It can be easily verified that even in the case of binary neutrosophic relations $[P_N(X, Y) \bullet Q_N(Y, Z)]^{-1} = Q^{-1}_N(Z, Y) \bullet P^{-1}_N(Y, X)$. $[P_N(X, Y) \bullet Q_N(Y, Z)] \bullet R_N(Z, W) = P_N(X, Y) \bullet [Q_N(Y, Z) \bullet R_N(Z, W)]$, that is, the standard (or max-min) composition is associative and its inverse is equal to the reverse composition of the inverse relation. However, the standard composition is not commutative, because $Q_N(Y, Z) \bullet P_N(X, Y)$ is not well defined when $X \neq Z$. Even if $X = Z$ and $Q_N (Y, Z) \circ P_N (X, Y)$ are well defined still we can have $P_N (X, Y) \circ Q (Y, Z) \neq Q (Y, Z) \circ P (X, Y)$.

Compositions of binary neutrosophic relation can be performed conveniently in terms of membership matrices of the relations. Let $P = [p_{ik}]$, $Q = [q_{kj}]$ and $R = [r_{ij}]$ be membership matrices of binary relations such that $R = P \circ Q$. We write this using matrix notation

$$[r_{ij}] = [p_{ik}] \text{ o } [q_{kj}]$$

where $r_{ij} = \max_{k} \min (p_{ik}, q_{kj})$.

A similar operation on two binary relations, which differs from the composition in that it yields triples instead of pairs, is known as the relational join. For neutrosophic relation $P_N (X, Y)$ and $Q_N (Y, Z)$ the relational join P * Q corresponding to the neutrosophic standard max-min composition is a ternary relation $R_N (X, Y, Z)$ defined by $R_N (x, y, z) = [P * Q]_N (x, y, z) = \min [P_N (x, y), Q_N (y, z)]$ for each $x \in X$, $y \in Y$ and $z \in Z$.

This is illustrated by the following Figure 3.7.2.

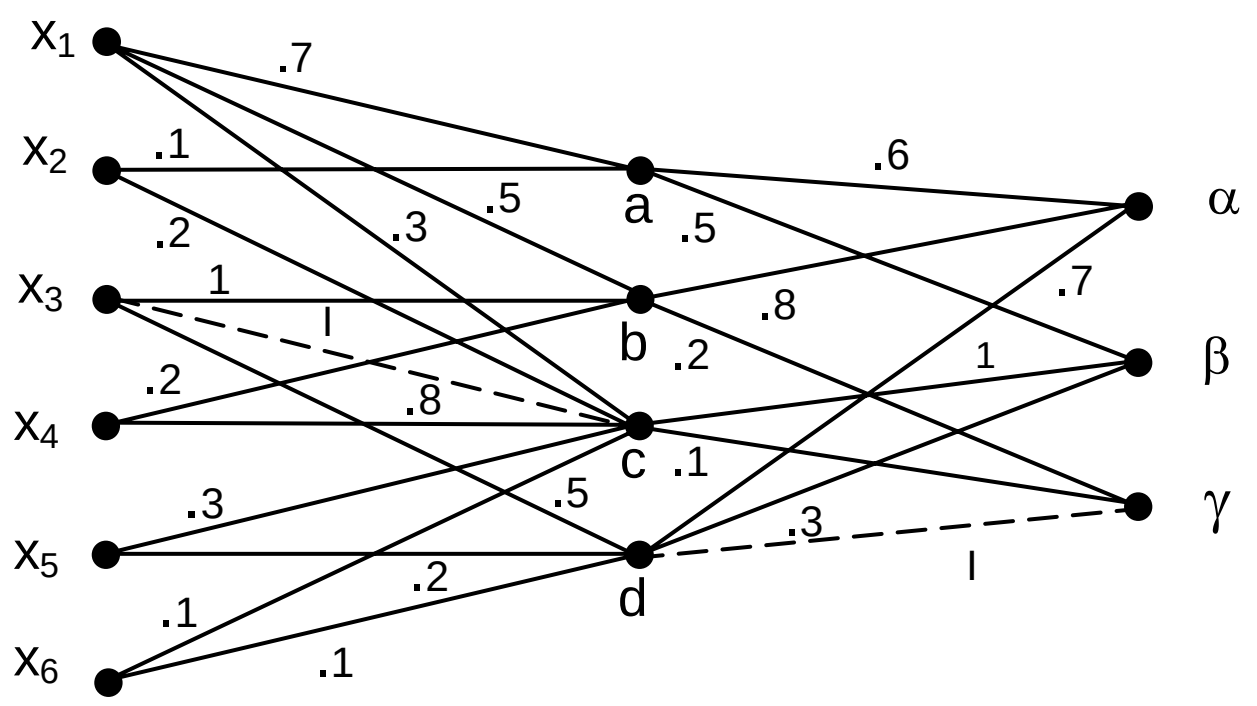

FIGURE: 3.7.2

In addition to defining a neutrosophic binary relation there exists between two different sets, it is also possible to define neutrosophic binary relation among the elements of a single set X.

A neutrosophic binary relation of this type is denoted by $R_N(X, X)$ or $R_N(X^2)$ and is a subset of $X \times X = X^2$.

These relations are often referred to as neutrosophic directed graphs or neutrosophic digraphs.

Neutrosophic binary relations $R_N$ (X, X) can be expressed by the same forms as general neutrosophic binary relations. However they can be conveniently expressed in terms of simple diagrams with the following properties:

i. Each element of the set X is represented by a single node in the diagram.
ii. Directed connections between nodes indicate pairs of elements of X for which the grade of membership in R is non zero or indeterminate.
iii. Each connection in the diagram is labeled by the actual membership grade of the corresponding pair in R or in indeterminacy of the relationship between those pairs.

The neutrosophic membership matrix and the neutrosophic sagittal diagram is as follows for any set X = {a, b, c, d, e}.

$$
\begin{array}{c} \\ a \\ b \\ c \\ d \\ e \end{array}
\begin{array}{c} \begin{array}{ccccc} a & b & c & d & e \end{array} \\
\begin{bmatrix} 0 & I & .3 & .2 & 0 \\ 1 & 0 & I & 0 & .3 \\ I & .2 & 0 & 0 & 0 \\ 0 & .6 & 0 & .3 & I \\ 0 & 0 & 0 & I & .2 \end{bmatrix} \end{array}
$$

Neutrosophic membership matrix for x is given above and the neutrosophic sagittal diagram is given below:

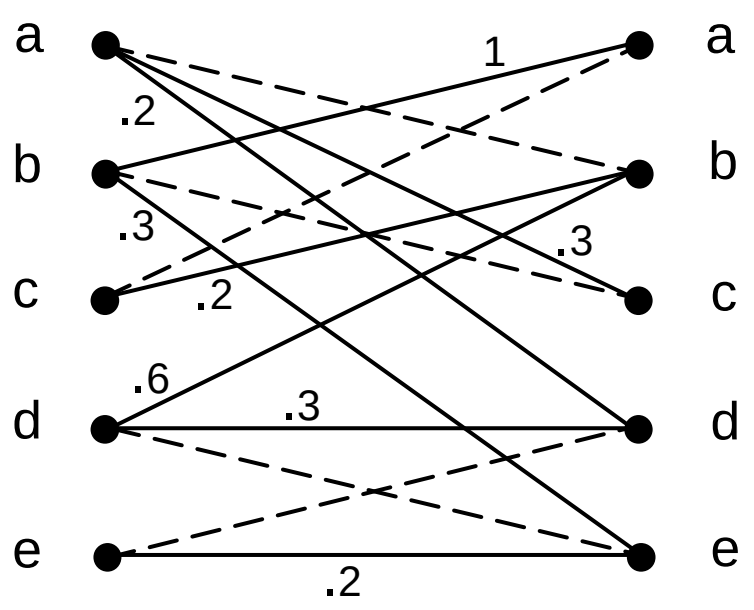

FIGURE: 3.7.3

Neutrosophic diagram or graph is left for the reader as an exercise.

The notion of reflexivity, symmetry and transitivity can be extended for neutrosophic relations $R_N$ (X, Y) by defining them in terms of the membership functions or indeterminacy relation.

*Thus $R_N$ (X, X) is reflexive if and only if $R_N (x, x) = 1$ for all $x \in X$.* If *this is not the case for some $x \in X$ the relation is irreflexive.* *A weaker form of reflexivity, if for no x in X; $R_N(x, x) = 1$ then we call the relation to be anti-reflexive referred to as $\in$-reflexivity, is sometimes defined by requiring that*

$$R_N (x, x) \geq \in \text{ where } 0 < \in < 1.$$

*A neutrosophic relation is symmetric if and only if*

$$R_N (x, y) = R_N (y, x) \text{ for all } x, y, \in X.$$

*Whenever this relation is not true for some $x, y \in X$ the relation is called asymmetric. Furthermore when $R_N (x, y) > 0$ and $R_N (y, x) > 0$ implies that $x = y$ for all $x, y \in X$ the relation $R_N(X, Y)$ is called anti-symmetric.*

*A fuzzy relation $R_N$ (X, X) is transitive (or more specifically max-min transitive) if*

$$R_N(x, z) \geq \max_{y \in Y} \min [R_N(x, y), R_N(y, z)]$$

*is satisfied for each pair (x, z) $\in X^2$. A relation failing to satisfy the above inequality for some members of X is called non-transitive and if $R_N(x, x) < \max_{y \in Y} \min [R_N(x, y), R_N(y, z)]$ for all (x, x) $\in X^2$, then the relation is called anti-transitive.*

*Given a relation $R_N$(X, X) its transitive closure $\overline{R}_{NT}$ (X, X) can be analyzed in the following way.*

The transitive closure on a crisp relation $R_N$ (X, X) is defined as the relation that is transitive, contains

$$R_N(X, X) < \max_{y \in Y} \min [R_N(x, y)\, R_N(y, z)]$$

for all (x, x) $\in X^2$, then the relation is called anti-transitive. Given a relation $R_N$ (X, X) its transitive closure $\overline{R}_{NT}$ (X, X) can be analyzed in the following way.

The transitive closure on a crisp relation $R_N$ (X, X) is defined as the relation that is transitive, contains $R_N$ and has the fewest possible members. For neutrosophic relations the last requirement is generalized such that the elements of transitive closure have the smallest possible membership grades, that still allow the first two requirements to be met.

Given a relation $R_N$ (X, X) its transitive closure $\overline{R}_{NT}$ (X, X) can be determined by a simple algorithm.

**DEFINITION 3.7.1:** *A Neutrosophic Cognitive Map (NCM) is a neutrosophic directed graph with concepts like policies, events etc. as nodes and causalities or indeterminates as edges. It represents the causal relationship between concepts.*

*Let $C_1, C_2, \ldots, C_n$ denote n nodes, further we assume each node is a neutrosophic vector from neutrosophic vector space V. So a node $C_i$ will be represented by $(x_1, \ldots, x_n)$ where $x_k$'s are zero or one or I (I is the indeterminate) and $x_k = 1$ means that the node $C_k$ is in the on state and $x_k = 0$ means the node is in the off state and $x_k = I$ means the nodes state is an indeterminate at that time or in that situation.*

*Let $C_i$ and $C_j$ denote the two nodes of the NCM. The directed edge from $C_i$ to $C_j$ denotes the causality of $C_i$ on $C_j$ called connections. Every edge in the NCM is weighted with a number in the set {-1, 0, 1, I}. Let $e_{ij}$ be the weight of the directed edge $C_iC_j$, $e_{ij} \in \{-1, 0, 1, I\}$. $e_{ij} = 0$ if $C_i$ does not have any effect on $C_j$, $e_{ij} = 1$ if increase (or decrease) in $C_i$ causes increase (or decreases) in $C_j$, $e_{ij} = -1$ if increase (or decrease) in $C_i$ causes decrease (or increase) in $C_j$ . $e_{ij} = I$ if the relation or effect of $C_i$ on $C_j$ is an indeterminate.*

**DEFINITION 3.7.2:** *NCMs with edge weight from {-1, 0, 1, I} are called simple NCMs.*

**DEFINITION 3.7.3**: *Let $C_1, C_2, \ldots, C_n$ be nodes of a NCM. Let the neutrosophic matrix N(E) be defined as N(E) = ($e_{ij}$) where $e_{ij}$ is the weight of the directed edge $C_i C_j$, where $e_{ij} \in \{0, 1, -1, I\}$. N(E) is called the neutrosophic adjacency matrix of the NCM.*

**DEFINITION 3.7.4**: *Let $C_1, C_2, \ldots, C_n$ be the nodes of the NCM. Let $A = (a_1, a_2, \ldots, a_n)$ where $a_i \in \{0, 1, I\}$. A is called the instantaneous state neutrosophic vector and it denotes the on – off – indeterminate state position of the node at an instant*

$a_i$ = 0 if $a_i$ is off (no effect)
$a_i$ = 1 if $a_i$ is on (has effect)
$a_i$ = I if $a_i$ is indeterminate(effect cannot be determined)

*for i = 1, 2,…, n.*

***Example 3.7.1:*** The child labor problem prevalent in India is modeled in this example using NCMs.

Let us consider the child labor problem with the following conceptual nodes:

$C_1$ - Child Labor
$C_2$ - Political Leaders
$C_3$ - Good Teachers
$C_4$ - Poverty
$C_5$ - Industrialists
$C_6$ - Public practicing/encouraging Child Labor
$C_7$ - Good Non-Governmental Organizations (NGOs)

$C_1$ - Child labor, it includes all types of labor of children below 14 years which include domestic workers, rag pickers, working in restaurants / hotels, bars etc. (It can be part time or fulltime).

$C_2$ - We include political leaders with the following motivation: Children are not vote banks so political leaders are not directly concerned with child labor but they indirectly help in the flourishing of it as industrialists who utilize child laborers or cheap labor are the decision makers for the winning or losing of the political leaders. Also industrialists financially control political interests. So we are forced to include political leaders as a node in this problem.

$C_3$ - Teachers are taken as a node because mainly school dropouts or children who have never attended the school are child laborers. So if the motivation by the teacher is very good, there would be less school dropouts and therefore there would be a decrease in child laborers.

$C_4$ - Poverty which is the most responsible reason for child labor.

$C_5$ - Industrialists – when we say industrialists we include one and all starting from a match factory or beedi factory, bars, hotels etc.

$C_6$ - Public who promote child labor as domestic servants, sweepers etc.

$C_7$ - We qualify the NGOs as good for some NGOs may not take up the issue fearing the rich and the powerful. Here "good NGOs" means NGOs who try to stop or prevent child labor.

Now we give the directed graph as well as the neutrosophic graph of two experts in the following Figures 3.7.4 and 3.7.5:

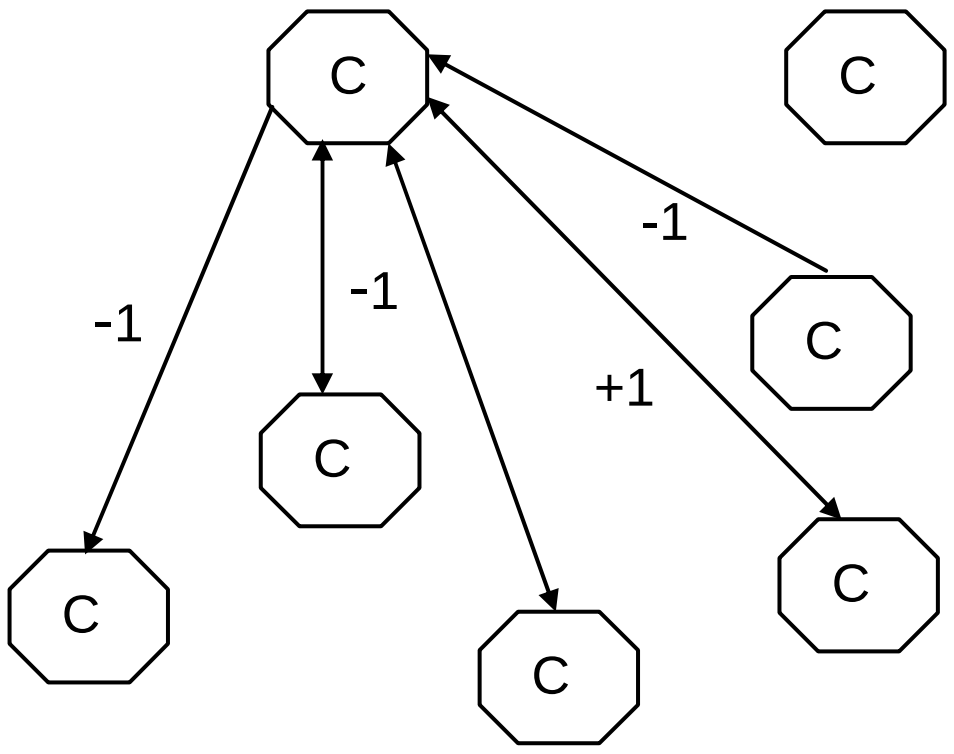

FIGURE: 3.7.4

Figure 3.7.4 gives the directed graph with $C_1$, $C_2$, …, $C_7$ as nodes and Figure 3.7.5 gives the neutrosophic directed graph with the same nodes.

The connection matrix E related to the graph in Figure 3.7.4 is given below:

$$E = \begin{bmatrix} 0 & 0 & 0 & 1 & 1 & 1 & -1 \\ 0 & 0 & 0 & 0 & 0 & 0 & 0 \\ -1 & 0 & 0 & 0 & 0 & 0 & 0 \\ 1 & 0 & 0 & 0 & 0 & 0 & 0 \\ 1 & 0 & 0 & 0 & 0 & 0 & 0 \\ 0 & 0 & 0 & 0 & 0 & 0 & 0 \\ -1 & 0 & 0 & 0 & 0 & 0 & 0 \end{bmatrix}.$$

According to this expert no connection however exists between political leaders and industrialists.

Now we reformulate a different format of the questionnaire where we permit the expert to give answers like the relation between certain nodes is indeterminable or not known. Now based on the expert's opinion also about the notion of indeterminacy we obtain the following neutrosophic directed graph:

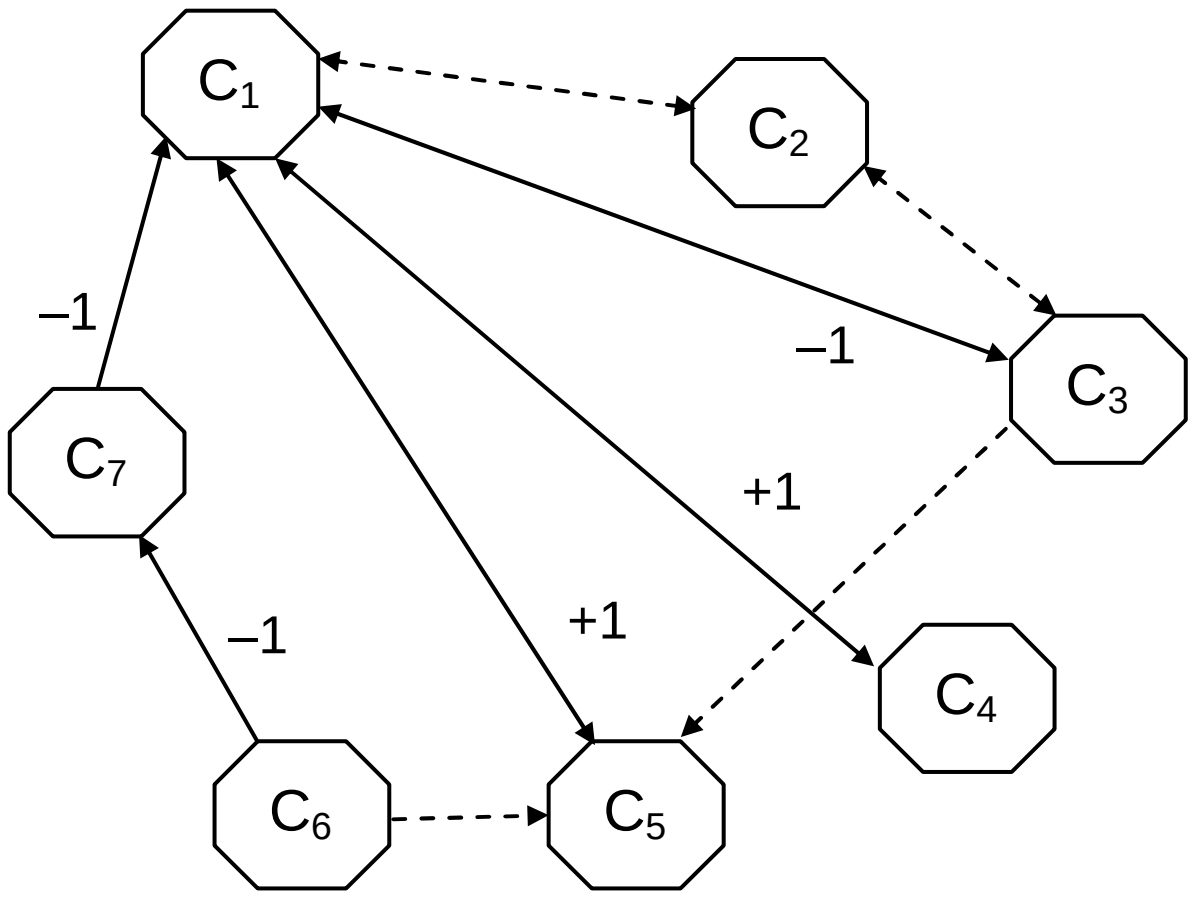

FIGURE: 3.7.5

The corresponding neutrosophic adjacency matrix N(E) related to the neutrosophic directed graph (Figure 3.7.5) is given below:

$$N(E) = \begin{bmatrix} 0 & I & -1 & 1 & 1 & 0 & 0 \\ I & 0 & I & 0 & 0 & 0 & 0 \\ -1 & I & 0 & 0 & I & 0 & 0 \\ 1 & 0 & 0 & 0 & 0 & 0 & 0 \\ 1 & 0 & 0 & 0 & 0 & 0 & 0 \\ 0 & 0 & 0 & 0 & I & 0 & -1 \\ -1 & 0 & 0 & 0 & 0 & 0 & 0 \end{bmatrix}$$

Suppose we take the state vector $A_1$ = (1 0 0 0 0 0 0). We will see the effect of $A_1$ on E and on N(E).

$A_1E$ = (0 0 0 1 1 1 –1) → (1 0 0 1 1 1 0) = $A_2$
$A_2$ E = (2 0 0 1 1 1 0) → (1 0 0 1 1 1 0) = $A_3 = A_2$.

Thus child labor flourishes with parents' poverty and industrialists' action. Public practicing child labor also flourish but good NGOs are absent in such a scenario. The state vector gives the fixed point.

Now we find the effect of $A_1$ = (1 0 0 0 0 0 0) on N(E).

$A_1$ N(E) = (0 I –1 1 1 0 0) → (1 I 0 1 1 0 0) = $A_2$

$A_2$ N(E) = (I + 2, I, –1+ I, 1 1 0 0) → (1 I 0 1 1 0 0) = $A_2$

Thus $A_2$ = (1 I 0 1 1 0 0), according to this expert the increase or the on state of child labor certainly increases with the poverty of parents and other factors are indeterminate to him. This mainly gives the indeterminates relating to political leaders and teachers in the neutrosophic cognitive model and the parents poverty and Industrialist become to on state.

However, the results by FCM give as if there is no effect by teachers and politicians for the increase in child labor. Actually the increase in school dropout increases the child labor hence certainly the role of teachers play a part. At least if it is termed as an indeterminate one would think or reflect about their (teachers) effect on child labor.

Also the node the role played by political leaders has a major part; for if the political leaders were stern about stopping the child labor, certainly it cannot flourish in the society. They are ignored for two reasons: First, if children were vote banks certainly their position would be better. The second reason is, industrialists who practice child labor, are a main source of help to politicians, and their victory/defeat depends on their (financial) support so the causes for politicians ignoring child labor is two-fold.

Now we seek the opinion of another expert who is first asked to give a FCM model and then a provocative questionnaire discussing about the indeterminacy of relation between nodes is suggested and he finally gives a neutrosophic version of his ideas.

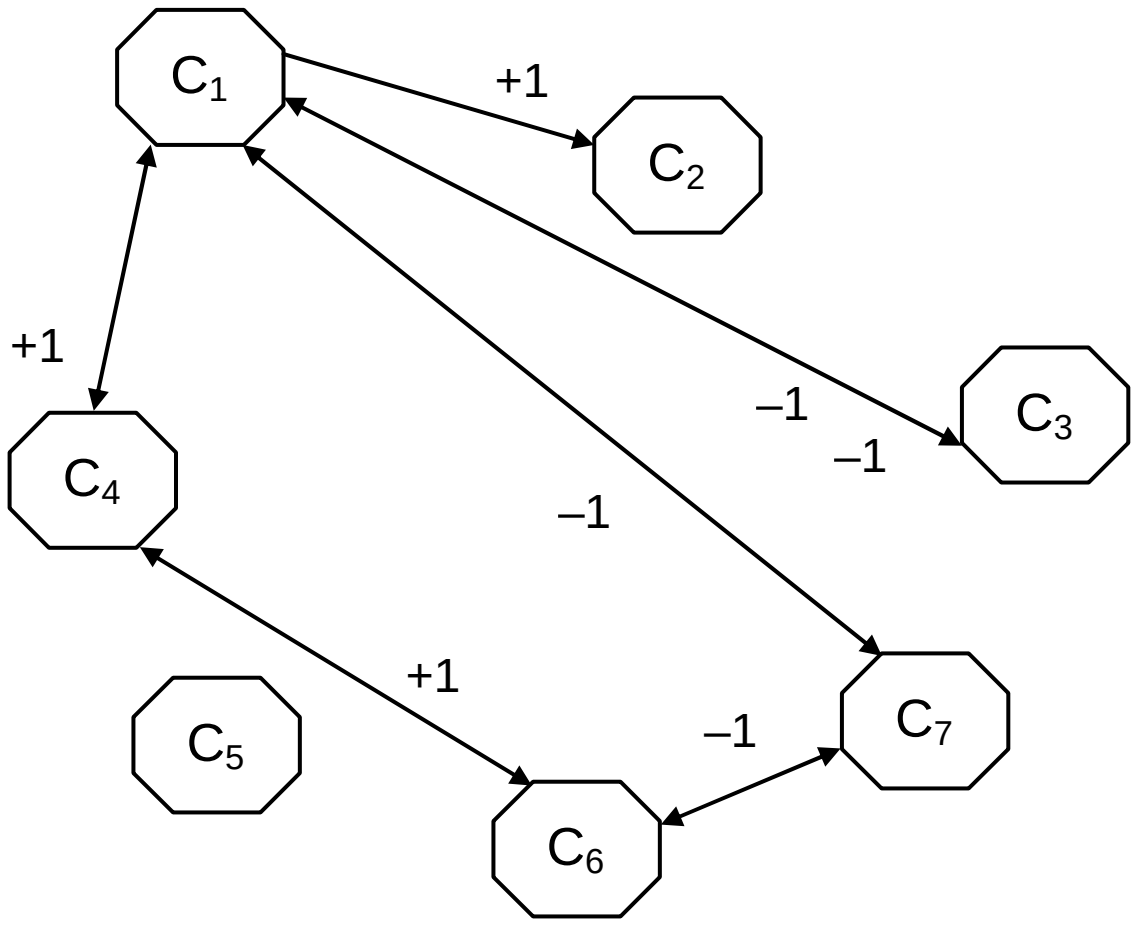

FIGURE: 3.7.6

Figure 3.7.6 is the directed graph of the expert. The related connection matrix $E_1$ is as follows:

$$E_1 = \begin{bmatrix} 0 & 1 & -1 & 1 & 0 & 0 & -1 \\ 0 & 0 & 0 & 0 & 0 & 0 & 0 \\ -1 & 0 & 0 & 0 & 0 & 0 & 0 \\ 1 & 0 & 0 & 0 & 0 & 1 & 0 \\ 0 & 0 & 0 & 0 & 0 & 0 & 0 \\ 0 & 0 & 0 & 1 & 0 & 0 & -1 \\ -1 & 0 & 0 & 0 & 0 & -1 & 0 \end{bmatrix}.$$

Take $A_1 = (1\ 0\ 0\ 0\ 0\ 0\ 0)$ the effect of $A_1$ on the system $E_1$ is

$A_1E_1 = (0\ 1\ –1\ 1\ 0\ 0\ –1) \to (1\ 1\ 0\ 1\ 0\ 0\ 0) = A_2$
$A_2E_2 = (1\ 1\ –1\ 1\ 0\ 1\ –1) \to (1\ 1\ 0\ 1\ 0\ 1\ 0) = A_3$
$A_3E_2 = (1\ 1\ –1\ 2\ 0\ 1\ –2) \to (1\ 1\ 0\ 1\ 0\ 1\ 0) = A_4 = A_3.$

Thus according to this expert child labor has direct effect on political leaders, no effect on good teachers, effect on poverty and industrialists and no-effect on the public who encourage child labor; and good NGOs.

The same person was now put with the neutrosophic questions i.e. terms like "can you find any relation between the nodes or are you not in apposition to decide any relation between two nodes and so on"; so that a idea of indeterminacy is introduced to them.

Now the neutrosophic directed graph is drawn using this experts opinion given in Figure 3.7.7.

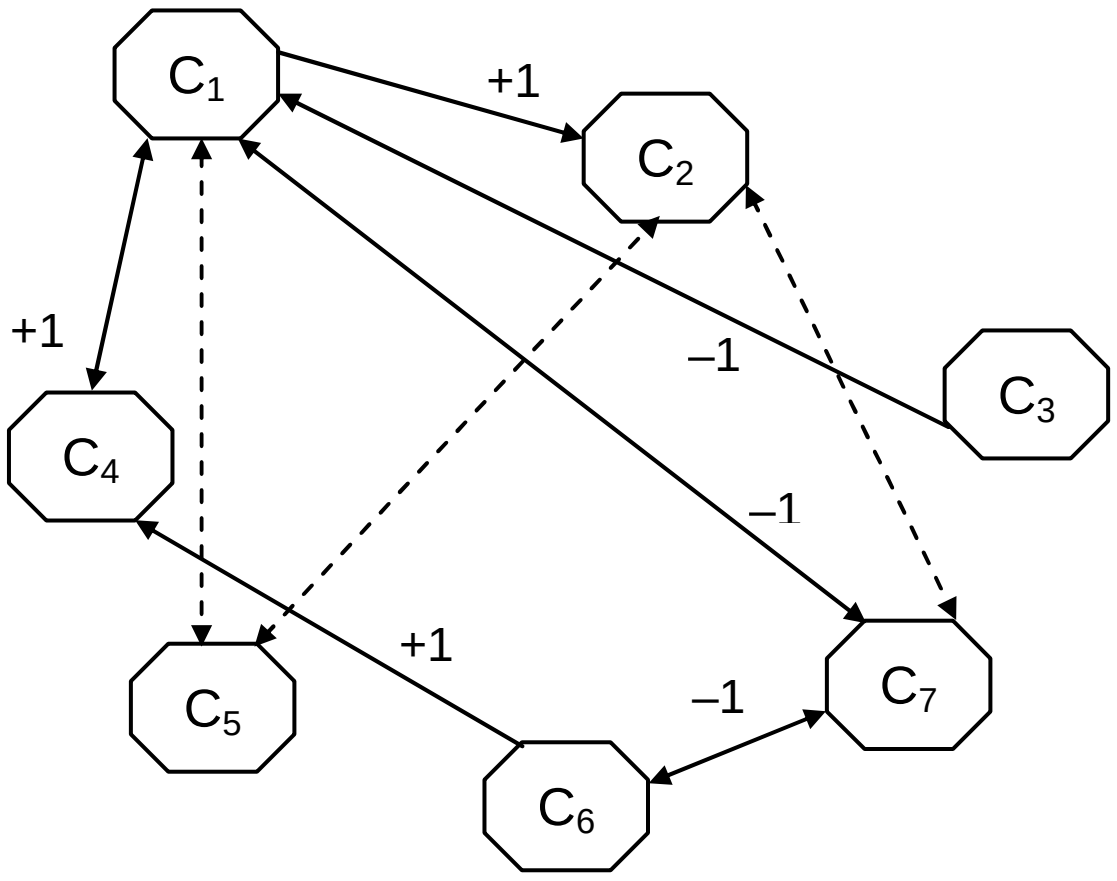

FIGURE: 3.7.7

The corresponding neutrosophic connection matrix $N(E_1)$ is as follows:

$$N(E_1) = \begin{bmatrix} 0 & 1 & -1 & 1 & I & 0 & -1 \\ 0 & 0 & 0 & 0 & I & 0 & I \\ -1 & 0 & 0 & 0 & 0 & 0 & 0 \\ 1 & 0 & 0 & 0 & 1 & 0 & 0 \\ I & I & 0 & 0 & 0 & 0 & 0 \\ 0 & 0 & 0 & 1 & 0 & 0 & -1 \\ -1 & 0 & 0 & 0 & 0 & -1 & I \end{bmatrix}.$$

Suppose $A_1 = (1, 0, 0, 0, 0, 0, 0)$ is the state vector whose effect on the neutrosophic system $N(E_1)$ is to be considered.

$A_1N(E_1) \;=\; (0\ 1\ –1\ 1\ I\ 0\ –1) \;\to\; (1\ 1\ 0\ 1\ I\ 0\ 0) \;=\; A_2$

$A_2N(E_1) \;=\; (1+I, 1+I, -1, 1, 2I+1, 0\ –1+I)$
$\to\; (1\ 1\ 0\ 1\ 1\ 0\ 0) \;=\; A_3$

$A_3N(E_1) \;=\; (1+I, 1+I, -1, 2\ I+1\ 0\ –1 + I)$
$\to\; (1\ 1\ 0\ 1\ 1\ 0\ 0\ ) \;=\; A_4.$

We see $A_2 = A_3$.

But according to the NCM when the conceptual node child labor is on it implies that the cause of it is political leaders, poverty and industrialists participation by employing children as laborers.

The reader is expected to compare the graphs of the NCM with FCMs for the same problem which is dealt earlier in this book as we have now indicated how a NCM works.

***Example 3.7.2:*** Application of NCM to study the Hacking of e-mail by students. One of the major problems of today's world of information technology that is faced by one and all is; How safe are the messages that are sent by e-mail? Is there enough privacy? For if a letter is sent by post one can by certain say that it cannot be read by any other person, other than the receiver. Even tapping or listening (over hearing) of phone calls from an alternate location / extension is only a very uncommon problem.

However compared to these modes of communication even though e-mail guarantees a lot of privacy it is a highly common practice to hack e-mail. Hacking is legally a cyber crime but is also one of the crimes that does not leave any trace. Hacking of another persons e-mail account can be carried out for a variety of purposes to study the factors, which are root-causes of such crimes we use NCM to analyze them. The following nodes are taken as the conceptual nodes.

$C_1$ - Curiosity
$C_2$ - Professional rivalry
$C_3$ - Jealousy/ enmity
$C_4$ - Sexual satisfaction
$C_5$ - Fun/pastime
$C_6$ - To satisfy ego
$C_7$ - Women students
$C_8$ - Breach of trust.

However more number of conceptual nodes can be added as felt by the expert or the investigator. The neutrosophic directed graph as given by an expert is given in Figure 3.7.8.

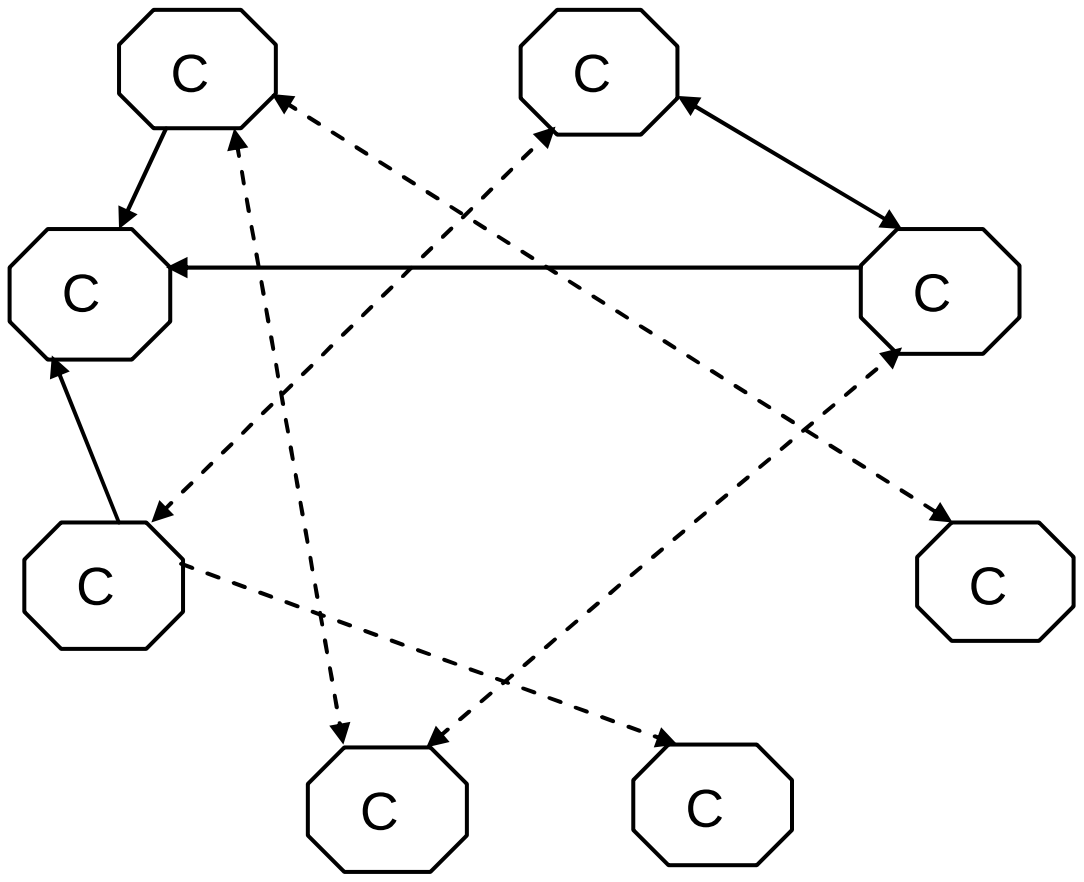

FIGURE: 3.7.8

The corresponding neutrosophic connection matrix N(E) is as follows:

$$N(E) = \begin{bmatrix} 0 & 0 & 0 & I & 0 & I & 0 & 1 \\ 0 & 0 & 1 & 0 & 0 & 0 & I & 0 \\ 0 & 1 & 0 & 0 & 0 & I & 0 & 1 \\ I & 0 & 0 & 0 & 0 & 0 & 0 & 0 \\ 0 & 0 & 0 & 0 & 0 & 0 & 0 & 0 \\ I & 0 & I & 0 & 0 & 0 & 0 & 0 \\ 0 & I & 0 & 0 & I & 0 & 0 & 1 \\ 0 & 0 & 0 & 0 & 0 & 0 & 0 & 0 \end{bmatrix}.$$

Suppose we take the instantaneous state vector $A_1$ = (0 0 0 0 0 0 1 0), women students node is in the on state then the effect of $A_1$ on the neutrosophic system N(E) is given by

$$\begin{aligned}
A_1N(E) &= (0\ I\ 0\ 0\ I\ 0\ 0\ 1) \rightarrow (0\ I\ 0\ 0\ I\ 0\ 1\ 1) = A_2 \\
A_2N(E) &= (0\ I\ I\ 0\ I\ 0\ I\ 1) \rightarrow (0\ I\ I\ 0\ I\ 0\ 1\ 1) = A_3 \\
A_3N(E) &= (0\ I\ I\ 0\ I\ I\ I\ 1{+}I) \rightarrow (0\ I, I, 0, I, I, 1, 1) = A_4 \\
A_4N(E) &= (I, I\ I\ 0\ I\ I\ I\ I\ {+}1) \rightarrow (I\ I\ I\ 0\ I\ I\ 1\ 1) = A_5 \\
A_5N(E) &= (I\ I\ I\ 0\ I\ I\ I\ 1) \rightarrow (I\ I\ I\ 0\ I\ I\ 1\ 1) = A_6 = A_5.
\end{aligned}$$

So in case the node "women students" is in the on state node we see curiosity is an indeterminate, professional rivalry is an indeterminate, jealousy/ enmity is an indeterminate sexual satisfaction is in the off state, fun/ pastime is an indeterminate, to satisfy ego is an indeterminate and breach of trust is in the on state whereas if the '$I$ 's are removed and N(E) is used as a usual FCM matrix then the effect $B_1$ = (0 0 0 0 0 0 1 0) in the on state when passed through the system we get B = (0 0 0 0 1 0 1 1) implies to satisfy ego becomes the on state and Breach of trust is in the on state. Thus we see other sates are in the off state.

The reader is expected to work with other coordinates and compare with FCMs which is got by replacing all $I$ 's in the neutrosophic connection matrix N(E) by 0.

Several other examples can be shown using the method of NCM. We give some more application of NCM. The study of application of FCMs are given in chapter 2 of this book. Here in this section we apply NCMs only to some of the illustrations mentioned in that section.

***Example 3.7.3:*** Here Analysis of Strategic Planning Simulation based on NCMs Knowledge and Differential Game is given. FCM has been used in the study of differential game we use same map of FCM but after discussing with an expert convert it into an NCM by adjoining the edges which are indeterminate, and this is mainly carried out for easy comparison.

Now according to this expert, competitiveness and market demand is an indeterminate. Also sales price and economic condition is an indeterminate. Also according to him the productivity and market share is an indeterminate whether a relation exists directly cannot be said but he is not able to state that there is no relation between these concepts so he says let it be an indeterminate.

Also according to him quality control and market share is an indeterminate. Thus on the whole the market share is an FCM with a lot of indeterminacy so is best fit with an NCM model.

Figure 3.7.9 gives the NCM and obtain analysis and conclusions using NCMs and compare it with FCM. Thus obtain the initial version of NCM matrix and refined version of NCM matrix, also give the corresponding comment. Study the factor of indeterminacy and prove the result is nearer to truth for finding solutions to the market share problem. Compare FCM and NCM in the case of market share problem.

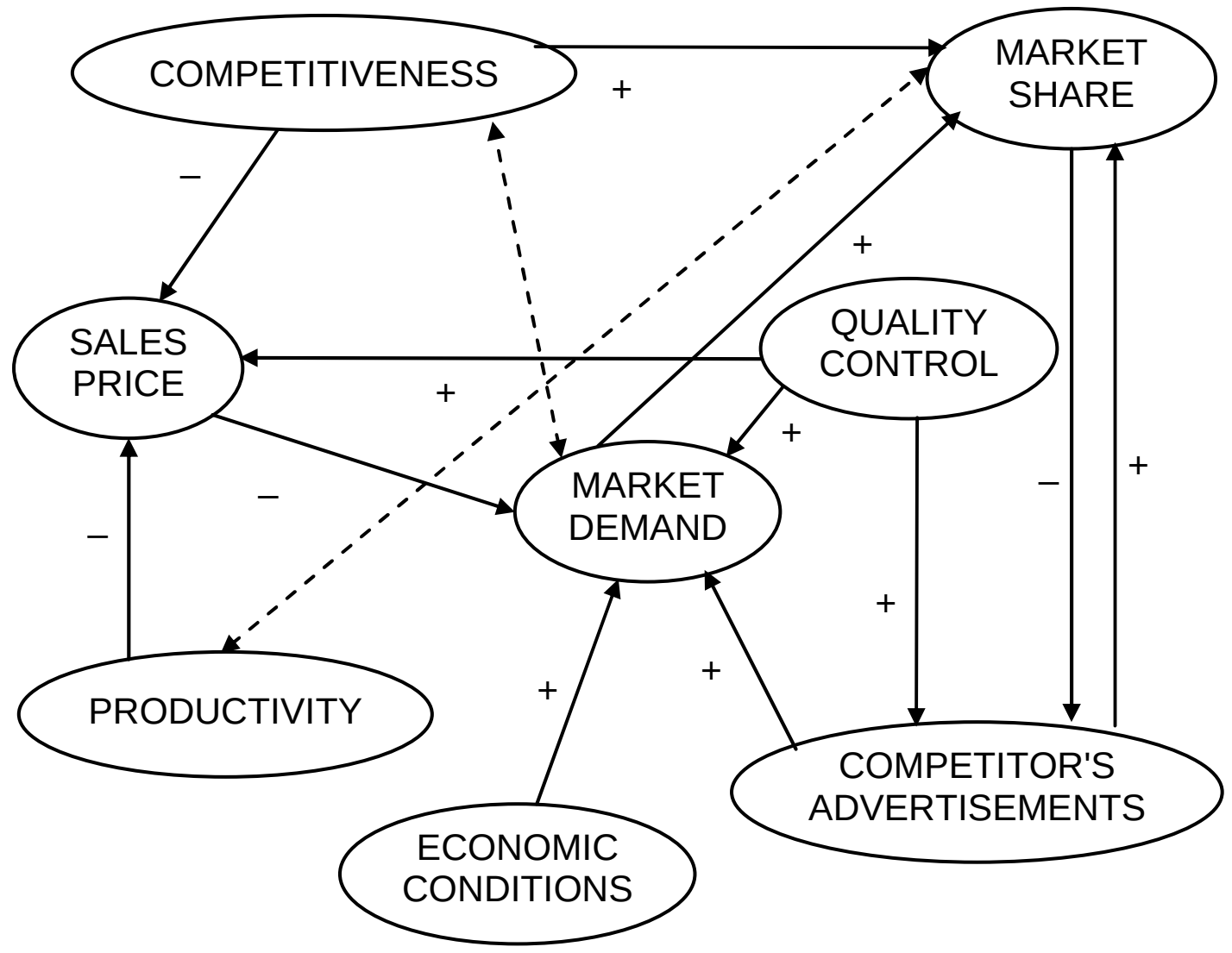

FIGURE: 3.7.9

It is pertinent to mention that K.C. Lee et al [38] have suggested that for better results more refined FCM with edge values can be used. In our opinion NCM may be a better solution to Lee's problems suggested in [38].

***Example 3.7.4:*** The application or use of adaptive fuzzy cognitive maps for hyper knowledge representation in strategy formation process was already has been carried out by [9]. For more about it please refer Carlsson and Fuller [9].

Adaptive Fuzzy Cognitive Maps can learn the weight from historical data. Once the FCM is trained, it lets us play what-if games (eg. what if demand goes up and prices remain stable? i.e. we improve our market position). Likewise adaptive neutrosophic cognitive maps are fuzzy cognitive maps with an addition concept between two nodes when the relation between them is an indeterminate. Once the NCM and the related FCM which is got when I’s are replaced by zeros is trained, let us play what-if games to predict the future in a realistic way.

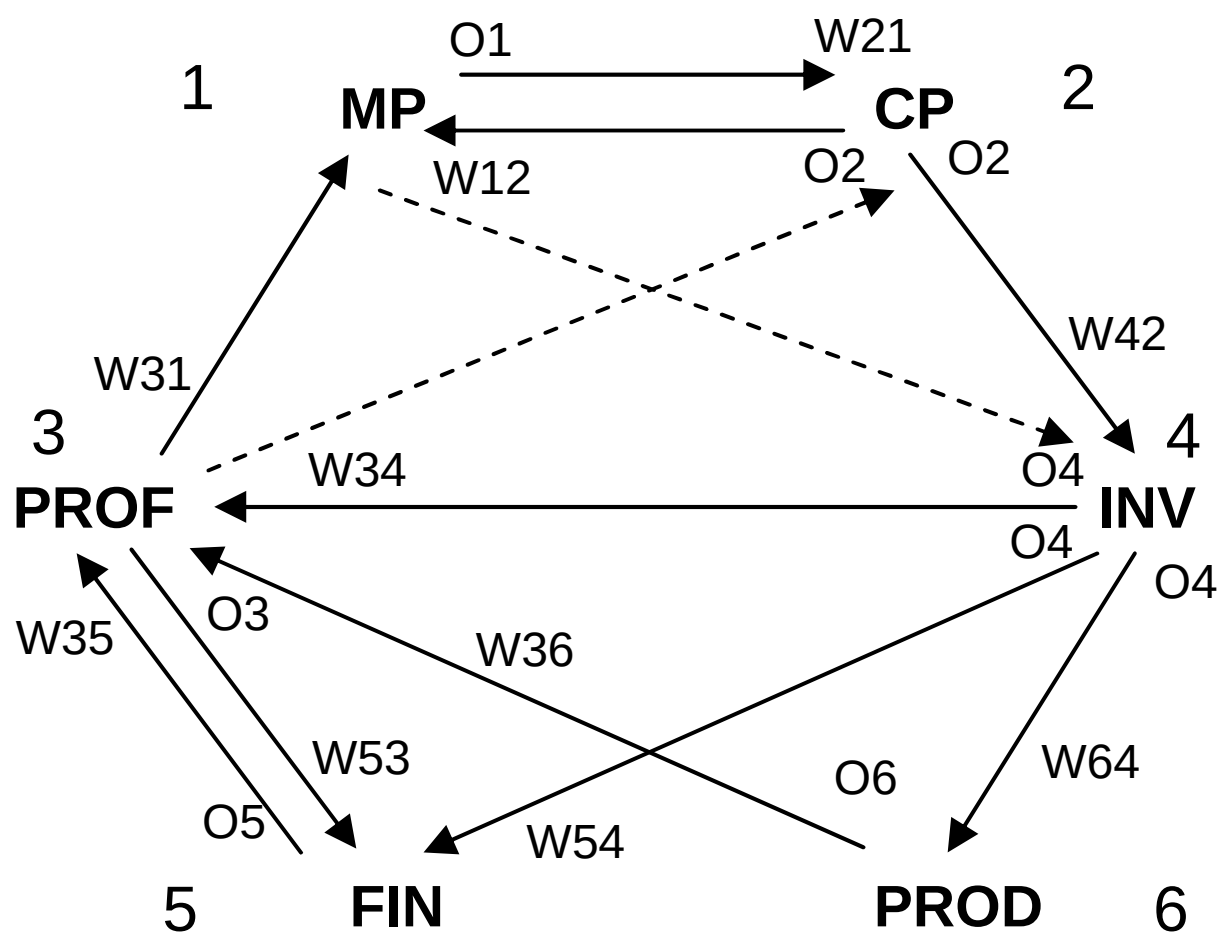

FIGURE: 3.7.10

Adaptive Neutrosophic Cognitive Map for the strategy formation process. The neutrosophic matrix N(w) is obtained using Figure 3.7.10. The reader is advised to study analyze and compare the results of NCM with FCM, for more about Adaptive FCM please refer [9]. From the Figure 3.7.10 we see according to the expert the relations between MP and INV is indeterminate also that between CP and PROF is an indeterminate.

Now we introduce the notion of a New Balance Degree for Neutrosophic Cognitive Maps (NCMs)

We first define when is an NCM imbalanced. In the opposite case, we say the NCM is balanced analogous to balance degree of FCM given by Tsadiras et al [61].

**DEFINITION 3.7.5:** *An NCM is imbalanced if we can find two paths between the same two nodes that create causal relations of different sign. In the opposite case the NCM is balanced. The term 'balanced' neutrosophic digraph is used in the following sense that is in a imbalanced NCM we cannot determine the sign or the presence of indeterminacy of the total effect of a concept to another.*

*Now on similar lines based on the idea that as the length of the path increases, the indirect causal relations become weakened the total effect should have the sign of the shortest path between two nodes.*

The degree to which a neutrosophic digraph of the NCM is balanced or imbalanced is given by the balance degree of the neutrosophic digraph. There are as in case of graph various types of balance degree in case of neutrosophic graphs also defined purely in an analogous way. An interested reader can obtain nice results on Balanced Degree of Neutrosophic Digraphs using results from [61].

***Example 3.7.5:*** Illustration of Neutrosophic Cognitive State maps of users web behavior is described. Searching for information in general is complex, with lot of indeterminacies and it is an uncertain process for it depends on the search engine; number of key words, sensitivity of search, seriousness of search etc. Hence we can see several of the factors will remain as indeterminates for the $C_1, C_2, \ldots, C_7$ as given by [42] we can remodel using NCM.

The NCM modeling of the users web behavior is given by the following neutrosophic digraph (Figure 3.7.11) and the corresponding N(E) built using an expert opinion is given by the following neutrosophic matrix:

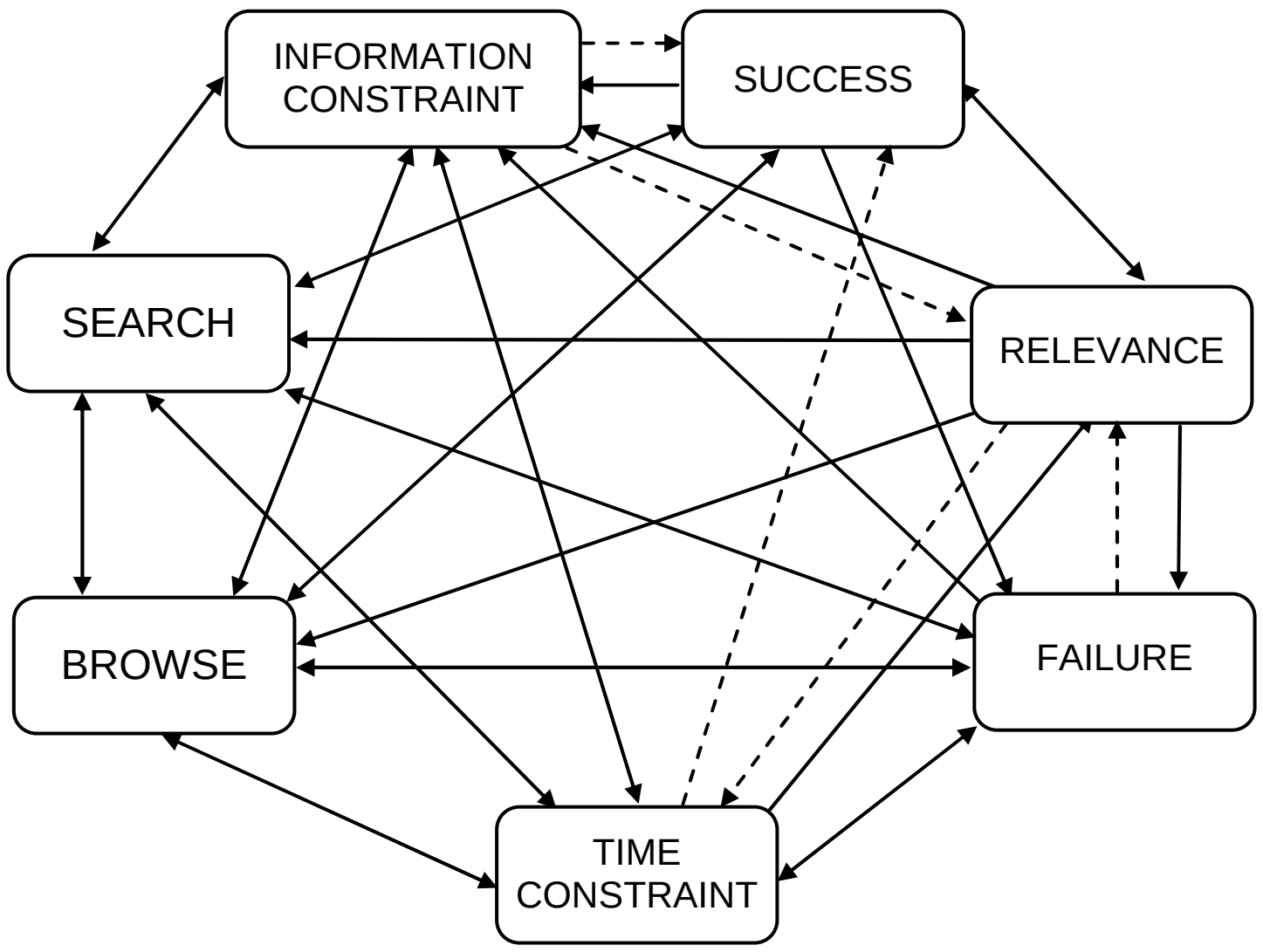

FIGURE: 3.7.11

$$N(E) = \begin{bmatrix} 0 & -1 & -1 & 1 & -1 & -1 & 1 \\ 1 & 0 & -1 & -1 & -1 & -1 & 1 \\ -1 & -1 & 0 & -1 & I & 1 & 1 \\ -1 & -1 & -1 & 0 & 0 & I & I \\ 1 & 1 & 1 & 1 & 0 & 1 & -1 \\ 1 & 1 & I & 1 & 1 & 0 & -1 \\ 1 & 1 & 1 & 1 & 0 & I & 0 \end{bmatrix}$$

Several results and conclusions can be derived for each of the state vectors.

The reader is given the work of comparing Fuzzy Cognitive State Map with the Neutrosophic Cognitive State Map given in Figure 3.7.11.

***Example 3.7.6:***

Now we describe the use of NCMs in Robotics. While it has been argued that FCMs are preferred for usage in robotics and applications of intimate technologies, owing to their ability to handle contradictory inputs, NCMs would be the more viable tool, for not only are they capable of handling contradictory inputs, but they can also handle indeterminacy.

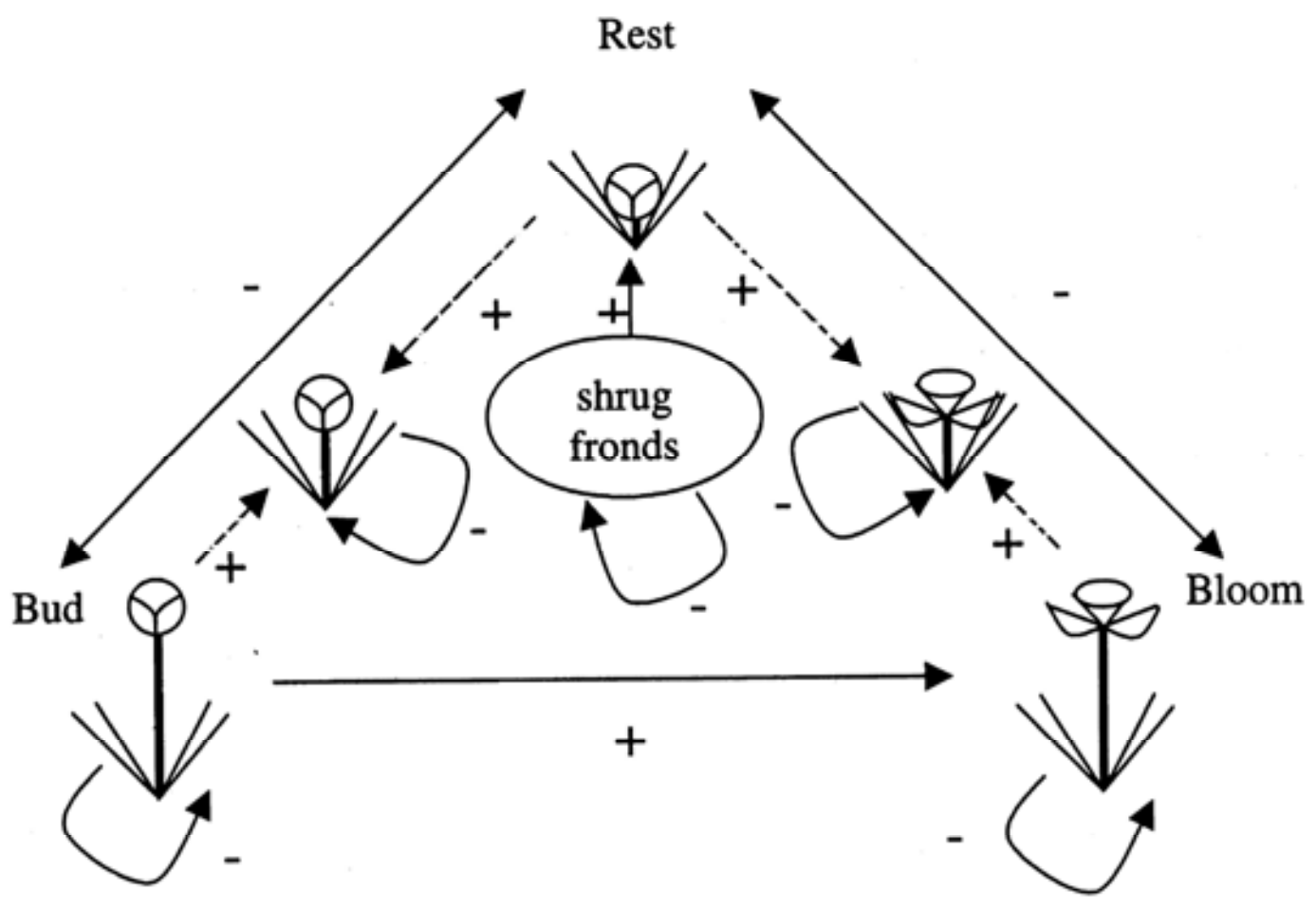

FIGURE: 3.7.12

Further FCMs have been used to model the Office Plant #1 to analyze the types of e-mails [9]. It is once again left as an exercise for the reader to use NCM in the place of FCM in this study. For categorically one cannot always divide the e-mails as official / non official, friendly / business like and so on for some can be termed as indeterminate, semi-friendly and semiofficial or so on and so forth.

So NCM can be adopted in mobile robots like Office Plant #1 and the study can be carried out as a maiden effort. A description of use of FCM can be derived even using this NCM.

For example we see at each stage the relation would be indeterminate if the email received has an over-lapping attributes in which case the section of the node may be indeterminate. Thus in the behaviour of the office plant, the dotted arrows ought to be adopted in situations where there is indeterminacy.

Thus in this case the number of indeterminate edges will be varying with time i.e., as in the case of the correspondence.

***Example 3.7.7:*** Here the Adaptation of NCMs to Model and Analyze Business Performance Assessment is given. The use of FCMs to model and analyze business performance assessment [30] was dealt.

For the study of NCM, introduce on the FCM the NCM structure so that one is in a position to analyze and construct a model with the FCM We give a FCM model given by [30].

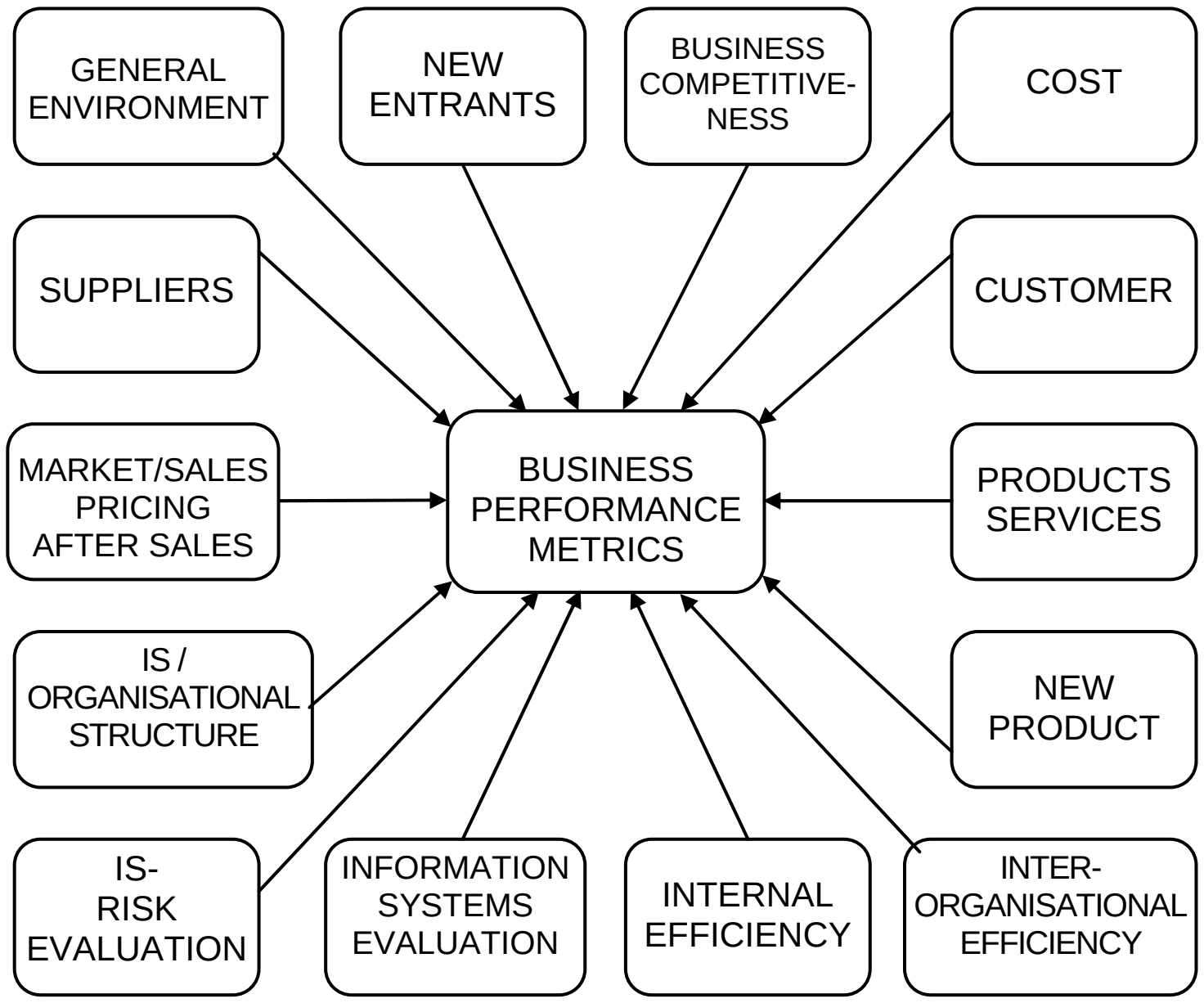

FIGURE: 3.7.13

The reader to requested to implement NCMs and draw conclusions based on the introduction of NCM to this model. Figure 3.7.13 describes the Business performance metrics given by [30].

It is suggested that in the NCM model in which some of the nodes can be considered as an indeterminate one leading to a strong neutrosophic graph be analyzed. However one modifies form is given in figure 3.7.14.

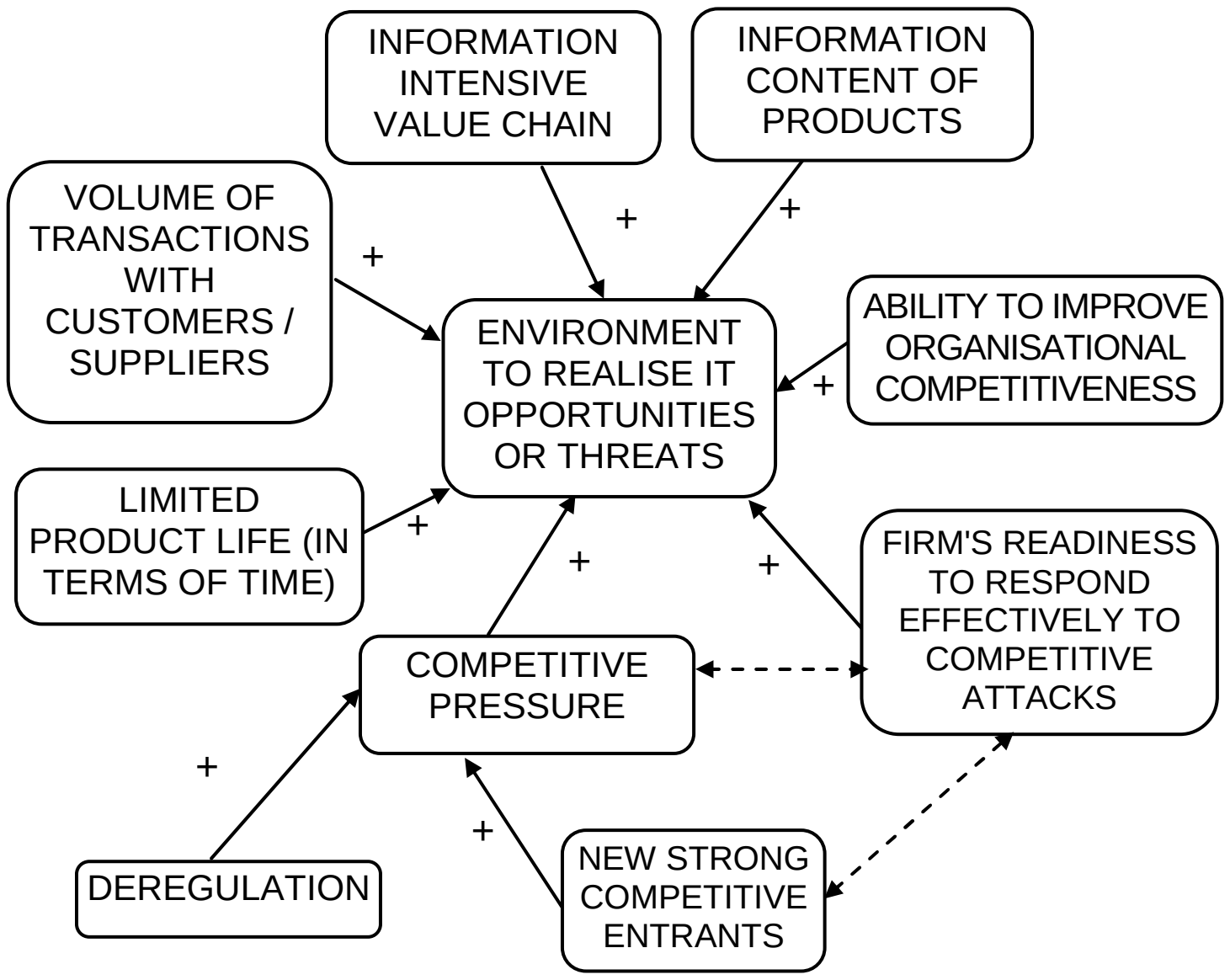

FIGURE: 3.7.14

***Example 3.7.8:*** Use of NCMs in legal rules

The implementation of legal rules using FCMs are very well studied by [1] and we have discussed [1].

Study NCM using Figure 3.7.15 and analyze and compare it with FCM given by [1]. Now we first show NCMs are better tools than the FCMs as in FCM we do not have the concept of indeterminacy. Only in case of NCMs we can say the relation between two concepts / attributes / nodes can be indeterminable also. For especially in criminal cases the concerned may not be able draw conclusions based on the data provided to him. Very many relations can be indeterminable so NCMs should be a better fit than the FCMs. Further instead of saying no relation between two nodes exist but still if the feelings exist with some doubt we cannot represent it in terms of FCMs but can easily implement NCMs so that while spelling out the judgment, the court can be very careful and give due weightages to the indeterminable relations.

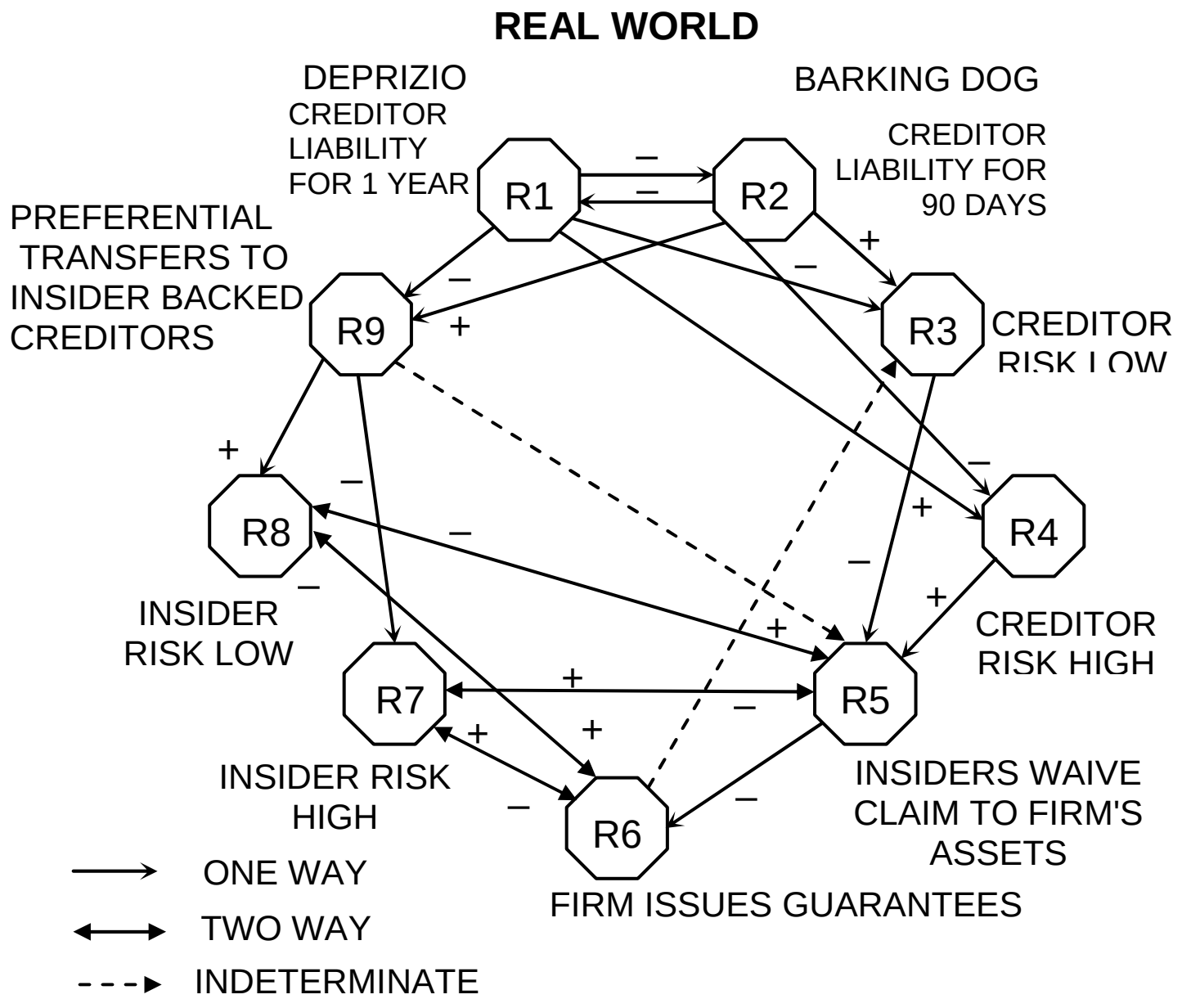

FIGURE: 3.7.15

The reader is given the task of implementing NCMs to the problem discussed in [1] using FCMs. They can also compare and contrast the FCM and NCM for this problem by studying where FCM has already been used and compare it with NCM and derive conclusions.

Now we illustrate how NCMs will play a major role in legal sides in several cases in India. We mention here a few types of them.

1. Encounter deaths with police.

2. Murder of political personalities which is very common in India

3. Custodial deaths in prison (Beaten to death, hanging, suicide etc.)

4. Undue delay in compensation cases in which government / private body is involved.

5. Means to punish intellectual harassment and torture and view it as more than the physical harassment and torture.

In all these five types of cases lot of indeterminacy is involved, so when NCM is applied certainly it will lead to better results.

***Example 3.7.9:*** Use of NCMs to find the driving speed in any one in freeway

Brubaker [8] used FCMs to create a model to find ones speed when driving in a California freeway. Thus FCM plays a major role in the study and analysis of transportation problem of all kind, for transportation problems are basically problems of decision-making. In our opinion, NCM can also be used to arrive at better results.

The concepts or nodes of the FCM are bad weather, freeway congestion, auto accidents, patrol frequency, own risk aversion, impatience and attitude.

Now if these are taken as nodes certainly we can have pairs of nodes for which the relation is indeterminate, for the concept of impatience and attitude with other concepts like bad weather or free way congestion and speed of others is an indeterminate in itself.

For fearing the bad weather one may be impatient and drive fast due to the fear that the weather may become worst or some other may fear bad weather (and consequent accidents) and be obsessed with fear and drive slow; the minute the nature of impatience or fear dominates a person certainly one cannot predict the speed, hence a lot of uncertainty and indeterminacy is involved.

So the adaptation of NCM may yield a better understanding and modeling of the problem than FCMs.

Thus we request the reader to model this problem using NCM and compare it when only FCM is applied; complete discussions using FCM is given in [65] as taken from [8]. For a slight suggestion of how to go about with this NCM, we have for the reader's sake provided a possible NCM graph in Figure 3.7.16.

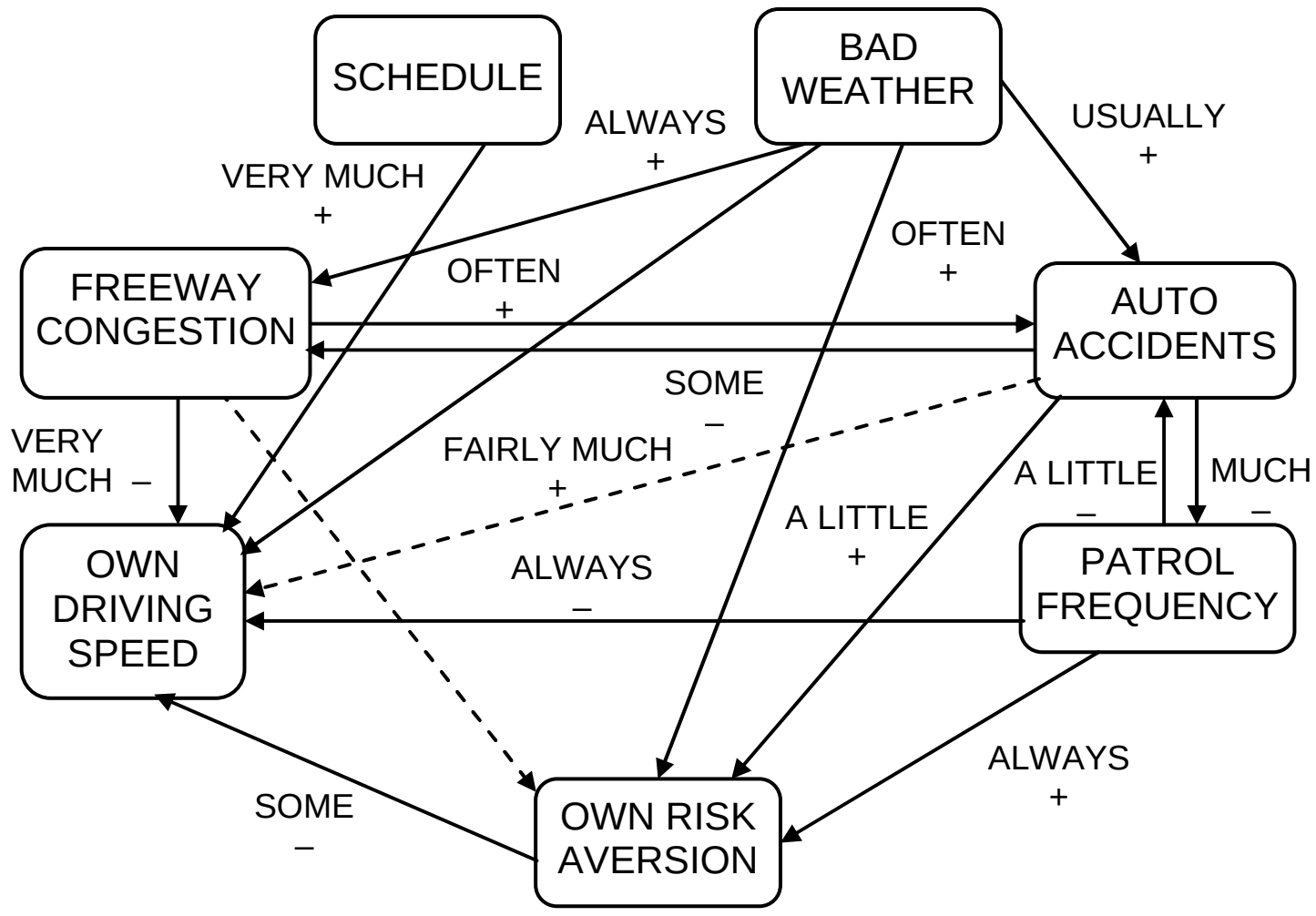

FIGURE: 3.7.16

As our main aim is to motivate researchers to use NCMs in place of FCMs whenever applicable and apply them to real world problems, we just give justification for the use of NCM and leave the work of constructing an NCM model to the reader. Even in this problem it is not only speed of one who drives but also several others factors like the speed of others, congestion etc. may or may not play an indeterminable role. Accidents are very common in the countries like India, where many other factors like bad roads; reckless driving by others; drunken driving etc. wreak havoc on the number of accident deaths. We can say in conclusion that problems related to traffic and transportation can be very efficiently handled with tools like the FCM and the NCM.

Next we shall see how best NCM can be applied in medical diagnosis. FCMs have found applications in many medical diagnostics including symptom disease model (in homeopathy), studying the depression of terminally ill patients and studies like death wish of terminally ill patients, etc. Since we can have these models where indeterminacy can exist between two nodes we can as well apply NCM in the place of FCMs.

**DEFINITION 3.7.6:** *Let D be the domain space and R be the range space with $D_1,\ldots, D_n$ the conceptual nodes of the domain space D and $R_1,\ldots, R_m$ be the conceptual nodes of the range space R such that they form a disjoint class i.e. $D \cap R = \phi$. Suppose there is a*

*FRM relating D and R and if at least a edge relating a $D_i$ $R_j$ is an indeterminate then we call the FRM as the Neutrosophic relational maps. i.e. NRMs.*

***Note:*** In everyday occurrences we see that if we are studying a model built using an unsupervised data we need not always have some edge relating the nodes of a domain space and a range space or there does not exist any relation between two nodes, it can very well happen that for any two nodes one may not be always in a position to say that the existence or nonexistence of a relation, but we may say that the relation between two nodes is an indeterminate or cannot be decided.

Thus to the best of our knowledge indeterminacy models can be built using neutrosophy. One model already discussed is the Neutrosophic Cognitive Model. The other being the Neutrosophic Relational Maps model, which are a further generalization of Fuzzy Relational Maps.

It is not essential when a study/ prediction/ investigation is made we are always in a position to find a complete answer. This is not always possible (sometimes or many a times) it is almost all models built using unsupervised data, we may have the factor of indeterminacy to play a role. Such study is possible only by using the Neutrosophic logic.

***Example 3.7.10:*** Female infanticide (the practice of killing female children at birth or shortly thereafter) is prevalent in India from the early vedic times, as women were (and still are) considered as a property. As long as a woman is treated as a property/ object the practice of female infanticide will continue in India.

In India, social factors play a major role in female infanticide. Even when the government recognized the girl child as a critical issue for the country's development, India continues to have an adverse ratio of women to men. Other reasons being torture of the in-laws may also result in cruel death of a girl child. This is mainly due to the fact that men are considered superior to women. Also they take into account that fact that men are breadwinners for the family. Even if women work like men, parents think that her efforts is going to end once she is married and enters a new family.

Studies have consistently shown that girl babies in India are denied the same and equal food and medical care that the boy babies receive. Girl babies die more often than boy babies even though medical research has long ago established that girls are generally biologically stronger as newborns than boys. The birth of a male child is a time for celebration, but the birth of female child is often viewed as a crisis. Thus the female infanticide cannot be attributed to single reason it is highly dependent on the feeling of individuals ranging from social stigma, monetary waste, social status etc.

Suppose we take the conceptual nodes for the unsupervised data relating to the study of female infanticide. We take the status of the people as the domain space D

$D_1$ – very rich
$D_2$ – rich
$D_3$ – upper middle class
$D_4$ – middle class
$D_5$ – lower middle class
$D_6$ – Poor
$D_7$ – Very poor.

The nodes of the range space R are taken as

$R_1$ – Number of female children - a problem
$R_2$ – Social stigma of having female children
$R_3$ – Torture by in-laws for having only female children
$R_4$ – Economic loss / burden due to female children
$R_5$ – Insecurity due to having only female children (They will marry and enter different homes thereby leaving their parents, so no one would be able to take care of them in later days.)

Keeping these as nodes of the range space and the domain space experts opinion were drawn which is given the following Figure 3.7.17:

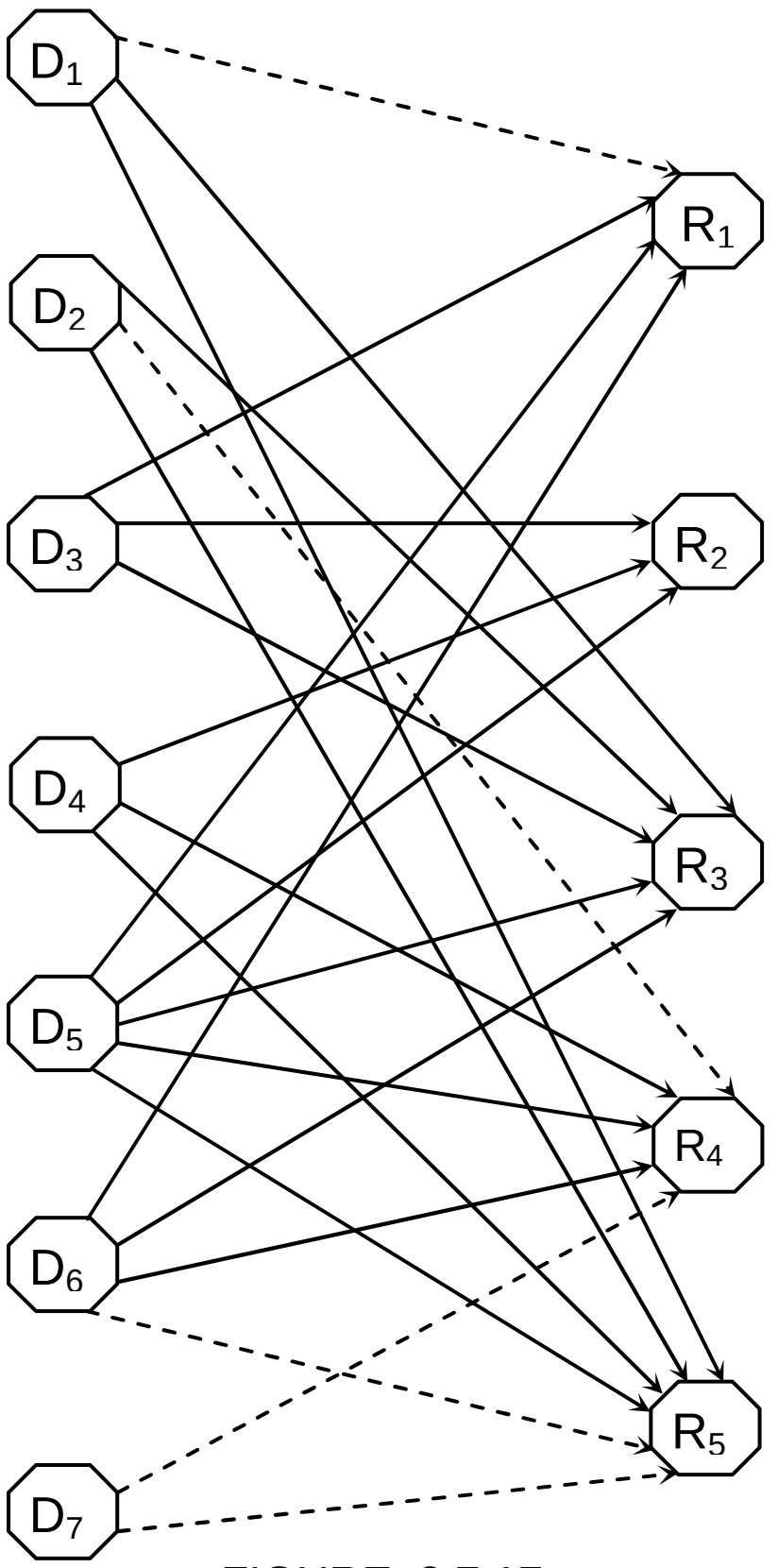

FIGURE: 3.7.17

Figure 3.7.17 is the neutrosophic directed graph of the NRM.

The corresponding neutrosophic relational matrix $N(R)^T$ is given below:

$$N(R)^T = \begin{bmatrix} I & 0 & 1 & 0 & 1 & 1 & 0 \\ 0 & 0 & 1 & 1 & 1 & 0 & 0 \\ 1 & 1 & 1 & 1 & 1 & 1 & 0 \\ 0 & I & 0 & 0 & 1 & 1 & I \\ 1 & 1 & 0 & 1 & 1 & I & I \end{bmatrix} \text{ and } N(R) = \begin{bmatrix} I & 0 & 1 & 0 & 1 \\ 0 & 0 & 1 & I & 1 \\ 1 & 1 & 1 & 0 & 0 \\ 0 & 1 & 1 & 0 & 1 \\ 1 & 1 & 1 & 1 & 1 \\ 1 & 0 & 1 & 1 & I \\ 0 & 0 & 0 & I & I \end{bmatrix}.$$

One can apply NRM using the neutrosophic directed graph given by the experts and analyze the problem

**DEFINITION 3.7.7:** *Let us assume that we are analyzing some nodes / concepts which are divided into three disjoint units. Suppose we have 3 spaces say P, Q and R. We say some m nodes in the space P, some n nodes in the space Q and some r nodes in the space R.*

*We can directly find NRMs or directed neutrosophic graphs relating P and Q and Q and R. But we are not in a position to link or get a relation between P and R directly but in fact there exists a hidden link between them which cannot be easily weighted, in such cases we use linked NRM.*

*Thus pairwise linked NRMs are those NRMs connecting three distinct spaces P, Q and R in such a way that using the pair of NRM, we obtain a NRM relating P and R. Thus if* $E_1$ *is the connection matrix relating P and Q then* $E_1$ *is a* $m \times n$ *matrix and* $E_2$ *is the connection matrix relating Q and R which is a* $n \times r$ *matrix.*

*Now consider P and R we are not in a position to link P and R directly by any directed graph but the product matrix* $E_1 E_2$ *gives a neutrosophic connection matrix between P and R. Also* $E^T_2 E^T_1$ *gives the neutrosophic connection matrix between R and P. When we have such a situation we call the NRMs as the pairwise linked NRMs.*

We will illustrate this definition explicitly by an example.

***Example 3.7.11:*** We just recall the example in [65] where the study of child labor is carried out using linked FRM. Now instead of FRM we instruct the experts that they need not always state the presence or absence of relation between any two nodes but they can also spell out the indeterminacy of any relation between two nodes, with these additional instruction to the experts, the opinions are taken.

The spaces under study are

G – the concepts / attributes associated with the government policies preventing / helping child labor.

C – attributes or concepts associated with children working as child laborers and P – attributes associated with public awareness and support of child labor.

G – Concepts associated with government policies:

- $G_1$ - Children do not form vote bank
- $G_2$ - Business men/industrialists who practice child labor are the main source of vote bank and finance
- $G_3$ - Free and compulsory education for children
- $G_4$ - No proper punishment given by Government to those who practice child labor.

Now we list out some of the attributes / concepts associated with the children working as laborers – C:

- $C_1$ - Abolition of child labor
- $C_2$ - Uneducated parents
- $C_3$ - School dropouts / children who never attended school
- $C_4$ - Social status of child laborers
- $C_5$ - Poverty / sources of living
- $C_6$ - Orphans runways parents are beggars or in prison.
- $C_7$ - Habits like smoking cinema, drinking etc.

Now we list out the attributes / concepts associated with public awareness or public supporting or exploiting the existence of child labor – P:

- $P_1$ - Cheap and long hours of labor with less pay
- $P_2$ - Children as domestic servants
- $P_3$ - Sympathetic public
- $P_4$ - Motivation by teachers to children to pursue education
- $P_5$ - Perpetuating slavery and caste bias.

Taking the experts opinion we first give the directed neutrosophic graph relating to child labor and the government policies in Figure 3.7.18.

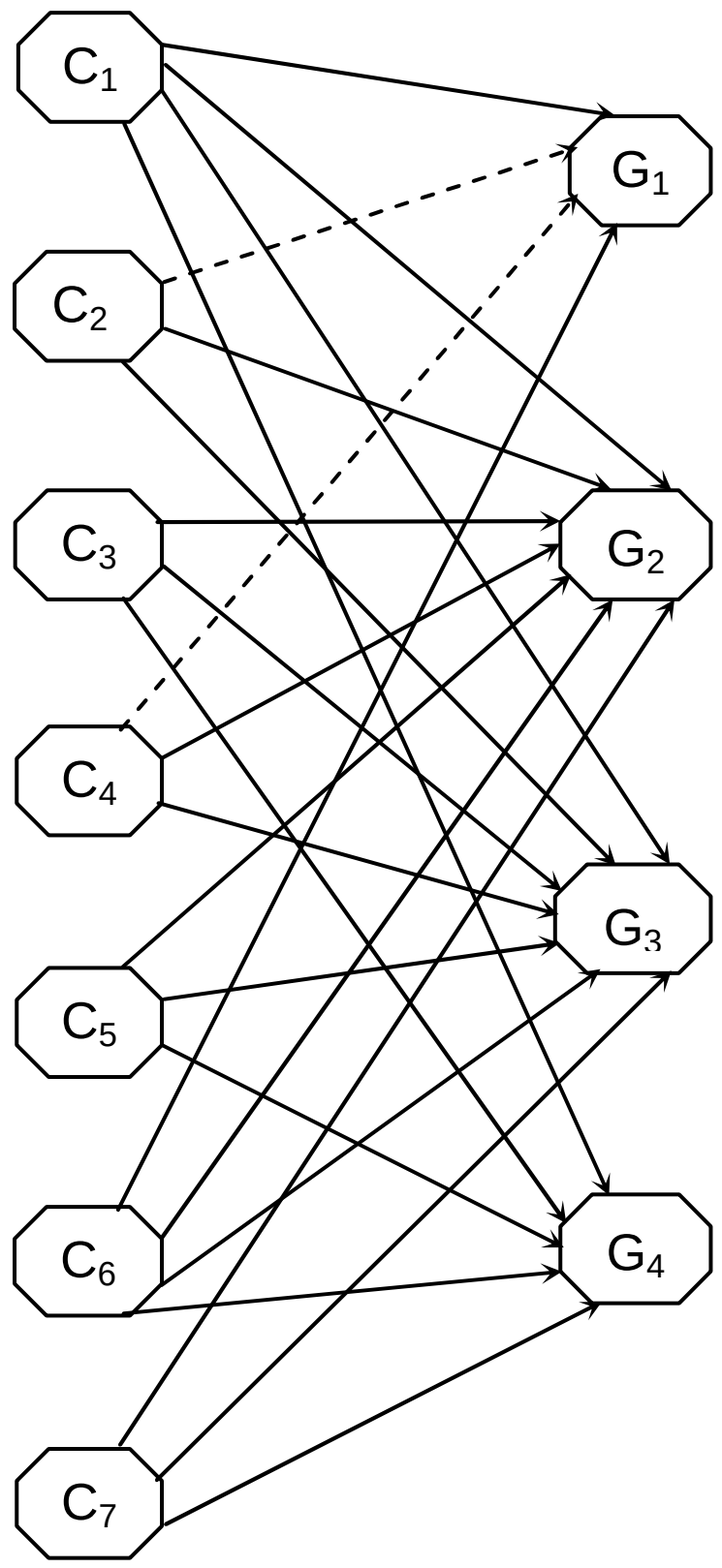

FIGURE: 3.7.18

The related neutrosophic connection matrix $N(E_1)$ is as follows:

$$N(E_1) = \begin{bmatrix} 1 & -1 & 1 & -1 \\ I & 1 & -1 & 0 \\ 0 & 1 & -1 & 1 \\ I & 1 & 1 & 0 \\ 0 & 1 & -1 & 1 \\ 1 & 1 & -1 & 1 \\ 0 & 1 & -1 & 1 \end{bmatrix}.$$

Now we are not interested in seeing the effect of instantaneous state vector on the neutrosophic dynamical system $N(E_1)$ but we are more interested in the illustration of how the model interconnects two spaces which have no direct relation, the same expert opinion is sought connecting neutrosophically the flourishing of child labor and the role played by the public. The neutrosophic directed graph relating the child labor and the public supporting child labor is given in Figure 3.7.19.

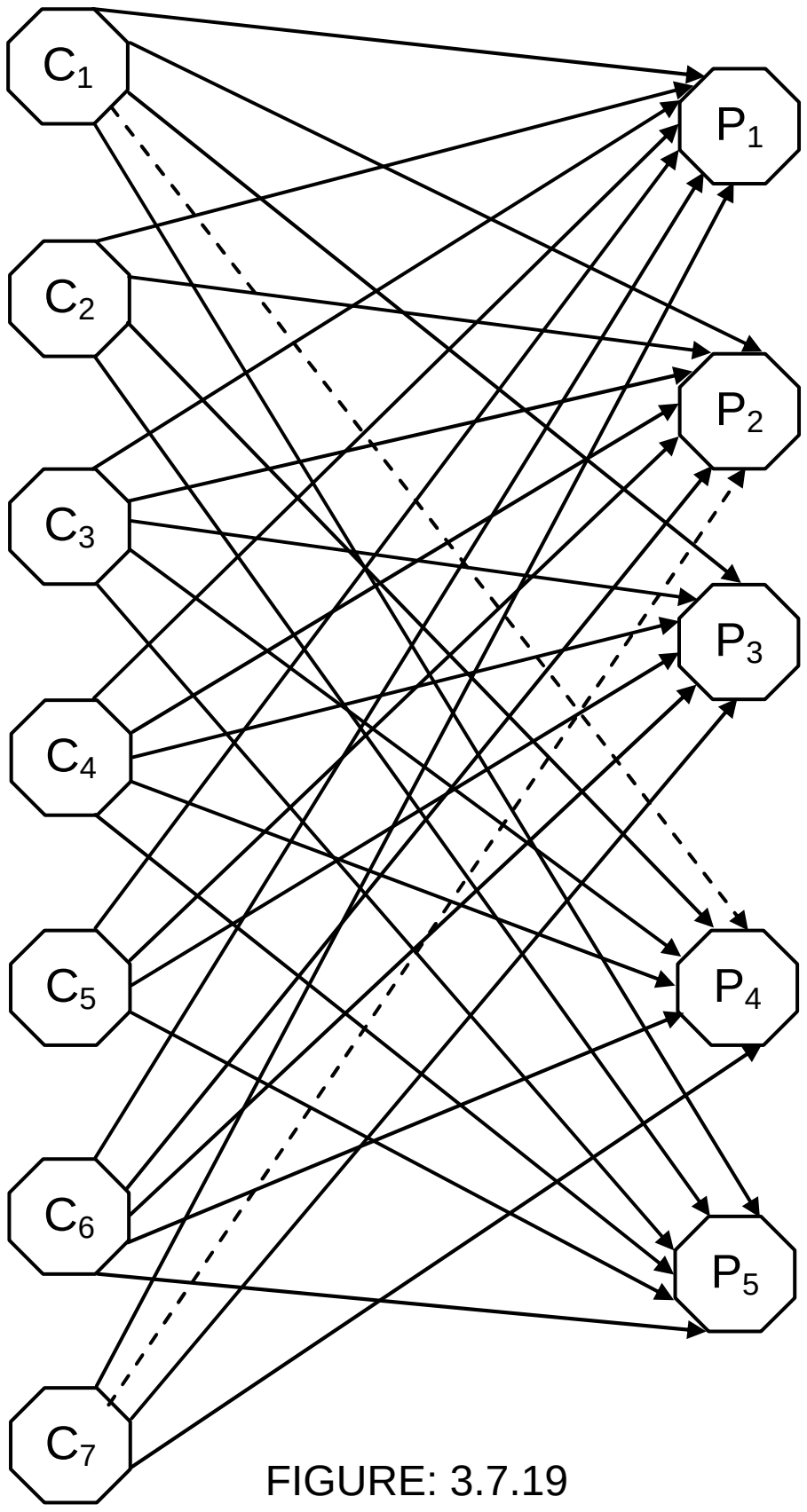

FIGURE: 3.7.19

The corresponding neutrosophic connection matrix related to the neutrosophic directed graph is given by $N(E_2)$.

$$N(E_2) = \begin{bmatrix} -1 & -1 & 1 & I & -1 \\ -1 & 1 & 0 & -1 & 1 \\ 1 & 1 & -1 & -1 & 1 \\ -1 & -1 & 1 & 1 & -1 \\ 1 & 1 & 1 & 0 & 1 \\ 1 & 1 & -1 & -1 & 1 \\ 1 & I & -1 & -1 & 0 \end{bmatrix}.$$

Thus $N(E_1)$ is a $7 \times 4$ matrix and $N(E_2)$ is a $7 \times 5$ matrix. Without the aid of the expert we can get the expert opinion relating the public and the government concerning child labor.

The neutrosophic relational matrix $[N(E_1)]^T N(E_2)$ is given as

$$N(E_2)]^T N(E_1) = \begin{bmatrix} 1 & 1 & 1 & -1 & 0 \\ 1 & 1 & -1 & -1 & 1 \\ -1 & -1 & 1 & 1 & -1 \\ 1 & 1 & -1 & -1 & 1 \end{bmatrix}.$$

Thus this maps variedly differs from the linked FRMs (discussed in [65]).

The reader is expected to model real world problems and apply linked NRMs.

**DEFINITION 3.7.8:** *Let P, Q and R be three spaces with nodes related to the same problem. Suppose P and Q are related by an NRM and Q and R are related by an NRM then the indirectly calculated NRM relating P and R got as the product of the neutrosophic connection matrices related with the NRMs. P and Q and Q and R is called the hidden neutrosophic connection matrix and the neutrosophic directed graph drawn using the Hidden neutrosophic matrix is called the Hidden neutrosophic directed graph of the pairwise linked NRMs.*

Now we proceed on to define on similar lines the three linked NRMs four linked NRMs, and in general n-linked NRMs.

**DEFINITION 3.7.9:** *Suppose we are analyzing a data for which the nodes / concepts are divided into four disjoint classes say A, B, C, D where A has n concepts / nodes B has m nodes, C has p concepts / nodes and D has q nodes / concepts.*

*Suppose the data under study is an unsupervised one and the expert is able to relate A and B, B and C and C and D through neutrosophic directed graphs that is the related NRMs and suppose the expert is not in a clear way to interrelate the existing relation between A and D then by using the neutrosophic connection matrixes we can get the neutrosophic connection matrix between A and D using the product of the three known matrix resulting in a* $n \times q$ *neutrosophic matrix which will be known as the hidden neutrosophic connection matrix and the related neutrosophic graph is called the hidden neutrosophic graph.*

*This structure is called as the 3-linked NRM.*

On similar lines we can by using 5 disjoint concepts or nodes of a problem define the 4-linked NRM. Thus if we have the nodes / concepts be divided into n + 1 disjoint classes of nodes then we can define the n-linked NRMs as in case of FRMs.

However we leave the task of constructing examples and applications of linked NRMs to the reader.

## 3.8 Neutrosophic Models versus Fuzzy models

We know in all mathematical analysis of the unsupervised data, not only uncertainty is predominant but also the concept of indeterminacy is in abundance. But in mathematical logic no one used the concept of indeterminacy till the year 1995.

Only in the year 1995 Florentine Smarandache introduced and studied the notion of indeterminacy, creating a further generalization of fuzzy logic, which was termed by him as Neutrosophic logic.

So whenever we introduce the notion of indeterminacy in our mathematical analysis we name the structure as Neutrosophic structure. Very recently (1999-2000) the notion of Fuzzy Relational Maps (FRM) was introduced. FRMs were a special particularization of FCMs when the data under study can be divided into disjoint sets. In this section we bring the comparison between FRMs and FCMs and the comparison of NRMs versus FRMs. We show that in the analysis of data, NRMs give a better and a realistic prediction than FRMs.

Here we give a comparison between FRMs and FCMs.

### 3.8.1 Fuzzy Relational Maps versus Fuzzy Cognitive Maps

FRMs are best suited when the data under consideration has its attributes or nodes to be divided into two (or more) disjoint classes.

FRMs cannot be applied when the nodes or attributes under consideration in a data cannot be divided into disjoint classes. Whenever the notion of FRMs are applied the main advantage is it minimizes labor and we need to work only with a smaller rectangular matrix which saves time and energy.

FCMs give a single hidden pattern arising from a fixed point or a limit cycle. On the other hand in the case of FRMs we have two hidden patterns for a given input vector one given by the range space and the other related to the domain space. Thus one is able to study the effect of the instantaneous vector considered in the domain space not only on the domain space but also on the range space.

The directed graph obtained from an FCM may or may not be bigraph. If the directed graph of the FCM is made into a bigraph then it implies and is implied that the FCM can be made into a FRM. All directed graphs of the FRMs are bigraphs.

If the nodes of the FCM is such that the related directed graph can never be made into a bigraph then it automatically implies that the FCM has nodes which cannot be made into two disjoint classes; so it cannot be made into FRMs.

Thus when FCMs can be converted into FRMs, they always enjoy better and a sensitive resultant apart from being economic and time saving.

### 3.8.2: Neutrosophic Relational Maps versus Neutrosophic Cognitive Maps

The NCMs are the Neutrosophic Cognitive Maps which are FCMs in which at least one of the directed edge which is an indeterminacy i.e. the directed graph of the NCM is a directed Neutrosophic graph. NRMs are the Neutrosophic Relational Maps i.e. the Neutrosophic directed graph is a Neutrosophic bigraph. The study of NRMs like FRMs is also economic and less time consuming.

The NRMs can be applied only in case the nodes / attributes of the unsupervised data can be divided into two disjoint sets; otherwise the data cannot have NRMs to be applied on it. The same notion in the language of Neutrosophic graphs is as follows:

If the associated Neutrosophic graph of an NCM can be made or is a Neutrosophic bigraph then certainly we can apply NRMs to the given data. If the data i.e. the Neutrosophic graph is never a Neutrosophic bigraph then certainly for that data NRMs can never be applied.

All NRMs can be easily changed into NCMs but NCMs in general cannot be always converted into NRMs. Just like FRMs, NRMs also gives two hidden patterns for a given instantaneous state vector a fixed point or a limit cycle one pertaining to the domain space and one fixed point or a limit cycle pertaining to the range space.

Also NRMs helps the influence of an instantaneous state vector on the space where it is taken and also in the other space.

Thus we can always say when the data is such that it is possible to apply NRMs it certainly yields a better conclusion than the NCMs.

Finally we give the comparison of NRMs and FRMs.

### 3.8.3: Neutrosophic Relation Maps versus Fuzzy Relational Maps

Now we recollect how NRMs are better than FRMs. In the first place in reality one cannot always say that relations between a node in the domain space is related with a node in the range space or not related with node in the range space. It may so happen that the existing relation between two nodes may not always be determinable by an expert. In FRMs there is no scope for such statement or such analysis, we can have a relation or no relation but this will not always be true in case of real world problems that too in case of unsupervised data, the relation can be an indeterminate, in such cases only NRMs are better disposed than FRMs. Thus NRMs play a better role and give a sensitive result than the FRMs.

Fuzzy world is about fuzzy data and fuzzy membership but it has no capacity to deal with indeterminate concepts, only Neutrosophy helps us to treat the notion of indeterminacy as a concept and work with it. Thus whenever in the resultant data we get the indeterminacy i.e. the symbol I the person who analyze the data can deal with more caution their by getting sensitive results than treating the nonexistence or associating 0 to that co-ordinate.

Thus from our study we have made it very clear that NRMs and NCMs are better tools yielding sensitive and truer results than FRMs and FCMs.

## Chapter Four

# Suggested Problems

In this chapter we have just suggested 32 problems most of them easy for any researcher only a few may be a little difficult. The main motivation for giving problems is to make the reader accustomed to the notion of neutrosophic graphs and its application to neutrosophic models

1. State and prove modified form / direct Euler theorem in case of

    a. neutrosophic vertex graphs
    b. neutrosophic edge graphs.
    c. neutrosophic graphs

2. Does every connected neutrosophic graph contain a neutrosophic spanning tree?

3. Give an example of a neutrosophic tree with 6 vertex and 5 indeterminate vertex i.e., $G_5$.

4. Give an example of a neutrosophic graph which is not a neutrosophic tree.

5. Can we say the number of neutrosophic edges of a tree with n indeterminate vertex is n – 1?

6. Will a connected neutrosophic graph with n indeterminate vertex and n – 1 neutrosophic edges be a neutrosophic tree?

7. Does every neutrosophic tree have a neutrosophic center consisting of either a single indeterminate vertex or two adjacent vertices which are indeterminates?

8. Obtain some applications of fuzzy model as graph models.

9. Can study of two model relation lead to their graph relations and vice versa?

10. For a neutrosophic connected graph, are the following statement equivalent.
    a. E is neutrosophic Eulerian
    b. The degree of each vertex of G is an even positive integer (Here the neutrosophic edge or indeterminate edge is also counted as an edge).
    c. G is the edge disjoint union of cycles. (Edge includes indeterminate edges also).

11. Is the following statements true?
    "A connected neutrosophic graph is neutrosophic Eulerian if and only if it admits a cycle decomposition?"
    or
    "A connected neutrosophic graph is neutrosophic Eulerian if and only if it admits a neutrosophic cycle?" Justify your answer!

12. A neutrosophic graph is strong neutrosophic Eulerian if and only if each edge e (edge does not include indeterminate edge) of G belongs to an odd number of cycles of G". Justify or refute the hypothesis!
    **Remark:** If in the above problem we replace "cycles" by "neutrosophic cycles" and the edge also includes indeterminate edge or neutrosophic edge, does the hypothesis given in the above problem true? Justify your claim.

13. Does the following statement hold true in case of strong neutrosophic Eulerian graph.
    "A neutrosophic graph is strong neutrosophic Eulerian if and only it has an odd number of neutrosophic cycle decomposition – Justify your claim.

14. A weak neutrosophic digraph D is neutrosophic Eulerian if and only if every point of D has equal indegree and outdegree. (Here we assume all vertices are real and none of them are in determinates)

15. Obtain an analogous result in case of neutrosophic Eulerian digraph.
Is it true in case of neutrosophic Eulerian digraph the number of Eulerian neutrosophic trials is

$$C.\prod_{i=1}^{p} (di - 1)\,!$$

where di = id ($v_i$ and C is the common value of all the co factors).

16. Find $\chi_N$ (G) for the graph G in the following figure.

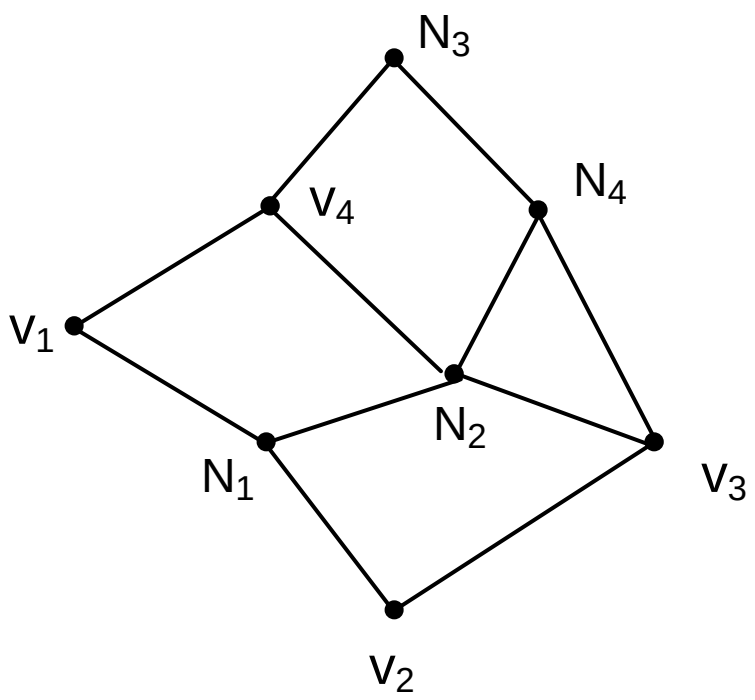

FIGURE: 4.1

Is $\chi$(G) = $\chi_N$ (G) ?

17.

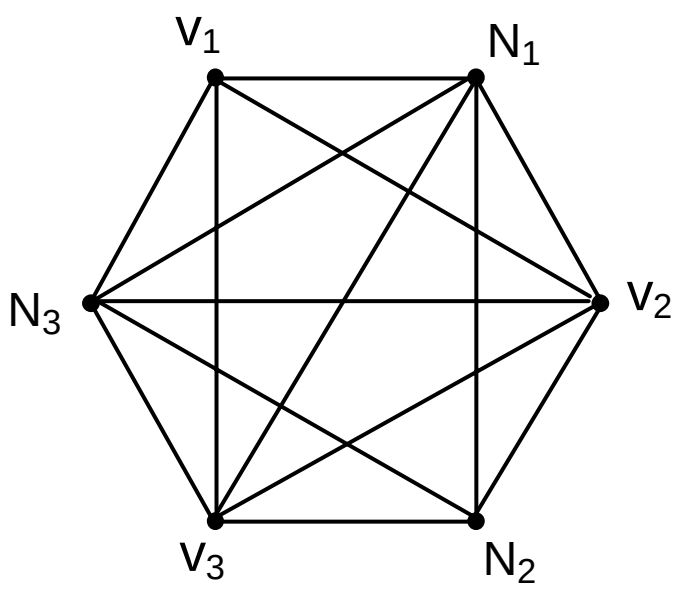

FIGURE: 4.2

For this neutrosophic graph find $\chi_N$ (G).

18. Compare the value of $\chi_N$ (G) and $\chi$ (H) where H is the graph

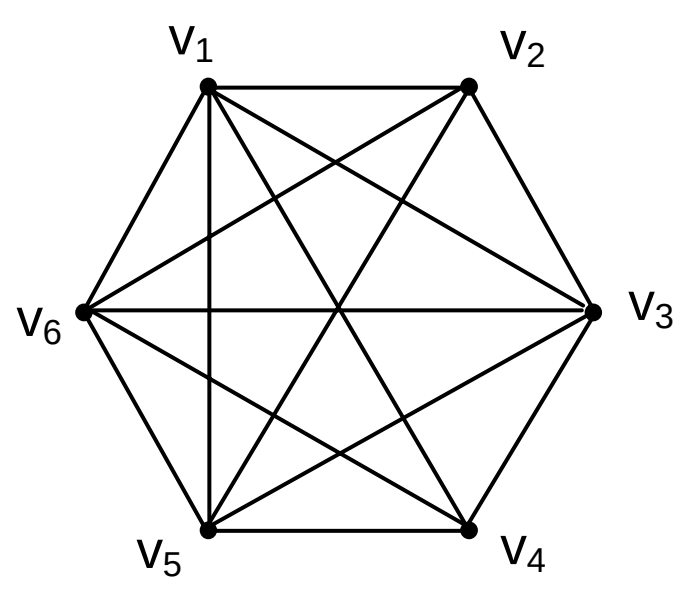

FIGURE: 4.3

and for the neutrosophic graph G given in the problem 17. for which $\{v_1\ v_2\ v_3\}$ are real vertices and $N_1$, $N_2$, $N_3$ are indeterminate vertices.

19. What is $\chi_N$ (G' ) for the neutrosophic graph G'.

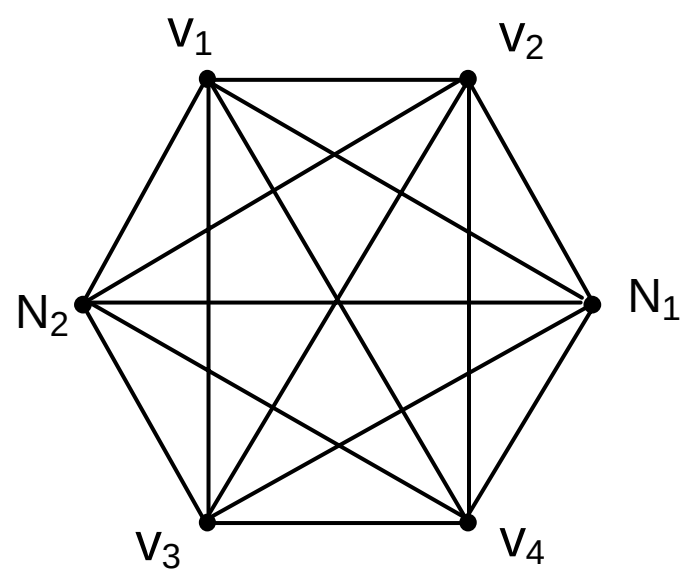

FIGURE: 4.4

a. Compare $\chi_N$(G) with $\chi_N$(G') G given in (17).
b. Compare $\chi$ (H) with $\chi_N$ (G'), H given in (18).

20. Find $\chi_N$ (G) when G is a neutrosophic graph, which is bipartite.

21. Can we find all neutrosophic graphs in which $\chi_N$(G) = 2?

22. Obtain interesting results relating $\chi_N$ (G) and special types of neutrosophic graphs.

23. Is in a loopless bipartite neutrosophic graph G $\chi'_N$ (G) = $\Delta$ (G)?

24. Find $\chi'_N$ (G) for the following graphs.

a.

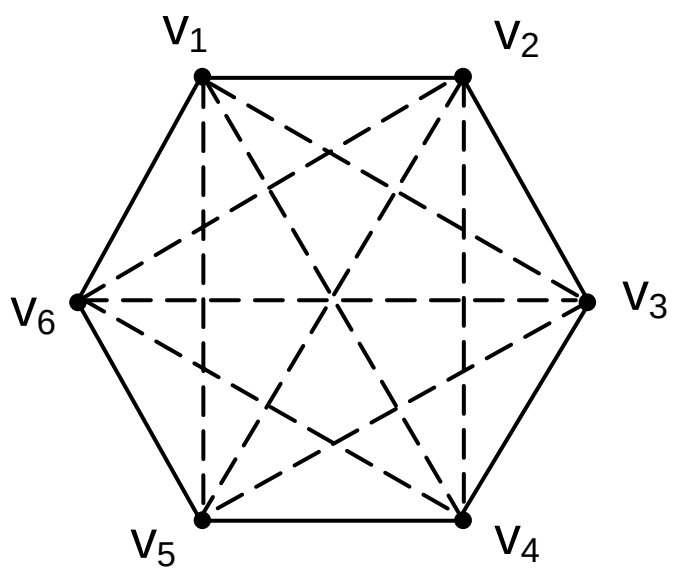

FIGURE: 4.5

b.

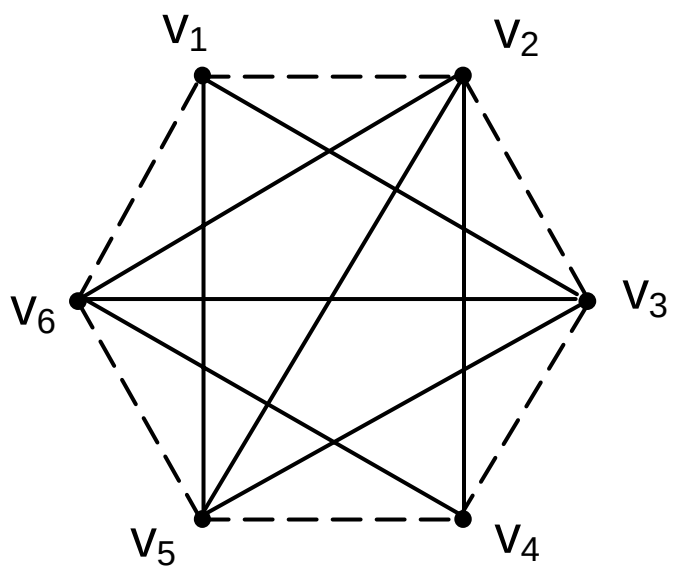

FIGURE: 4.6

c. Compare $\chi'_N$ (G) and $\chi'_N$ (G').

25. Find $\chi'$ (G) for the following two graphs B and B' given in figures 4.7 and 4.8 respectively and compare them.

a.

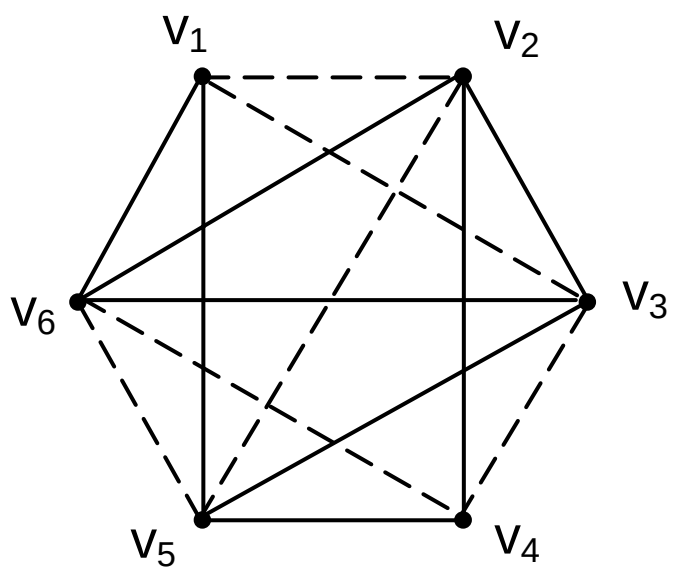

FIGURE: 4.7

Find $\chi'_N$ (B).

b.

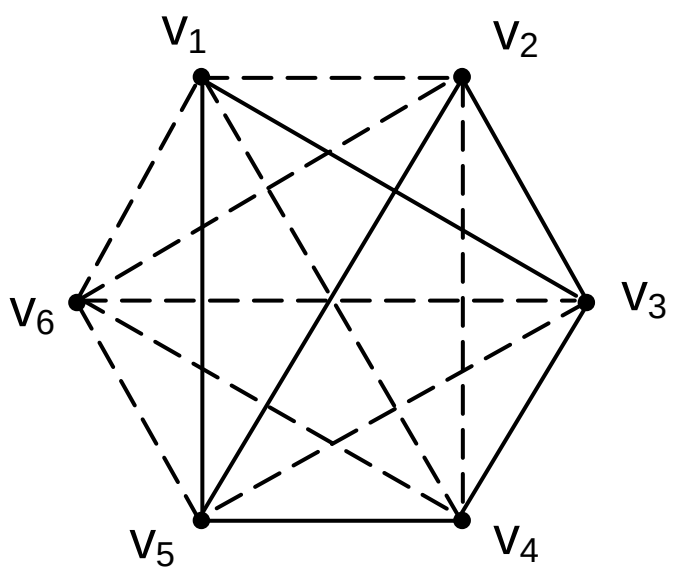

FIGURE: 4.8

Find $\chi'_N$ (B' ).

c. Compare $\chi'_N$ (B) with $\chi'_N$ (B' ).

26. Define neutrosophic simple graph and illustrate it with examples.

27. Find the total number of neutrosophic Petersen graphs.

28. Find the number of Petersen graphs which has only indeterminate vertices i.e., find all neutrosophic point graphs which are neutrosophic Petersen graphs.

29. Find the number of neutrosophic Petersen graphs in which all the 10 vertices are real.

30. Define and obtain some interesting properties about neutrosophic Eulerian graphs.

31. Obtain a analog of Harary – Nash Williams Theorem for a neutrosophic graph.

32. Define neutrosophic Pancyclic graph and obtain some interesting results about them.

# Bibliography

Here we give the basic materials, which we have used and also we do not promise to provide all existing references. The references, which we could get, are given.

1. **Adams, E.S., and D.A. Farber.** Beyond the Formalism Debate: Expert Reasoning, Fuzzy Logic and Complex Statutes, *Vanderbilt Law Review*, **52** (1999), 1243-1340. http://law.vanderbilt.edu/lawreview/vol525/adams.pdf

2. **Appel, K., and Haken, W.,** *The solution of the four color map problem*, Scientific American, 237 (4), (1977) 108-121.

3. **Balakrishnan, R., and K. Ranganathan**, *A textbook of Graph Theory*, Springer, 1999.

4. **Balakrishnan, R., and Paulraja, P.,** Line graphs of subdivision graphs, *J. Combin. Info. and Sys. Sci.*, 10, (1985) 33-35.

5. **Bondy, J.A., and Murthy, U.S.R.,** *Graph Theory and applications*, The Mac Milan Press 1976.

6. **Bondy, J.A.,** Pancyclic graphs, *J. Combin Theory Ser*. 11, (1971), 80-84.

7. **Brannback, M., L. Alback, T. Finne and R. Rantanen.** Cognitive Maps: An Attempt to Trace Mind and Attention in Decision Making, *in* C. Carlsson ed. *Cognitive Maps and Strategic Thinking*, Meddelanden Fran Ekonomisk Statsvetenskapliga Fakulteten vid Abo Akademi Ser. **A 442** (1995) 5-25.

8. **Brubaker, D.** Fuzzy Cognitive Maps, *EDN ACCESS*, 11 April,1996. http://www.e-insite.net/ednmag/archives/1996/041196/08column.htm

9. **Carlsson, C., and R. Fuller.** Adaptive Fuzzy Cognitive Maps for Hyper-knowledge Representation in Strategy Formation Process In *Proceedings of the International Panel Conference on Soft and Intelligent Computing*, Technical Univ. of Budapest, (1996) 43-50. http://www.abo.fi/~rfuller/asic96.pdf

10. **Carvalho, J.P., and Jose A.B. Tomè.** Rule-based Fuzzy Cognitive Maps and Fuzzy Cognitive Maps – a Comparative Study. In *Proc. of the 18$^{th}$ International Conference of the North American Fuzzy Information Processing Society*, by NAFIPS, New York, (1999) 115-119. http://digitais.ist.utl.pt/uke/papers/NAFIPS99rbfcm-fcm.pdf

11. **Chartrand, G. and Wall, C.E.,** On the Hamiltonian Index of a Graph, *Studia Sci. Math. Hungar*, 8 1973, 43-48.

12. **Chvatal, V., and Erdos, P.,** A note on Hamiltonian Circuits, Discrete Maths (2) 1972 111-113.

13. **Craiger, J.P., R.J. Weiss, D.F. Goodman, and A.A. Butler.** Simulating Organizational Behaviour with Fuzzy Cognitive Maps, *Int. J. of Computational Intelligence and Organization*, **1** (1996) 120-123.

14. **De Baets, B., and Kerre, E.E.**, Fuzzy Relational Compositions, *Fuzzy sets and Systems*, 60, (1993), 109-120.

15. **Dirac, G.A.,** Some Theorems on Abstract Graphs, *Proc. London Math Soc*., 2 (1952) 69-81.

16. **Fieischner, H.,** Eulerian Graphs and related topics, *Annals of Disc. Math*, 45, North Holland, 1990.

17. **Fiorini, S., and Wilson, R.J.**, Edge colorings of graphs, *In Research Notes in Mathematics*, 16, Pittman, London, 1971.

18. **Gibbons, A.,** *Algorithmic Graph Theory*, Cambridge. Univ. Press, 1985.

19. **Golambic, M.C.,** *Algorithmic Graph Theory and Perfect graphs*, Academic Press 1980.

20. **Goto, K. and T. Yamaguchi.** Fuzzy Associative Memory Application to a Plant Modeling, *in Proc. of the International Conference on Artificial Neural Networks*, Espoo, Finland, (1991) 1245-1248.

21. **Hall Jr, M.,** Combinatorial Theory, Blaisdell, Waltham, 1967.

22. **Harary, F., and C.St J.A. Nash-Williams**, On Eulerian and Hamiltonian graphs and line graphs, *Canada Math Bull*, 8 (1965) 701-710.

23. **Harary, F.**, *Graph Theory*, Addison-Wesley, 1969.

24. **Higashi, M A Di Nola S. Sessa and Pedrycz, W.**, Ordering fuzzy sets by consensus concept and fuzzy relational equations, *Intern. J of General Systems,* 10, (1984), 47-56.

25. **Holton, D.A., and Sheehan, J.,** The Petersen Graph, *Australian Math Soc. Lecture Series 7,* Cambridge Univ Press 1993.

26. **Jaeger, F.,** A note on Sub Eulerian Graphs, *J. Graph Theory*, 3 (1979), 91-93.

27. **Jensen T.R., and Toft, B.,** Graph Colouring Problems, *John Wiley & sons*, 1995.

28. **Junger, M., Pulley Blank, W.R., and Reinelt, G.,** On Partitioning the edges of graphs into connected subgraphs, *J. Graph Theory*, 9 (1985) 539-549.

29. **Kainen, P.C. and Saaty, T.L.,** The Four Color Problems, *Dover Pub. Inc.*, 1977.

30. **Kardaras, D., and G. Mentzas.** Using fuzzy cognitive maps to model and analyze business performance assessment, In *Prof. of Int. Conf. on Advances in Industrial Engineering – Applications and Practice II*, Jacob Chen and Anil Milal (eds.), (1997) 63-68.

31. **Kim, H.S., and K. C. Lee.** Fuzzy Implications of Fuzzy Cognitive Maps with Emphasis on Fuzzy Causal Relations and Fuzzy Partially Causal Relationship, *Fuzzy Sets and Systems*, **97** (1998) 303-313.

32. **Klir, G.J., and Yuan, B**., *Fuzzy Sets and Fuzzy Logic: Theory and Applications*, Prentice-Hall, Englewood Cliffs NJ, 1995.

33. **Kosko, B.** Fuzzy Cognitive Maps, *Int. J. of Man-Machine Studies*, **24** (1986) 65-75.

34. **Kosko, B.**, *Neural Networks and Fuzzy Systems: Dynamical Approach to Machine Intelligence*, Prentice-Hall, Englewood Cliffs NJ, 1992.

35. **Kundu, S.,** Bounds on the number of disjoint spanning trees, *J. Combin. Theory* 17, (1974), 199-203.

36. **Lee, K.C., Kim, J.S., Chang, N.H., and Kwon, S.J.,** Fuzzy Cognitive Map Approach to Web-mining Inference Amplification, *Expert Systems with Applications*, **22** (2002) 197-211.

37. **Lee, K.C., S. C. Chu and S.H. Kim.** Fuzzy Cognitive Map-based Knowledge Acquisition Algorithm: Applications to Stock Investment Analysis, in W.Cheng, Ed., *Selected Essays on Decision Science* (Dept. of Decision Science and Managerial Economics), The Chinese University of Hong Kong, (1993), 129-142.

38. **Lee, K.C., W.J. Lee, O.B. Kwon, J.H. Han, P.I. Yu.** A Strategic Planning Simulation Based on Fuzzy Cognitive Map Knowledge and Differential Game, *Simulation*, **71** (1998) 316-327.

39. **Liu, F., and F. Smarandache.** Logic: A Misleading Concept. A Contradiction Study toward Agent's Logic, in *Proceedings of the First International Conference on Neutrosophy, Neutrosophic Logic, Neutrosophic Set, Neutrosophic Probability and Statistics,* Florentin Smarandache editor, Xiquan, Phoenix, ISBN: 1-931233-55-1, 147 p., 2002, *also published in* "Libertas Mathematica", University of Texas at Arlington, **22**

(2002) 175-187. http://lanl.arxiv.org/ftp/math/papers/0211/0211465.pdf

40. **Liu, F., and Smarandache, F.**, *Intentionally and Unintentionally. On Both, A and Non-A, in Neutrosophy*. http://lanl.arxiv.org/ftp/math/papers/0201/0201009.pdf

41. **Lovasz, L. and Plummer, M.D.,** *Matching Theory*, Annals of Discrete Mathematics, 27, (1986), 121.

42. **Meghabghab, G.** Fuzzy Cognitive State Map vs. Markovian Modeling of User's Web Behaviour, Invited Paper, *International Journal of Computation Cognition*, **1** (Sept. 2003), 51-92. Article published electronically on December5,2002). (http://www.YangSky.com/yangijcc.htm)

43. **Moon, I.W.**, *Topics in Tournaments*, Holt Rinehart and Winston Inc., New York, 1968.

44. **Mycielski, J.,** Surle Coloriage des graph, *Colloq math*, 3 (1955) 161-162.

45. **Nebesky, L.,** On the line graph of square and the square of the line graph of a connected graph, *Casopis Pset Mat*, 98 (1973) 285-287.

46. **Oberly, D.J., and Summer, P.P.,** Every connected locally connected non trivial graph with no induced claw is Hamiltonian, *J. Graph Theory*, 3 (1979) 351-356.

47. **Ore, O.,** *The Four Colour Problem*, Academic Press New York 1967.

48. **Ozesmi, U.** Ecosystems in the Mind: Fuzzy Cognitive Maps of the Kizilirmak Delta Wetlands in Turkey, Ph.D. Dissertation titled *Conservation Strategies for Sustainable Resource use in the Kizilirmak Delta-Turkey*, University of Minnesota, (1999) 144-185. http://env.erciyes.edu.tr/Kizilirmak/UODissertation/uozesmi5.pdf

49. **Pelaez, C.E., and J.B. Bowles.** Applying Fuzzy Cognitive Maps Knowledge Representation to Failure Modes Effects Analysis, In *Proc. of the IEEE Annual*

*Symposium on Reliability and Maintainability*, (1995) 450-456.

50. **Pelaez, C.E., and J.B. Bowles.** Using Fuzzy Cognitive Maps as a System Model for Failure Modes and Effects Analysis, *Information Sciences*, **88** (1996) 177-199.

51. **Serre, J.P.**, *Trees*, Springer Vertage 1980.

52. **Siraj, A., S.M. Bridges, and R.B. Vaughn.** *Fuzzy cognitive maps for decision support in an intelligent intrusion detection systems*, www.cs.msstate.edu/~bridges/papers/nafips2001.pdf

53. **Smarandache, F.**, *A Unifying Field in Logics: Neutrosophic Logic. Neutrosophy, Neutrosophic Set, Neutrosophic Probability and Statistics*, third edition, Xiquan, Phoenix, 2003.

54. **Smarandache, F.**, Definitions Derived from Neutrosophics, In Proceedings of the *First International Conference on Neutrosophy, Neutrosophic Logic, Neutrosophic Set, Neutrosophic Probability and Statistics*, University of New Mexico, Gallup, 1-3 December 2001

55. **Styblinski, M.A., and B.D. Meyer.** Fuzzy Cognitive Maps, Signal Flow Graphs, and Qualitative Circuit Analysis, in *Proc. of the 2nd IEEE International Conference on Neural Networks (ICNN – 87)*, San Diego, California (1988) 549-556.

56. **Styblinski, M.A., and B.D. Meyer.** Signal Flow Graphs versus Fuzzy Cognitive Maps in Applications to Qualitative Circuit Analysis, *Int. J. of Man-machine Studies*, **18** (1991) 175-186.

57. **Stylios, C.D., and P.P. Groumpos.** The Challenge of Modeling Supervisory Systems using Fuzzy Cognitive Maps, *J. of Intelligent Manufacturing*, **9** (1998) 339-345.

58. **Summer, D.P.,** Graphs with 1-factors, *Proc. Amer Math Soc*, 42 (1974) 8-12.

59. **Taber W. R.** Fuzzy Cognitive Maps Model Social Systems, *Artificial Intelligence Expert*, **9** (1994) 18-23.

60. **Toida, S.,** Properties of an Euler graph, *J. Franklin Inst.*, **295** (1973) 346-346.

61. **Tsadiras, A.K., and K.G. Margaritis.** *A New Balance Degree for Fuzzy Cognitive Maps*, http://www.erudit.de/erudit/events/esit99/12594_p.pdf

62. **Vasantha Kandasamy, W.B., and M. Mary John.** Fuzzy Analysis to Study the Pollution and the Disease Caused by Hazardous Waste From Textile Industries, *Ultra Sci*, **14** (2002) 248-251.

63. **Vasantha Kandasamy, W.B., and R. Praseetha.** New Fuzzy Relation Equations to Estimate the Peak Hours of the Day for Transport Systems, *J. of Bihar Math. Soc.*, **20** (2000) 1-14.

64. **Vasantha Kandasamy, W.B., and Smarandache, F.**, *Analysis of social aspects of migrant labourers living with HIV/AIDS using fuzzy theory and neutrosophic cognitive maps*, Xiquan, Phoenix, 2004.

65. **Vasantha Kandasamy, W.B., and Smarandache, F.**, *Fuzzy Cognitive Maps and Neutrosophic Cognitive Maps*, Xiquan, Phoenix, 2004 USA.

66. **Vasantha Kandasamy, W.B., and Smarandache, F.**, *Fuzzy Relational Equations and Neutrosophic Relational Equations*, Neutrosophic Book Series 3, *HEXIS, Church Rock*, USA, 2004.

67. **Vasantha Kandasamy, W.B., and V. Indra.** Applications of Fuzzy Cognitive Maps to Determine the Maximum Utility of a Route, *J. of Fuzzy Maths*, publ. by the Int. fuzzy Mat. Inst., **8** (2000) 65-77.

68. **Vasantha Kandasamy, W.B., and Yasmin Sultana,** FRM to Analyse the Employee-Employer Relationship Model, *J. Bihar Math. Soc.*, **21** (2001) 25-34.

69. **Vasantha Kandasamy, W.B., and Yasmin Sultana,** Knowledge Processing Using Fuzzy Relational Maps, *Ultra Sci.*, **12** (2000) 242-245.

70. **Vasantha Kandasamy, W.B., N.R. Neelakantan and S. Ramathilagam.** Maximize the Production of Cement Industries by the Maximum Satisfaction of Employees using Fuzzy Matrix, *Ultra Science*, **15** (2003) 45-56.

71. **Vysoký, P.** *Fuzzy Cognitive Maps and their Applications in Medical Diagnostics*. http://www.cbmi.cvut.cz/lab/publikace/30/Vys98_11.doc

72. **Yasmin Sultana,** *Construction of Employee-Employee Relationship Model using Fuzzy Relational Maps*, Masters Dissertation, Guide: Dr. W. B. Vasantha Kandasamy, Department of Mathematics, Indian Institute of Technology, April 2000.

73. **Yuan, Miao and Zhi-Qiang Liu.** On Causal Inference in Fuzzy Cognitive Maps, *IEEE Transactions on Fuzzy Systems*, **81** (2000) 107-119.

74. **Zadeh, L.A.**, A Theory of Approximate Reasoning, *Machine Intelligence*, 9 (1979) 149- 194.

75. **Zhang, W.R., and S. Chen.** A Logical Architecture for Cognitive Maps, *Proceedings of the 2nd IEEE Conference on Neural Networks* (ICNN-88), **1** (1988) 231-238.

76. **Zimmermann, H.J.**, *Fuzzy Set Theory and its Applications*, Kluwer, Boston, 1988.

# INDEX

## D

## E

## F

**G**

**H**

**I**

**J**

**K**

**L**

**M**

**N**

**O**

**P**

**R**

**S**

**T**

**U**

**V**

**W**

# About the Authors

**Dr.W.B.Vasantha Kandasamy** is an Associate Professor in the Department of Mathematics, Indian Institute of Technology Madras, Chennai, where she lives with her husband Dr.K.Kandasamy and daughters Meena and Kama. Her current interests include Smarandache algebraic structures, fuzzy theory, coding/ communication theory. In the past decade she has guided eight Ph.D. scholars in the different fields of non-associative algebras, algebraic coding theory, transportation theory, fuzzy groups, and applications of fuzzy theory of the problems faced in chemical industries and cement industries. Currently, six Ph.D. scholars are working under her guidance. She has to her credit 255 research papers of which 203 are individually authored. Apart from this, she and her students have presented around 294 papers in national and international conferences. She teaches both undergraduate and post-graduate students and has guided over 41 M.Sc. and M.Tech. projects. She has worked in collaboration projects with the Indian Space Research Organization and with the Tamil Nadu State AIDS Control Society. She has authored a Book Series, consisting of ten research books on the topic of Smarandache Algebraic Structures which were published by the American Research Press.

She can be contacted at vasantha@iitm.ac.in
You can visit her work on the web at: http://mat.iitm.ac.in/~wbv

---

**Dr.Florentin Smarandache** is an Associate Professor of Mathematics at the University of New Mexico, Gallup Campus, USA. He published over 60 books and 80 papers and notes in mathematics, philosophy, literature, rebus. In mathematics his research papers are in number theory, non-Euclidean geometry, synthetic geometry, algebraic structures, statistics, and multiple valued logic (fuzzy logic and fuzzy set, neutrosophic logic and neutrosophic set, neutrosophic probability). He contributed with proposed problems and solutions to the Student Mathematical Competitions. His latest interest is in information fusion were he works with Dr.Jean Dezert from ONERA (French National Establishment for Aerospace Research in Paris) in creasing a new theory of plausible and paradoxical reasoning (DSmT).

He can be contacted at smarand@unm.edu